\documentclass[12pt]{article}
\usepackage{xcolor,caption,subcaption,graphicx}

\usepackage{xchicago}

\newcommand {\citeAY} [1] {\citeNP {#1}}%
\newcommand {\citeAPY}[1] {\citeN  {#1}}%
\renewcommand {\showoriginalref}[1]{}
\newcommand\eq[1] {(\ref{#1})}

\newcommand\sect[1] {\ref{sec:#1}}

\newcommand\labsect[1] {\label{sec:#1}}

\newcommand{\bfm}[1]{\mbox{\boldmath ${#1}$}}
\newcommand{\nonum}{\nonumber \\}
\newcommand{\beqa}{\begin{eqnarray}}
\newcommand{\eeqa}[1]{\label{#1}\end{eqnarray}}
\newcommand{\beq}{\begin{equation}}
\newcommand{\eeq}[1]{\label{#1}\end{equation}}

\newcommand{\Ga}{\alpha}
\newcommand{\Gb}{\beta}
\newcommand{\Gd}{\delta}

\newcommand{\Gg}{\gamma}

\newcommand{\Gs}{\sigma}

\newcommand{\Gy}{\psi}

\newcommand{\GY}{\Psi}

\newcommand{\BGv}{\bfm\nu}

\newcommand{\BGF}{\bfm\Phi}
\newcommand{\BGG}{\bfm\Gamma}
\newcommand{\BGL}{\bfm\Lambda}
\newcommand{\BGP}{\bfm\Pi}

\newcommand{\CE}{{\cal E}}

\newcommand{\CH}{{\cal H}}

\newcommand{\CJ}{{\cal J}}
\newcommand{\CK}{{\cal K}}

\newcommand{\CP}{{\cal P}}

\newcommand{\CU}{{\cal U}}
\newcommand{\CV}{{\cal V}}
\newcommand{\CW}{{\cal W}}

\def\Be{{\bf e}}

\def\Bj{{\bf j}}

\def\Bp{{\bf p}}

\def\Bu{{\bf u}}
\def\Bv{{\bf v}}
\def\Bw{{\bf w}}

\def\BA{{\bf A}}
\def\BB{{\bf B}}
\def\BC{{\bf C}}

\def\BE{{\bf E}}
\def\BF{{\bf F}}

\def\BH{{\bf H}}
\def\BI{{\bf I}}
\def\BJ{{\bf J}}
\def\BK{{\bf K}}
\def\BL{{\bf L}}
\def\BM{{\bf M}}

\def\BP{{\bf P}}
\def\BQ{{\bf Q}}

\def\BS{{\bf S}}
\def\BT{{\bf T}}
\def\BU{{\bf U}}
\def\BV{{\bf V}}
\def\BW{{\bf W}}

\def\BY{{\bf Y}}
\def\BZ{{\bf Z}}

\def \ba {\begin{array}}
\def \ea {\end{array}}
\newtheorem {Thm} {Theorem} [section]

\newtheorem {Adef} [Thm] {Definition}

\newtheorem {Arem} [Thm] {Remark}

\newtheorem {Aexa} [Thm] {Example}

\newtheorem {Anot} [Thm] {Notation}

\def \refe #1.{(\ref{#1})}
\def \reff #1.{figure~\ref{#1}}
\def \refs #1.{section~\ref{#1}}
\def \refss #1.{subsection~\ref{#1}}
\def \refD #1.{Definition~\ref{#1}}
\def \refT #1.{Theorem~\ref{#1}}
\def \refL #1.{Lemma~\ref{#1}}
\def \refC #1.{Corollary~\ref{#1}}
\def \refP #1.{Proposition~\ref{#1}}
\def \refR #1.{Remark~\ref{#1}}
\def \refE #1.{Example~\ref{#1}}
\def \refN #1.{Notation~\ref{#1}}

\usepackage{graphics}
\usepackage{amssymb}
\usepackage{amsmath}
\begin{document}
\vspace{-1in}
\title{Superfunctions and the algebra of subspace collections and their association with rational functions of several complex variables}
\author{Graeme W. Milton}
\date{\small{Department of Mathematics, University of Utah, Salt Lake City, UT 84112, USA}}
\maketitle
\begin{abstract}
A natural connection between rational functions of several real or complex variables, and subspace collections
is explored. A new class of function, superfunctions,
are introduced which are the
counterpart to functions at the level of subspace collections. Operations on subspace collections are
found to correspond to various operations on rational functions, such as addition, multiplication and substitution. 
It is established that every rational matrix valued function which is homogeneous of degree 1 can be generated from
an appropriate, but not necessarily unique, subspace collection: the mapping from subspace collections to rational functions is onto, but not one to one. For some applications superfunctions may be more important than functions, as they incorporate more information about the physical problem, yet can be manipulated in much the same way as functions. Previously subspace collections had been introduced when there was an inner product on the vector (or Hilbert) space, and appropriate subspaces were mutually 
orthogonal. In that setting certain normalization and reduction operations on subspace collections led to a continued fraction expansion of the associated function, 
which allowed one to bound the function in terms of a set of weight matrices and normalization matrices that are derived from series expansions. Here we also
initiate the theory of normalization and reduction operations, appropriate when there is no inner product on the space.   

\end{abstract}
\vskip2mm

\section{Introduction}
\setcounter{equation}{0}

This Chapter 7 of the book ''Extending the Theory of Composites to other Areas of Science'', edited by Graeme W. Milton, is concerned
with developing the theory of subspace collections, particularly nonorthogonal subspace collections. 
Subspace collections have a rich algebraic structure, and a close connection with rational functions of several real or complex variables. 
Here we are interested in three types of subspace collections:%
\index{subspace collections}
first, finite dimensional vector spaces $\CH$ that have the decomposition
\beq \CH=\CU\oplus\CE\oplus\CJ=\CP_1\oplus\CP_2\oplus\cdots\oplus\CP_n,
\eeq{I1.2}
which we call a $Z(n)$ subspace collection;%
\index{subspace collections!Z@$Z(n)$}
second finite dimensional vector spaces $\CK$ (over the real or complex numbers) that have the decomposition
\beq \CK=\CE\oplus\CJ=\CV\oplus\CP_1\oplus\CP_2\oplus\cdots\oplus\CP_n,
\eeq{I1.1}
which we call a $Y(n)$ subspace collection,%
\index{subspace collections!Y@$Y(n)$}
where the $\CE$ and $\CJ$ entering \eq{I1.1} are not to be confused with the subspaces $\CE$ and $\CJ$ entering
\eq{I1.2};
and third finite dimensional vector spaces $\CK$ (over the real or complex numbers) 
that have the decomposition
\beq \CK=\CE\oplus\CJ=\CV^I\oplus\CV^O\oplus\CP_1\oplus\CP_2\oplus\cdots\oplus\CP_n,
\eeq{I1.2a}
which we call a superfunction%
\index{superfunctions}
$F^s(n)$. In a superfunction the space $\CV^I$ and the space $\CV^O$ are called the input and output subspaces respectively, and they have the same dimension.  For superfunctions we
require the technical condition that for any choice of vectors $\BE^I, \BJ^I \in\CV^I$ and  $\BE^O, \BJ^O \in\CV^O$ there exist vectors $\BE\in\CE$ and $\BJ\in\CJ$ such that
\beq \BE^I=\BGP^I\BE, \quad \BJ^I=\BGP^I\BJ, \quad \BE^O=\BGP^O\BE, \quad \BJ^O=\BGP^O\BJ. \eeq{I1.2b}

\begin{figure}[ht]
\centering
\includegraphics[width=0.75\textwidth]{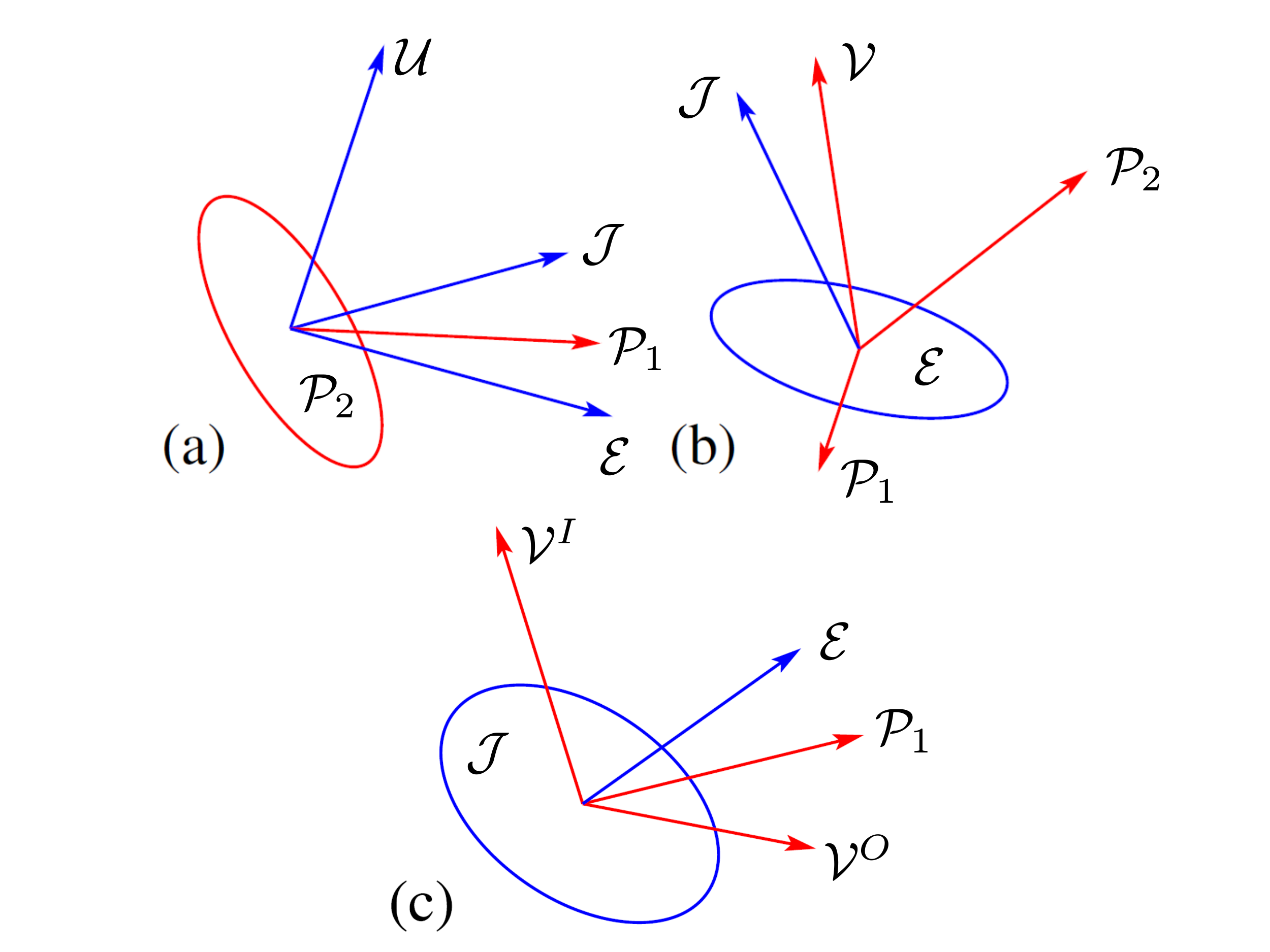}
\caption[Examples of $Z(3)$ and $Y(2)$ subspace collections, and a superfunction $F^s(1)$.]
        {Shown in (a) is an example of a $Z(3)$ subspace collection, in (b) a $Y(2)$ subspace collection, and in (c) a superfunction $F^s(1)$. 
The rays denote one-dimensional subspaces: they should really be drawn
as lines, but for clarity they are drawn as rays and should be extended in the opposite direction as the ray. The circles, which look like ellipses as they are tilted, represent
two-dimensional subspaces.}
\label{IS-1}
\end{figure}

As we will see there is a very close direct connection between a superfunction $F^s(n)$ and a $Y(n)$ subspace collections, 
and also many connections between them and $Z(n)$ subspace collections. All are intertwined 
and that is the beauty of the theory. $Z(3)$ and $Y(2)$ subspace collections, and superfunctions $F^s(1)$ can be visualized in $3$-dimensional space, 
and examples of these are given Figure~\ref{IS-1}.

One reason $Y(n)$ subspace collections, $Z(n)$ subspaces collections, and superfunctions  $F^s(n)$ 
are important is because they arise in many physical problems. For examples in network theory and 
in the theory of the effective moduli of composite materials, see the review in Chapter 2 of this book (\citeAY{Milton:2016:ETC})
and \citeAPY{Milton:2002:TOC}. 
There are also many other physical problems where subspace collections arise as is apparent in Chapters 1,3,8,9,
12, 13, and 14
of this book (\citeAY{Milton:2016:ETC}). In physics applications the subspaces are usually orthogonal with respect to some inner product on the space $\CH$ or $\CK$ but as this chapter shows the theory of them can be developed without reference to an inner product. This generalization is important to make contact between general rational functions of complex variables, thus extending the notion of a function: hence the name superfunction. The generalization is also important for applications, such as speeding up numerical methods for calculating
the fields that solve the problem: we will see an example of this in the next chapter.

\begin{figure}[ht]
\centering
\includegraphics[width=0.75\textwidth]{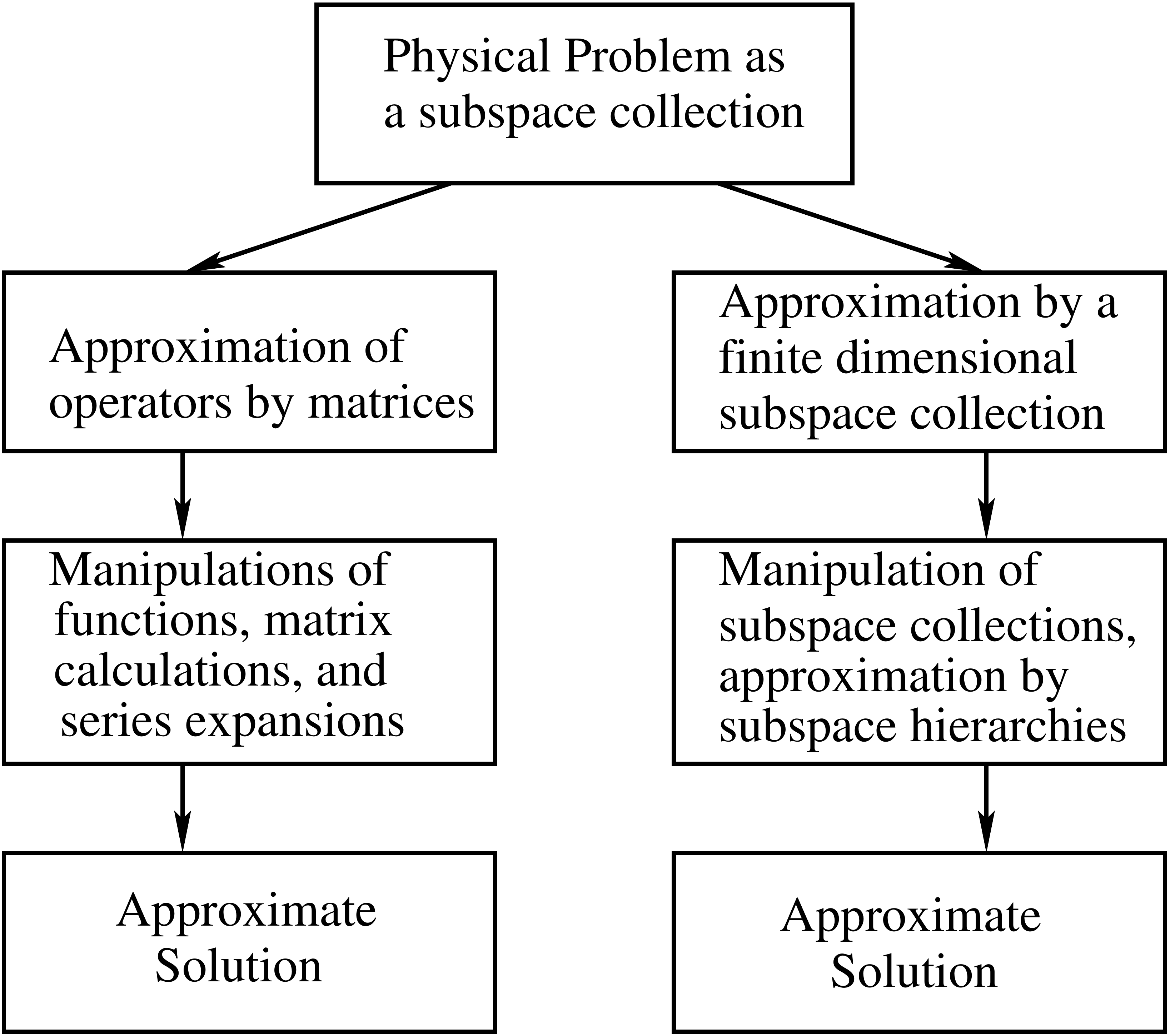}
\caption{Two routes to solving a physical problem formulated in terms of subspace collections. It is suggested that the route on the right may result in a better
approximation as more information is kept.}\label{IS0}
\end{figure}

It may very well be the case that superfunctions become more important than functions in some applications, as suggested
by the flow chart of Figure~\ref{IS0}. The reason is that
when one extracts the function from a superfunction, which we will see how to do shortly, one generally
loses information that is contained in a superfunction. For example, in the context of physical problems 
where there is an inner product on the space this information may came in the form of a series expansion for the fields up to a given order, 
and from this series expansion one can extract the ``weight matrices''%
\index{weight matrices}
and ``normalization matrices'',
\index{normalization matrices}
introduced by \citeAPY{Milton:1985:TCC} and Milton (\citeyearNP{Milton:1987:MCEa}, \citeyearNP{Milton:1987:MCEb}). 
These matrices basically encode the information about the ``angles'' between the various subspaces (when there is an inner product). One can then develop a continued fraction expansion%
\index{continued fraction expansion}
 for the function associated with the superfunction, with the normalization factors and weight matrices that enter it at each level having the property that they are positive semidefinite, with the weight matrices summing to one. Truncating the continued fraction gives approximations to the function, that are similar in some respects the diagonal Pad\'e approximants,%
\index{Pad\'e approximants}
and in fact give bounds%
\index{bounds!truncated continued fractions}
on the function if the truncation is done appropriately. The information contained in the 
weight matrices and normalization matrices, cannot in general be recovered (at least when $n\geq 4$) from the series expansion of the associated function.
(Although one can potentially determine these matrices from the series expansion of the functions associated with coupled field problems,%
\index{coupled field problems}
as shown in Chapter 9 of this book (\citeAY{Milton:2016:ETC})).      
This theory was established by \citeAPY{Milton:1985:TCC} 
and Milton (\citeyearNP{Milton:1987:MCEa}, \citeyearNP{Milton:1987:MCEb}, \citeyearNP{Milton:1991:FER}).
(see also Chapters 19, 20 and 29 in \citeAPY{Milton:2002:TOC}) for the case of 
$Z(n)$ subspace collections, for any integer $n\geq 1$. In this paper develop the basic theory
of subspace collections in the case where there is no inner product on the vector space $\CH$ or $\CK$. We also make the
first steps towards generating continued fraction expansions in the case where there is 
no inner product on the vector space $\CH$ or $\CK$.  

Let us first suppose $\CV$ and $\CU$ are one-dimensional. We will see that there are generally homogeneous (of degree $1$)
rational functions $Y(z_1, z_2,\ldots,z_n)$ and $Z(z_1,z_2,\ldots,z_n)$ (over the real or complex numbers) of 
degree 1 that are associated respectively with these $Y(n)$ and $Z(n)$ subspace
collections, where  $Z(z_1,z_2,\ldots,z_n)$ satisfies the additional constraint that $Z(1,1,\ldots,1)=1$. Conversely, we will see that given any
rational functions $Y(z_1, z_2,\ldots,z_n)$ and $Z(z_1,z_2,\ldots,z_n)$ with these properties, then there exists at least one subspace collection realizing these functions
as its associated function. There are also operations on these subspace collections that correspond to operations on the associated function, such as substitution.%
\index{subspace collections!substitution}

For superfunctions the simplest case is when the input and output spaces%
\index{superfunctions!input and output spaces}
$\CV^I$ and $\CV^O$ are 
one-dimensional. Then with a specific basis for $\CV^I$ and $\CV^O$ 
the corresponding function $\BF(z_1,z_2,\ldots,z_n)$ is 2 by 2 matrix valued with the elements $F_{11}(z_1,z_2,\ldots,z_n)$ and $F_{22}(z_1,z_2,\ldots,z_n)$ 
being homogeneous of degree zero, the element $F_{12}(z_1,z_2,\ldots,z_n)$ being homogeneous of degree minus 1, and $F_{21}(z_1,z_2,\ldots,z_n)$ being homogeneous of degree 1.
There are operations on superfunctions that correspond to addition, multiplication and forming an inverse (and hence division) of the associated functions.
So superfunctions form an algebra.%
\index{superfunctions!algebra}
Also one can do substitutions at the level of subspace collections. Actually the 
operation of addition of superfunctions are naturally done with the associated $Y$-problem, although one could equally do them with the associated inverse $Y$-problem (where the spaces $\CE$ and $\CJ$
are interchanged). Thus there is an inherent ambiguity of how one wants to define addition of superfunctions. The definitions of
addition, multiplication. and substitution of subspace collections may seem a little complicated and abstract, yet they are the exact
counterpart of similar operations one may do on multiterminal electrical networks, and they do produce the corresponding action on
the associated functions. (In fact it was thinking about electrical circuits%
\index{electrical circuit}
which guided the construction of these operations in a more general setting). 

When $\CV$ and $\CU$ have dimension greater than 1, then $Y(z_1, z_2,\ldots,z_n)$ and $Z(z_1,z_2,\ldots,z_n)$ get replaced by linear operator valued functions
$\BY(z_1, z_2,\ldots,z_n)$ and $\BZ(z_1,z_2,\ldots,z_n)$ which map $\CV$ to $\CV$ and $\CU$ to $\CU$ respectively. Similarly, the function $\BF(z_1,z_2,\ldots,z_n)$ should
really be thought of as a linear operator mapping $\CV^I$ to $\CV^O$

The original motivation%
\index{subspace collections!motivation}
for studying subspace collections, and their associated functions, arose from the study of the effective conductivity tensor $\BZ$ of periodic composite materials. For a 
composite with $n$ isotropic phases, with scalar conductivities $z_1, z_2, \ldots, z_n$, the effective conductivity tensor%
\index{effective conductivity tensor}
was found to be a homogeneous (of degree $1$) analytic function 
$\BZ(z_1,z_2,\ldots,z_n)$ of the component conductivities with positive definite imaginary part when the component conductivities have positive imaginary part 
[\citeAY{Bergman:1978:DCC}; Milton \citeyearNP{Milton:1979:TST}, \citeyearNP{Milton:1981:BCP}, \citeAY{Golden:1983:BEP}] (see also Chapter 18 of \citeAPY{Milton:2002:TOC}). 
It was also recognized (Milton \citeyearNP{Milton:1987:MCEa}, \citeyearNP{Milton:1990:CSP}) that the problem of determining the effective conductivity function could be formulated in terms of three mutually orthogonal spaces in the Hilbert space $\CH$ of square integrable functions: namely the space $\CU$ of constant fields, the space $\CE$ of periodic square integrable electric fields (having zero curl), and the space $\CJ$ of square integrable
current fields (having zero divergence), and if the composite had $n$ isotropic phases, with conductivities $z_1, z_2, \ldots, z_n$, then it was also natural to decompose $\CH$ into
the direct sum of $n$ mutually orthogonal subspaces $\CP_1, \CP_2,\ldots, \CP_n$ where $\CP_i$ consists of those square integrable fields which are nonzero only within component $i$.
This formulation, in terms of a $Z(n)$ subspace collection,%
\index{subspace collections!Z@$Z(n)$}
evolved out of earlier Hilbert space formulations%
\index{Hilbert space formulations}
of the problem (\citeAY{Fokin:1982:IMT}; \citeAY{Kohler:1982:BEC}; 
\citeAY{Papanicolaou:1982:BVP}; \citeAY{Golden:1983:BEP}; \citeAY{Kantor:1984:IRB}; \citeAY{DellAntonio:1986:ATO}) and can easily be extended to the elastic, thermoelastic, piezoelectric, and poroelastic equations of multiphase and polycrystalline materials (see, for example, Chapter 12
in \citeAPY{Milton:2002:TOC}). The formulation has proved to be particularly important in the theory of exact relations%
\index{composites!exact relations}
of composite materials
(\citeAY{Grabovsky:1998:EREa}; \citeAY{Grabovsky:1998:EREb}:  \citeAY{Grabovsky:1998:ERC}; \citeAY{Grabovsky:2000:ERE}; \citeAY{Grabovsky:2004:AGC})
(see also Chapter 17 in \citeAPY{Milton:2002:TOC}) where one seeks microstructure independent relations%
\index{microstructure independent relations}
 satisfied by effective tensors. For two-dimensional polycrystals%
\index{polycrystals}
a complete 
correspondence was established between subspace collections and a representative class of multiple rank laminate polycrystal geometries (\citeAY{Clark:1994:MEC}), thus showing that the
subspace collection of any two-dimensional polycrystal, with any configuration of crystal grains, could be approximated arbitrarily closely by the subspace collection of one of these 
multiple rank laminate polycrystal%
\index{polycrystal!multiple rank laminate}
geometries.

Curiously the connection between $Z(n)$ subspace collections and the effective conductivity allowed the effective conductivity function $\BZ(z_1,z_2,\ldots,z_n)$ to be expanded as a new type of continued fraction,%
\index{continued fraction expansion}
involving matrices of increasing dimension as one proceeds down the continued fraction 
when $n>2$ (Milton \citeyearNP{Milton:1987:MCEa}, \citeyearNP{Milton:1987:MCEb}, \citeyearNP{Milton:1991:FER}; see also Chapters 19, 20 and 29 in \citeAY{Milton:2002:TOC}). The coefficients in the weight%
\index{weight matrices}
and normalization matrices%
\index{normalization matrices}
entering the continued fraction can be expressed in terms of inner products%
\index{inner product}
between fields that enter the series expansion%
\index{series expansions}
of the
solution field in a nearly homogeneous medium (with all the conductivities $z_1, z_2, \ldots, z_n$ being close to one another).  One application of the continued fraction expansion
has been to obtain bounds%
\index{bounds!three-phase composites}
on the diagonal elements of the complex effective conductivity tensor of a three phase conducting composite, with
complex conductivities $z_1$, $z_2$ and $z_3$, that were tighter than bounds obtained by any other method (see figure 4 in \citeAY{Milton:1987:MCEb}). 
This procedure essentially extended to multivariate functions the procedure, using successive fractional linear transformations, 
that was used to obtain bounds (\citeAY{Baker:1969:BEB}) on the values in the complex plane that Stieltjes functions!
\index{Stieltjes functions}
can take 
when a finite number of Taylor series coefficients are known (see also \citeAY{Golden:1983:BEP}; \citeAY{Bergman:1986:EDC}) 
where essentially the same transformation is used to derive bounds
on the complex dielectric constant of two component media using series expansion coefficients, as noted in the appendix in \citeAY{Milton:1986:MPC}, 
and see \citeAY{Milton:1981:BTO}, where the same set of bounds is derived using a different procedure, namely the method of variation of poles and zeros.)%
\index{variation of poles and zeros method}

In the case $n=2$ the continued fraction reduces to a usual continued fraction expansion, like those continued fractions 
associated with Pad\'e approximants%
\index{Pad\'e approximants}
(see Chapter 4 of Part I of \citeAY{Baker:1981:PAB}). $Y(n)$ subspace collections
enter, for example, if one eliminates from the Hilbert space the constant fields and then reformulates the conductivity equations in terms of the remaining fields: 
the driving fields are then fields which are constant in each phase, but have zero average value (see Chapter 19 in \citeAY{Milton:2002:TOC} and references therein). 
The interrelationship between $Z(n)$ subspace collections and $Y(n)$ subspace collections is what gives rise to these novel continued fractions. 

\begin{figure}[ht]
\centering
\includegraphics[width=0.75\textwidth]{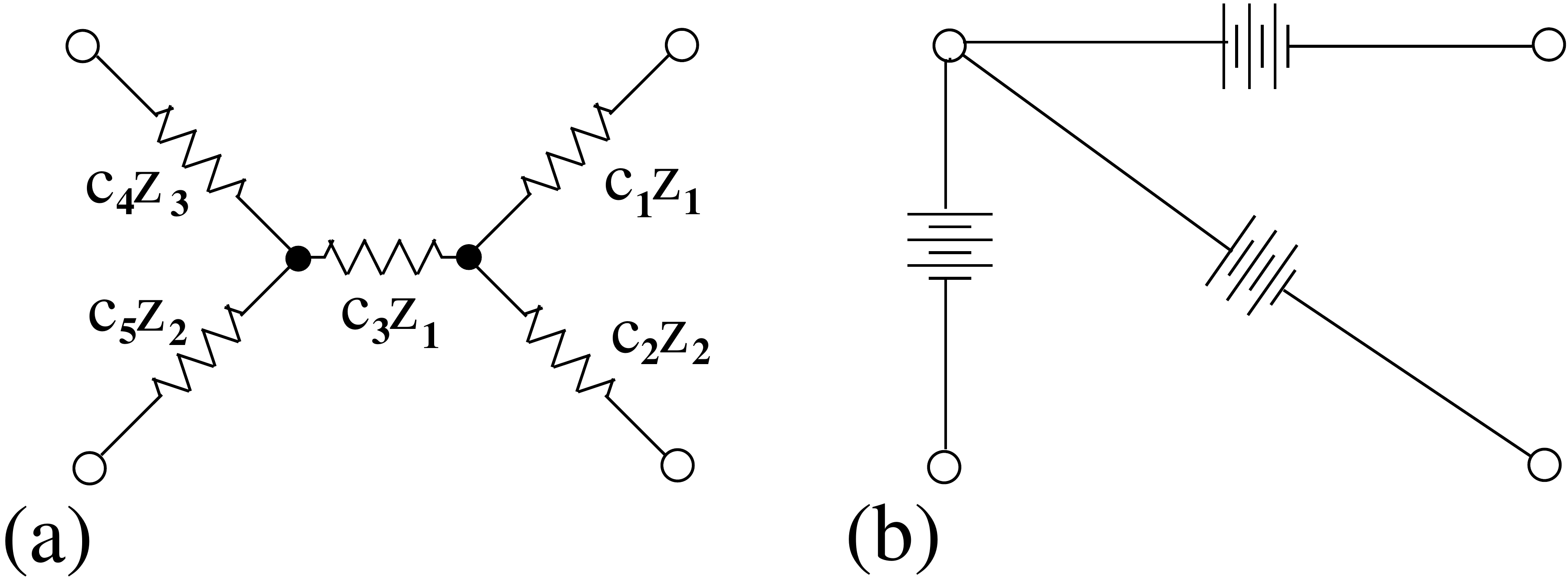}
\caption[Electrical circuits as representative of $Z(n)$ subspace collections, and of $Y(n)$ subspace collections when the batteries are added.]
        {Shown in (a) is a 4 terminal electrical network, which is representative of a $Z(3)$ subspace collection. Here the
$c_i$ are real positive scaling constants: the conductance of each element is $c_jz_k$ where $z_k$ is real or complex (when 
$z$ is complex we should refer to $c_jz_k$ as an  admittance rather than as a conductance). Complex values of $z$ are appropriate
when the applied potentials vary sinusoidally with time, and some of the impedence elements are capacitors or inductors.
Figure (b) shows the batteries on the back side of the circuit board, representing the space $\CV$, which combined with
the resistors on the front side is representative of a $Y(3)$ subspace collection.
The $Y$-function $\BY(z_1,z_2,z_3)$ 
gives the current going through the three batteries, in response to the voltages across them.}\label{IS1}
\end{figure}

Finite dimensional $Z(n)$ and $Y(n)$ subspace collections also arise naturally in the study of the effective resistance of electrical circuits%
\index{electrical circuit}
constructed from $n$ types of
resistors having conductances $z_1$, $z_2$, $\ldots z_n$ (see Chapter 20 in \citeAY{Milton:2002:TOC}). This is not surprising as periodic resistor networks%
\index{composite!approximated by resistor network}
can be seen as discrete approximations to conducting composite materials (see, for example, \citeAY{Milton:1981:BCP} and Figure~8.5(a) in this book
\citeAY{Milton:2016:ETC}). Figure~\ref{IS1} illustrates a discrete network of impedances, and gives an indication of the physical meaning of the $Z(n)$ and $Y(n)$ subspace collections in this context.

In this figure, the vector space $\CH$ is 6-dimensional, and is the direct sum of the two-dimensional space
$\CP_1$ consisting of fields that are nonzero only along the resistors $c_1z_1$ and $c_3z_1$;
the two-dimensional space $\CP_2$ consisting of fields that are nonzero only along the resistors $c_2z_2$ and $c_5z_2$;
and the one-dimensional space $\CP_3$ consisting of fields that are nonzero only along the resistor $c_4z_3$.
The response of the network, when one terminal is grounded (with zero voltage) 
is a $3\times 3$ matrix.  When it acts on the vector, having as elements the voltages at the three remaining terminals,
it gives the three currents flowing through these terminals. The $3\times 3$ matrix valued 
function $\BZ(z_1,z_2,z_3)$ gives the matrix valued response relative to the response when $z_1=z_2=z_3=1$.
Now, let us imagine all the resistors, or impedances, in (a) are on one side of the circuit board, with the terminals being conducting
posts that penetrate the board. On the other side of the board these posts are connected to a tree-like graph%
\index{tree-like graph}
of batteries (or alternating current sources if the fields vary sinusoidally in time) shown in (b). The three fields in these
batteries constitute the space $\CV$. The $Y(3)$ subspace collection 
contains fields on both sides of the board, in $\CK=\CH\oplus\CV$. The associated $3\times 3$ 
matrix valued $Y$-function $\BY(z_1,z_2,z_3)$ gives the current going through the three batteries, in response to the voltages
across them. Note that  $\BY(z_1,z_2,z_3)$ is not diagonal: a voltage across one battery, sends current through the other two batteries,
even when they have zero voltage across them.

\begin{figure}[ht]
\centering
\includegraphics[width=0.90\textwidth]{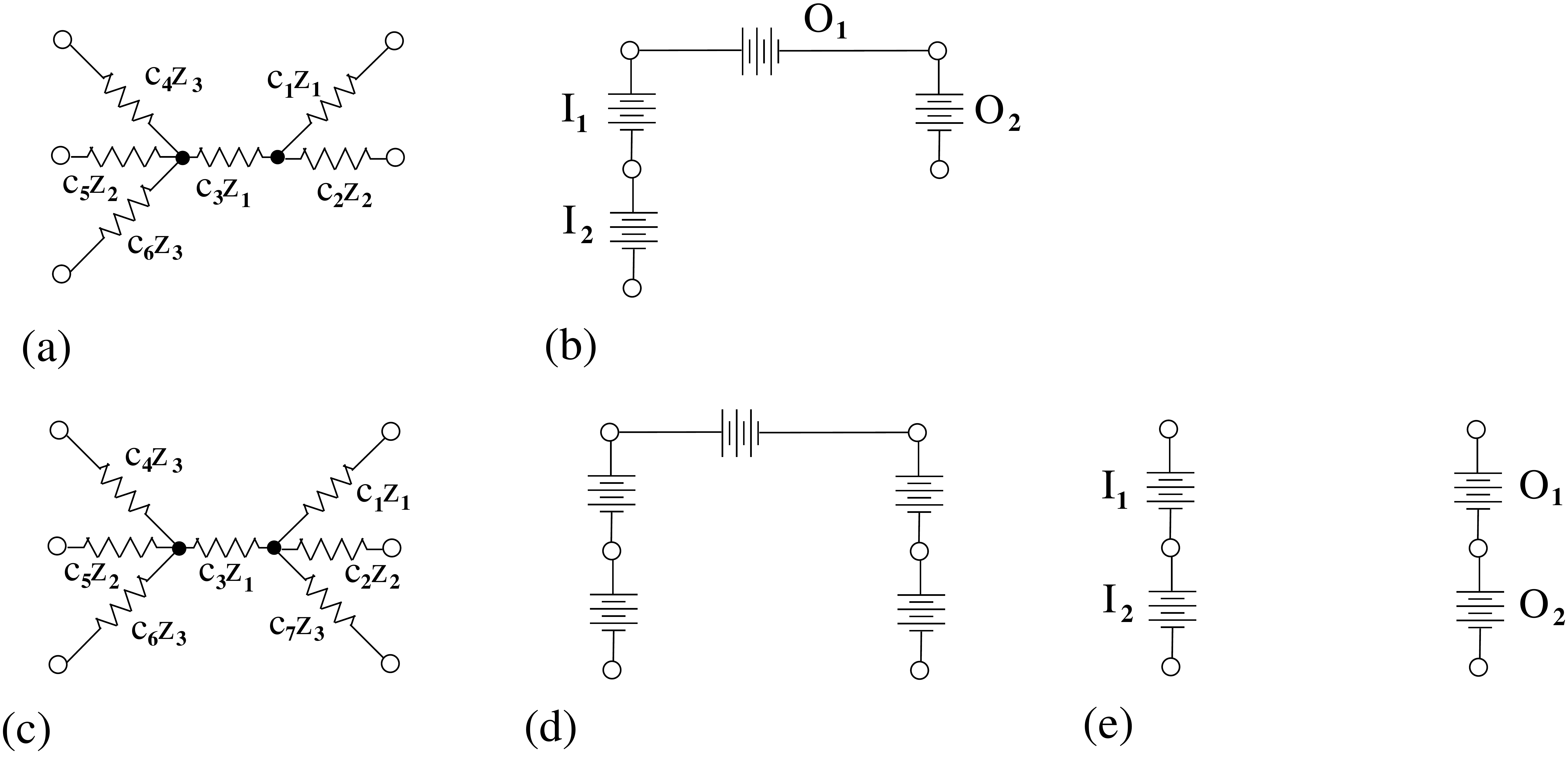}
\caption[Multiport electrical circuits as representative of superfunctions.]{Shown in (a) is a 5 terminal electrical network, 
which is representative of a $Z(3)$ subspace collection. Here the
$c_i$ are real positive scaling constants: the admittance of each element is $c_jz_k$ where $z_k$ is real or complex. 
Figure (b) shows the batteries on the back side of the circuit board, representing the space $\CV$, which is divided into the
input space $\CV^I$, consisting of those vectors in $\CK$ that are nonzero only in the batteries $I_1$ and $I_2$ and
the output space $\CV^O$, consisting of those vectors in $\CK$ that are nonzero only in the batteries/resistors $0_1$ and $0_2$.
Figure (c) shows a 6 terminal electrical network, and the naturally associated subspace $\CV$ represented by the batteries in Figure (d). To convert this
to a problem where the dimension of $\CV$ is even we remove the battery at the top, and accordingly reduce the dimension of both $\CV$ and $\CJ$ by one.  
Figure (e) shows the input space $\CV^I$, consisting of those vectors that are nonzero only in the batteries $I_1$ and $I_2$ and
the output space $\CV^O$, consisting of those vectors in $\CK$ that are nonzero only in the batteries/resistors $0_1$ and $0_2$.}\label{IS1.00}
\end{figure}

Superfunctions are a natural generalization of multiport electrical circuits with input ports and output ports, as illustrated in Figure~\ref{IS1.00}.
The function $\BF$ gives the currents and potential drops across the output batteries/resistors that are generated in response to currents and potential 
drops across the input batteries. Note that the networks associated with superfunctions automatically satisfy the ``port condition'' that the net flow of current
from the input terminals to the output terminals is zero.

In this chapter we show that the connection between finite dimensional $Z(n)$ and $Y(n)$ subspace collections and homogeneous (degree $1$) 
operator valued rational functions $\BY(z_1, z_2,\ldots,z_n)$ and $\BZ(z_1,z_2,\ldots,z_n)$ persists even when the subspaces in each 
decomposition are not necessarily mutually orthogonal, and indeed even
in the absence of an inner product (on the space $\CH$ or $\CK$). The results 
developed in (Milton, \citeyearNP{Milton:1987:MCEa}, \citeyearNP{Milton:1987:MCEb}, \citeyearNP{Milton:1991:FER} 
and in Chapters 19, 20 and 29 of Milton, \citeyearNP{Milton:2002:TOC}) are extended to the case where there is no inner product.
Accordingly some steps in the analysis, and some assumptions, need to be revised. In this more general setting we can generate, from an 
appropriate $Z(n)$ subspace collection, any desired scalar valued rational function $Z(z_1, z_2,\ldots,z_n)$ satisfying the homogeneity property $Z(1,1,\ldots,1)=1$.

It is to be emphasized that subspace collections, with the associated rules for addition, multiplication and substitution, are algebraic objects in their own right: there is no need
to think of the associated analytic functions (that are in general operator valued), except that the correspondence makes it easier to think about subspace collections. 
The resistor network examples of $Y(n)$ subspace collections made it possible for me to see how the operations of addition, multiplication and substitution of subspace collections
should be defined in the general case.

My belief is that the geometrical structure of subspace collections (and in particular superfunctions) will be reflected in the algebraic geometrical structure
of their associated rational functions. If this is the case, understanding the topological features of subspace collections%
\index{subspace collections!topological features}
might shed light on the geometrical features of algebraic varieties. While this paper does not directly address this issue, it sheds the 
first light on the relation between finite dimensional subspace collections and rational functions of several complex variables, 
in the case where the subspaces are not mutually orthogonal, and it introduces superfunctions. The functions derived from superfunctions are well studied and have 
widespread applications in signal processing,%
\index{signal processing}
control theory,%
\index{control theory}
network synthesis and design,%
\index{network synthesis and design}
and in optics, acoustics and elastodynamics (usually in layered media), where 
they are called a variety of names including transfer matrices,%
\index{transfer matrices}
transmission matrices,%
\index{transmission matrices}
transfer functions,
\index{transfer functions}
system functions,%
\index{system functions}
and network functions.%
\index{network functions}
In these contexts it is 
the function that is studied, but people do not think of the superfunction. 
I thank Aaron Welters and Mihai Putinar for drawing my attention to the connection between transfer functions and 
response functions (such the effective conductivity tensor of composites).   

We remark that for $Z(3)$ orthogonal subspace collections, with $\CU$ being one-dimensional, it is still an open and intriguing question as to whether
there could be a one-to-one correspondence between them (assuming they are pruned as described in Section \sect{Iprune}
and modulo trivial equivalences between subspace collections) and scalar functions $Z(z_1,z_2,z_3)$ satisfying the homogeneity,%
\index{homogeneity property}
Herglotz%
\index{Herglotz property}
and normalization properties.%
\index{normalization property}
 The $Z$-problem described the next section provides a nonlinear map from
the $Z(3)$ orthogonal subspace collection to an associated scalar function $Z(z_1,z_2,z_3)$ satisfying the homogeneity,
Herglotz and normalization properties, but the question is whether one can uniquely recover the pruned subspace collection, modulo
trivial equivalences, given only the function $Z(z_1,z_2,z_3)$? The intriguing counting argument given in Section 29.2 of \citeAPY{Milton:2002:TOC} suggests
the possibility of a one-to-one correspondence. There is a similar counting argument for 
nonorthogonal subspace collections given in Section \sect{Icount},
but in this case we will see in an explicit example that a one-to-one correspondence does not hold. 

\section{Subspace collections and their associated functions}
\labsect{Isub}
\setcounter{equation}{0}
Let $\CK$ be a vector space which has a decomposition into two different direct sums of subspaces
\beq \CK=\CE\oplus\CJ=\CV\oplus\CH, \eeq{I2.0}
where $\CH$ itself is a direct sum of $n$ subspaces
\beq \CH=\CP_1\oplus\CP_2\oplus\cdots\oplus\CP_n.
\eeq{I2.1}
Any vector $\BK\in\CK$ has a unique decomposition into component vectors,
\beq \BK=\BE+\BJ=\Bv+\BH,\quad \BH=\BP_1+\BP_2+\cdots+\BP_n, \eeq{I2.2}
each in the associated subspaces:
\beq \BE\in\CE,\quad\BJ\in\CJ,\quad\Bv\in\CV,\quad\BH\in\CH,\quad\BP_i\in\CP_i~{\rm for}~i=1,2,\ldots,n.
\eeq{I2.3}
This decomposition serves to define projection operators $\BGG_1$ and $\BGG_2$%
\index{projection operators!gammaa@$\BGG_1$ and $\BGG_2$}
onto $\CE$ and $\CJ$, projection operators $\BGP_1$ and $\BGP_2$%
\index{projection operators!P@$\BGP_1$ and $\BGP_2$}
onto $\CV$ and $\CH$, and projection 
operators $\BGL_i$ onto the subspaces $\CP_i$. By definition we have
\beq \BE=\BGG_1\BK,\quad\BJ=\BGG_2\BK,\quad\Bv=\BGP_1\BK,\quad\BH=\BGP_2\BK,\quad \BP_i=\BGL_i\BK. \eeq{I2.4}
Associated with this subspace collection is an linear operator valued function 
$\BY(z_1,z_2,\ldots,z_n)$ acting on the space $\CV$, which is a homogeneous function of degree 1 of the
$n$ complex variables $z_1,z_2,\ldots,z_n$. To obtain the function we take each field $\BE_1\in\CV$ and
look for vectors $\BJ$ and $\BE$ that solve the equations
\beq \BE\in\CE,~~~~\BJ\in\CJ,~~~~
\BJ_2=\BL\BE_2,~~~~{\rm where}~\BJ_2=\BGP_2\BJ,~~~\BE_2=\BGP_2\BE,
\eeq{I2.5}
with $\BE_1=\BGP_1\BE$, where
\beq \BL=\sum_{i=1}^nz_i\BGL_i. \eeq{I2.5a}
We call this problem the $Y$-problem.%
\index{Y@$Y$-problem}
The associated operator $\BY$, by definition, governs the linear relation
\beq \BJ_1=-\BY\BE_1,~~~~{\rm where}~~
\BJ_1=\BGP_1\BJ.
\eeq{I2.6}
A necessary condition for $\BJ_1$ to be uniquely defined given $\BE_1$ is that
\beq \CV\cap\CJ=0, \eeq{I2.6a}
since if $\BJ$ and $\BE$ solve \eq{I2.5} so too will $\BJ+\Bv$ and $\BE$, for any $\Bv\in\CV\cap\CJ$.
The inverse $Y$-problem%
\index{Y@$Y$-problem!inverse}
is to solve \eq{I2.5} for each field $\BJ_1=\BGP_1\BJ\in\CV$.
A necessary condition for $\BE_1$ to be uniquely defined given $\BJ_1$ is that
\beq \CV\cap\CE=0. \eeq{I2.6b}
If $\Bv_1,\Bv_2,\ldots,\Bv_m$ is a
basis of $\CV$, then the operator $\BY$ can be represented by a matrix, the $Y$-matrix, also denoted by $\BY$
with elements $Y_{ik}$ such that
\beq \BY\Bv_k=\sum_{i=1}^mY_{ik}\Bv_i. \eeq{I2.7}

If $m$ is even and $\CV$ has the decomposition
\beq \CV=\CV^I\oplus\CV^O, \eeq{I2.7a}
where $\CV^I$ and $\CV^O$ have the same dimension ($m/2$) then we have a superfunction%
\index{superfunctions}
$F^s$. The superfunction is the collection of subspaces and
there is a function $\BF$ associated with it. The fields $\BE_1$ and $\BJ_1$ have the unique decomposition  
\beq \BE_1=\BE^I+\BE^O,\quad \BJ_1=\BJ^I+\BJ^O, \eeq{I2.7aa}
with
\beq \BE^I,\BJ^I\in\CV^I,\quad \BE^O,\BJ^O\in\CV^O, \eeq{I2.7b}
where the superscripts $I$ and $O$ refer to input and output respectively. We write
\beq \BE^I=\BGP^I\BE_1,\quad \BE^O=\BGP^O\BE_1,\quad \BJ^I=\BGP^I\BJ_1,\quad  \BJ^O=\BGP^O\BJ_1, 
\eeq{I2.7ba}
which defines the projections $\BGP^I$ and $\BGP^O$%
\index{projection operators!piy@$\BGP^I$ and $\BGP^O$}
onto the input and output spaces.%
\index{superfunctions!input and output spaces}
Now the relation \eq{I2.6} can be written as
\beq \begin{pmatrix} \BJ^I  \\ \BJ^O \end{pmatrix}=
\begin{pmatrix} \BY^{II}  &  \BY^{IO} \\
               \BY^{OI}  &  \BY^{OO}\end{pmatrix}
\begin{pmatrix} \BE^I  \\ \BE^O \end{pmatrix},
\eeq{I2.7c}
and manipulated into the form
\beq  \begin{pmatrix} \BE^O  \\ \BJ^O \end{pmatrix}=\BF\begin{pmatrix} \BE^I  \\ \BJ^I \end{pmatrix},
\eeq{I2.7d}
which defines the linear operator valued function
\beq \BF=\begin{pmatrix} \BF^{EE}  &  \BF^{EJ} \\
               \BF^{JE}  &  \BF^{JJ}\end{pmatrix}
         =\begin{pmatrix} -(\BY^{IO})^{-1}\BY^{II}  & \quad -(\BY^{IO})^{-1} \\
                           [\BY^{OO}(\BY^{IO})^{-1}\BY^{II}-\BY^{OI}]  & \quad \BY^{OO}(\BY^{IO})^{-1}\end{pmatrix},
\eeq{I2.7e}
which, provided the operator $\BY^{IO}$ is nonsingular, is the function associated with the superfunction. This relation can be inverted to yield $\BY$ in terms of $\BF$,
\beq \BY=\begin{pmatrix} (\BF^{EJ})^{-1}\BF^{EE}  & \quad -(\BF^{EJ})^{-1} \\
                            [\BF^{JJ}(\BF^{EJ})^{-1}\BF^{EE}-\BF^{JE}]  & \quad -\BF^{JJ}(\BF^{EJ})^{-1}\end{pmatrix},
\eeq{I2.7f}
provided the operator $\BF^{EJ}$ can be inverted. The superfunction problem%
\index{superfunctions!problem}
is for given input fields $\BE^I$ and $\BJ^I$ to find
fields $\BE$ and $\BJ$ that solve the $Y$-problem \eq{I2.5} and \eq{I2.5a}, with $\BGP^I\BE=\BE^I$ and $\BGP^I\BJ=\BJ^I$. It may happen that the superfunction problem has a solution when the
$Y$-problem does not (this happens when $\BF^{EJ}$ is singular), and conversely the $Y$-problem may have a solution when the superfunction problem does not 
(this happens when $\BY^{IO}$ is singular). 

Another association between subspace collections and functions comes if a vector space $\CH$ has the
decomposition
\beq \CH=\CU\oplus\CE\oplus\CJ=\CP_1\oplus\CP_2\oplus\cdots\oplus\CP_n, \eeq{I2.8}
where $\CE$ and $\CJ$ are not to be confused with the spaces in \eq{I2.0}. Any vector $\BH\in\CH$ 
has a unique decomposition into component vectors,
\beq \BH=\Bu+\BE+\BJ=\BP_1+\BP_2+\cdots+\BP_n, \eeq{I2.9}
each in the associated subspaces:
\beq \Bu\in\CU,\quad\BE\in\CE,\quad\BJ\in\CJ,\quad\BP_i\in\CP_i~{\rm for}~i=1,2,\ldots,n.
\eeq{I2.10}
This decomposition serves to define projection operators $\BGG_0$, $\BGG_1$ and $\BGG_2$%
\index{projection operators!gamma@$\BGG_0$, $\BGG_1$, and $\BGG_2$} 
onto
$\CU$, $\CE$ and $\CJ$, and projection operators $\BGL_i$ onto the subspaces $\CP_i$.
Associated with this subspace collection is an linear operator valued function 
$\BZ(z_1,z_2,\ldots,z_n)$ acting on the space $\CU$, which is a homogeneous function of degree 1 of the
$n$ complex variables $z_1,z_2,\ldots,z_n$. To obtain the function we take each vector $\Be\in\CU$ and 
look for vectors $\Bj$, $\BJ$ and $\BE$ that solve the equations
\beq \Bj\in\CU,~~~~\BE\in\CE,~~~~\BJ\in\CJ,~~~~
\Bj+\BJ=\BL(\Be+\BE),~~~~{\rm where}~\BL=\sum_{i=1}^nz_i\BGL_i. \eeq{I2.11}
We call this problem the $Z$-problem.%
\index{Z@$Z$-problem}
The associated operator $\BZ$, by definition, governs the linear relation
\beq \Bj=\BZ\Be. \eeq{I2.12}
If $\Bu_1,\Bu_2,\ldots,\Bu_m$ is a
basis of $\CU$, then the operator $\BZ$ can be represented by a matrix, also denoted by $\BZ$
with elements $Z_{ik}$ such that
\beq \BZ\Bu_k=\sum_{i=1}^mZ_{ik}\Bu_i. \eeq{I2.13}
When $z_1=z_2=\cdots=z_n=1$ \eq{I2.11} has the trivial solution
\beq \Bj=\Be,\quad \BJ=\BE=0, \eeq{I2.14}
and so we deduce that
\beq \BZ(1,1,\ldots,1)=\BI. \eeq{I2.15}
The inverse $Z$-problem%
\index{Z@$Z$-problem!inverse}
is to solve the equations \eq{I2.11} for each given vector $\Bj\in\CU$.

\section{Some simple examples}
\setcounter{equation}{0}

Consider a $Y(n)$ subspace collection%
\index{subspace collections!Y@$Y(n)$}
\beq \CK=\CE\oplus\CJ=\CV\oplus\CP_1\oplus\CP_2\oplus\cdots\oplus\CP_n,
\eeq{I3.1}
where $\CE, \CV, \CP_1, \CP_2, \ldots \CP_n$ are all one-dimensional, and $\CJ$ 
is $n$-dimensional. Choose, as our basis for $\CK$, $n+1$ vectors $\Bp_0\in\CV$,
and $\Bp_i\in\CP_i$, $i=1,2,\ldots n$. Vectors $\BE\in\CE$ and $\BJ\in\CJ$ can be
expanded in this basis:
\beq \BE=\sum_{i=0}^nE_i\Bp_i,\quad \BJ=\sum_{i=0}^nJ_i\Bp_i. \eeq{I3.2}
The relation $\BGP_2\BJ=\BL\BGP_2\BE$ implies
\beq J_i=z_iE_1. \eeq{I3.3}

Let us suppose that $E_0=1$. Then $E_1$ and $E_2$ are determined by the orientation of the
one-dimensional subspace $\CE$ with respect to the subspaces $\CV, \CP_1, \CP_2, \ldots \CP_n$.
Also since $\CJ$ has  codimension 1, there exist constants $W_0, W_1, \ldots W_n$, determined
by the orientation of the $n$-dimensional subspace $\CJ$ with respect to the subspaces 
$\CV, \CP_1, \CP_2, \ldots \CP_n$ such that
\beq \sum_{i=0}^nW_iJ_i=0. \eeq{I3.4}
Let us suppose that $W_0=1$. Then we have
\beq J_0=-\sum_{i=1}^nW_iJ_i=-\sum_{i=1}^nW_iE_iz_i, \eeq{I3.5}
which since $E_0=1$ implies $J_0=-YE_0$, with
\beq Y=\sum_{i=1}^n\Ga_iz_i,\quad{\rm where}~\Ga_i=W_iE_i. \eeq{I3.6}
As the $E_i$ and $W_i$ are arbitrary constants, we see that $Y$ can be
any linear combination of the $z_i$. In particular, with $W_1E_1=1$ and $W_iE_i=0$ when $i\ne 1$
we obtain
\beq Y=z_1. \eeq{I3.6a}

As a second example consider a $Y(1)$ subspace collection
\beq \CK=\CE\oplus\CJ=\CV\oplus\CP_1,
\eeq{I3.6b}
where all the spaces $\CE$, $\CJ$, $\CV$, and $\CP_1$ are all two-dimensional. Choose as our basis for $\CK$ two vectors $\Bp_1$ and $\Bp_2$ in
$\CV$ and two vectors $\Bp_3$ and $\Bp_4$ in $\CP_1$. Then since $\CE$ is two-dimensional, there generically exist constants 
$e_{13}, e_{14}, e_{23}$ and $e_{24}$ such that 
\beq \Bp_1+e_{13}\Bp_3+e_{14}\Bp_4 \in \CE, \quad \Bp_2+e_{23}\Bp_3+e_{24}\Bp_4 \in \CE. \eeq{I3.6c}
Also since $\CJ$ is two-dimensional, there generically exist constants 
$j_{31}, j_{32}, j_{41}$ and $j_{42}$ such that
\beq \Bp_3+j_{31}\Bp_1+j_{32}\Bp_2 \in \CJ, \quad \Bp_4+j_{41}\Bp_1+e_{42}\Bp_2 \in \CJ. \eeq{I3.6d}
So the $Y$-problem is solved with vectors
\beqa \BE &= &\Bp_1+e_{13}\Bp_3+e_{14}\Bp_4, \nonum
      \BE_1& = & \Bp_1,\quad \BE_2=e_{13}\Bp_3+e_{14}\Bp_4, \nonum
      \BJ_2& = & z_1(e_{13}\Bp_3+e_{14}\Bp_4), \nonum
      \BJ &= &z_1[e_{13}(\Bp_3+j_{31}\Bp_1+j_{32}\Bp_2)+e_{14}(\Bp_4+j_{41}\Bp_1+e_{42}\Bp_2)], \nonum
      \BJ_1& = &z_1[(e_{13}j_{31}+ e_{14}j_{41})\Bp_1+(e_{13}j_{32}+ e_{42}j_{42})\Bp_2,
\eeqa{I3.6e}
and is also solved with vectors
\beqa \BE& = &\Bp_2+e_{23}\Bp_3+e_{24}\Bp_4, \nonum
      \BE_1& = & \Bp_2,\quad \BE_2=e_{23}\Bp_3+e_{24}\Bp_4, \nonum
      \BJ_2 & = & z_1(e_{23}\Bp_3+e_{24}\Bp_4), \nonum
      \BJ & = & z_1[e_{23}(\Bp_3+j_{31}\Bp_1+j_{32}\Bp_2)+e_{24}(\Bp_4+j_{41}\Bp_1+e_{42}\Bp_2)], \nonum
      \BJ_1& = & z_1[(e_{23}j_{31}+ e_{24}j_{41})\Bp_1+(e_{23}j_{32}+ e_{24}j_{42})\Bp_2.
\eeqa{I3.6f}
From these equations in follows that $\BY(z_1)$ in this basis is the 2 by 2 matrix
\beq \BY(z_1) = z_1\BA,\quad {\rm with}~~ \BA=\begin{pmatrix} a_{11}  &  a_{12} \\ a_{21}  &  a_{22} \end{pmatrix},
\eeq{I3.6g}
where
\beqa a_{11}=e_{13}j_{31}+ e_{14}j_{41}, \quad  a_{12}=e_{13}j_{32}+ e_{42}j_{42}, \nonum
      a_{21}=e_{23}j_{31}+ e_{24}j_{41}, \quad  a_{22}=e_{23}j_{32}+ e_{24}j_{42}.
\eeqa{I3.6h}
As the coefficients $e_{13}, e_{14}, e_{23}, e_{24}, j_{31}, j_{32}, j_{41}$ and $j_{42}$ can be any complex numbers we desire it follows that 
we can realize any desired complex matrix $\BA$. By taking $\CV^I$ to be the one-dimensional space spanned by $\Bp_1$ and
taking $\CV^O$ to be the one-dimensional space spanned by $\Bp_2$ we obtain a superfunction $Y^S$ where the associated function 
takes the form
\beq  \BF(z_1)=\begin{pmatrix} b_{11}  &  b_{12}/z_1 \\ b_{21}z_1  &  b_{22} \end{pmatrix},
\eeq{I3.6i}
in which the parameters $b_{11}, b_{12}, b_{21}$ and $b_{22}$ can be any complex numbers we choose.

As a third example consider a $Z(2)$ subspace collection%
\index{subspace collections!ZZ@$Z(2)$}
\beq \CH=\CU\oplus\CE\oplus\CJ=\CP_1\oplus\CP_2,
\eeq{I3.7}
where the subspaces  $\CU, \CE, \CJ$ and $\CP_2$ are all one-dimensional, while $\CP_1$ is 
two-dimensional. Choose, as our basis for $\CH$, $3$ vectors $\BU_0\in\CU$, $\BE_0\in\CE$ and
$\BJ_0\in\CJ$, and take a vector $\BP$ as a basis for $\CP_2$. The coefficients $P_U$, $P_E$
and $P_J$ in the expansion
\beq \BP=P_U\BU_0+P_E\BE_0+P_J\BJ_0 \eeq{I3.8}
determine the orientation of $\CP_2$ with respect to the subspaces $\CU$, $\CE$ and $\CJ$. 
In the basis $\BU_0$, $\BE_0$, and $\BJ_0$ the equations
\beq \Be+\BE=\BQ+\Ga\BP,\quad \Bj+\BJ=z_1\BQ+z_2\Ga\BP,\eeq{I3.9}
with
\beq
\Be,\Bj\in\CU,\quad \BE\in\CE,\quad\BJ\in\CJ,\quad\BQ\in\CP_1, \eeq{I3.10}
take the form
\beqa \begin{pmatrix} e \\ E \\ 0 \end{pmatrix}
& = &\begin{pmatrix} Q_U \\ Q_E \\ Q_J\end{pmatrix}
+ \Ga\begin{pmatrix} P_U \\ P_E \\ P_J \end{pmatrix}, \nonum
 \begin{pmatrix} j \\ 0 \\ J \end{pmatrix}
& = & z_1\begin{pmatrix} Q_U \\ Q_E \\ Q_J\end{pmatrix}
+z_2\Ga\begin{pmatrix} P_U \\ P_E \\ P_J \end{pmatrix},
\eeqa{I3.11}
and since $\BQ\in\CP_1$ there exist constants $W_U$, $W_E$ and $W_J$, which determine
the orientation of $\CP_1$ with respect to $\CU$, $\CE$ and $\CJ$, such that
\beq W_UQ_U+W_EQ_E+W_JQ_J=0. \eeq{I3.12}
Hence we obtain the equations
\beqa W_Ue+W_EE & = & \Ga(W_UP_U+W_EP_E+W_JP_J)\equiv\Ga\BW\cdot\BP, \nonum
      0 & = & z_1(E-\Ga P_E)+z_2\Ga P_E, \nonum
      j & = & z_1(e-\Ga P_U)+z_2\Ga P_U.
\eeqa{I3.13}
Eliminating $E$ and $\Ga$ from these equations gives $j=Z e$, with
\beq Z=z_1+\frac{(z_2-z_1)W_UP_U}{\BW\cdot\BP+W_EP_E(z_2-z_1)/z_1}.
\eeq{I3.14}
In particular if the subspaces are oriented so that
\beq \BW\cdot\BP=W_EP_E=-W_UP_U, \eeq{I3.15}
then \eq{I3.14} gives
\beq Z=z_1^2/z_2, \eeq{I3.16}
which with $z_2=1$ produces the function $z_1^2$ and with $z_1=1$ 
produces the function $1/z_2$. Also, with $W_EP_E=0$ we obtain
\beq Z=z_1+\frac{(z_2-z_1)W_UP_U}{\BW\cdot\BP}, \eeq{I3.17}
which is a ``weighted average'' of $z_1$ and $z_2$, $Z=w_1z_1+w_2z_2$ with ``weights''%
\index{weights}
$w_1$ and $w_2$ that sum to 1
but which are not necessarily positive, nor even real.

\section{Formulas for the associated functions}
\labsect{Isubformulas}
\setcounter{equation}{0}
Following Section 12.8 of \citeAPY{Milton:2002:TOC} a formula for the effective tensor $\BZ$ results by
applying the operator $\BGG_0+\BGG_2$ (which projects on the space
$\CU\oplus\CJ$) to both sides of
the constitutive law $\Be+\BE=\BL^{-1}(\Bj+\BJ)$. Solving the resulting equation,
\beq \Be=(\BGG_0+\BGG_2)\BL^{-1}(\BGG_0+\BGG_2)(\Bj+\BJ), \eeq{I4.8}
for $\Bj+\BJ$ gives
\beq \Bj+\BJ=[(\BGG_0+\BGG_2)\BL^{-1}(\BGG_0+\BGG_2)]^{-1}\Be,
\eeq{I4.9}
where the last inverse is to be taken on the subspace $\CU\oplus\CJ$.
By applying $\BGG_0$ to both sides of this equation we see that
\beq \BZ=\BGG_0[(\BGG_0+\BGG_2)\BL^{-1}(\BGG_0+\BGG_2)]^{-1}\BGG_0,
\eeq{I4.10}
which is the result given in (12.59) of \citeAPY{Milton:2002:TOC}.

Another formula for $\BZ$ follows from noting that for any arbitrary constant $z_0\ne 0$,
\beq [z_0\BI-\BGG_1(\BL-z_0\BI)](\Be+\BE)=z_0\Be+z_0\BE-\BGG_1\BJ-z_0\BGG_1\BE=z_0\Be.
\eeq{I4.11}
Solving this for $\Be+\BE$ gives
\beq \Be+\BE=z_0[z_0\BI-\BGG_1(\BL-z_0\BI)]^{-1}\Be, \eeq{I4.12}
and applying $\BGG_0\BL$ to both sides yields
\beq \Bj=z_0\BGG_0\BL[z_0\BI-\BGG_1(\BL-z_0\BI)]^{-1}\Be. \eeq{I4.13}
Thus we have a formula for the $\BZ$ operator,%
\index{Zo@$\BZ$-operator formula}
\beq \BZ=z_0\BGG_0\BL[z_0\BI-\BGG_1(\BL-z_0\BI)]^{-1}\BGG_0
=z_0\BGG_0+z_0\BGG_0(\BL-z_0\BI)[z_0\BI-\BGG_1(\BL-z_0\BI)]^{-1}\BGG_0,
\eeq{I4.14}
where we have used the identity
\beq \BGG_0=z_0\BGG_0[z_0\BI-\BGG_1(\BL-z_0\BI)]^{-1}\BGG_0, \eeq{I4.15}
obtained by applying $\BGG_0$ to both sides of \eq{I4.12}. This formula \eq{I4.14} is a special
case of the formula (12.60) given in \citeAPY{Milton:2002:TOC}, and is well known in different
contexts (\citeAY{Kroner:1977:BEE}).

To obtain a formula for $\BY$ notice that \eq{I2.5}
and \eq{I2.6} imply that
\beq 0=\BGG_2\BE'=\BGG_2\BE_1+\BGG_2\BE_2=\BGG_2\BE_1+
\BGG_2\BL^{-1}\BGP_2\BGG_2\BJ', \eeq{I4.3}
where the inverse of $\BL$ is to be taken on the subspace $\CH$.
Solving for $\BJ'$ gives
\beq \BJ'=-(\BGG_2\BL^{-1}\BGP_2\BGG_2)^{-1}\BGG_2\BE_1, \eeq{I4.4}
where the inverse is to be taken on the subspace $\CJ$. Then by applying $\BGP_1$
to both sides of this equation and equating $\BGP_1\BJ'=\BJ_1$ with $-\BY\BE_1$
we obtain the desired formula%
\index{Yo@$\BY$-operator formula}
\beq \BY=\BGP_1\BGG_2(\BGG_2\BL^{-1}\BGP_2\BGG_2)^{-1}\BGG_2\BGP_1,
\eeq{I4.5}
for $\BY$, as given in formula (19.29) of \citeAPY{Milton:2002:TOC}.

Another formula for $\BY$ is obtained by taking an arbitrary constant $z_0\ne 0$,
and defining
\beq \BP'=\BJ'-z_0\BE'. \eeq{I4.5a}
Applying $\BGG_1$ to both sides of \eq{I4.5a} gives
\beq \BGG_1\BP'=-z_0\BE'=-z_0(\BE_1+\BE_2), \eeq{I4.5b}
and applying $\BGP_2$ to both sides of \eq{I4.5b} gives
\beq \BGP_2\BP'=\BJ_2-z_0\BE_2=(\BL-z_0\BI)\BE_2. \eeq{I4.5c}
Combining these results we see that $\BP'$ satisfies
\beq [\BGG_1+z_0(\BL-z_0\BI)^{-1}\BGP_1]\BP'=-z_0\BE_1. \eeq{I4.5d}
Assuming that the operator $[\BGG_1+z_0(\BL-z_0\BI)^{-1}\BGP_1]$ is nonsingular
this gives 
\beq \BP'=-z_0[\BGG_1+z_0(\BL-z_0\BI)^{-1}\BGP_1]^{-1}\BE_1. \eeq{I4.5e}
Applying $\BGP_1=\BI-\BGP_2$ to both sides yields
\beq \BJ_1-z_0\BE_1=-(\BY+z_0\BI)\BE_1=-z_0\BGG_1[\BGG_1+z_0(\BL-z_0\BI)^{-1}\BGP_1]^{-1}\BE_1
\eeq{I4.5f}
As this holds for all $\BE_1\in\CV$ we obtain the formula%
\index{Yo@$\BY$-operator formula}
\beq \BY=-z_0\BGP_1+z_0\BGG_1[\BGG_1+z_0(\BL-z_0\BI)^{-1}\BGP_1]^{-1}\BGP_1 \eeq{I4.5g}
which is a special case of the formula (19.37) obtained in Section 19.5 of \citeAPY{Milton:2002:TOC}.

\section{Multiplying superfunctions}
\begin{figure}[!ht]
\centering
\includegraphics[width=0.75\textwidth]{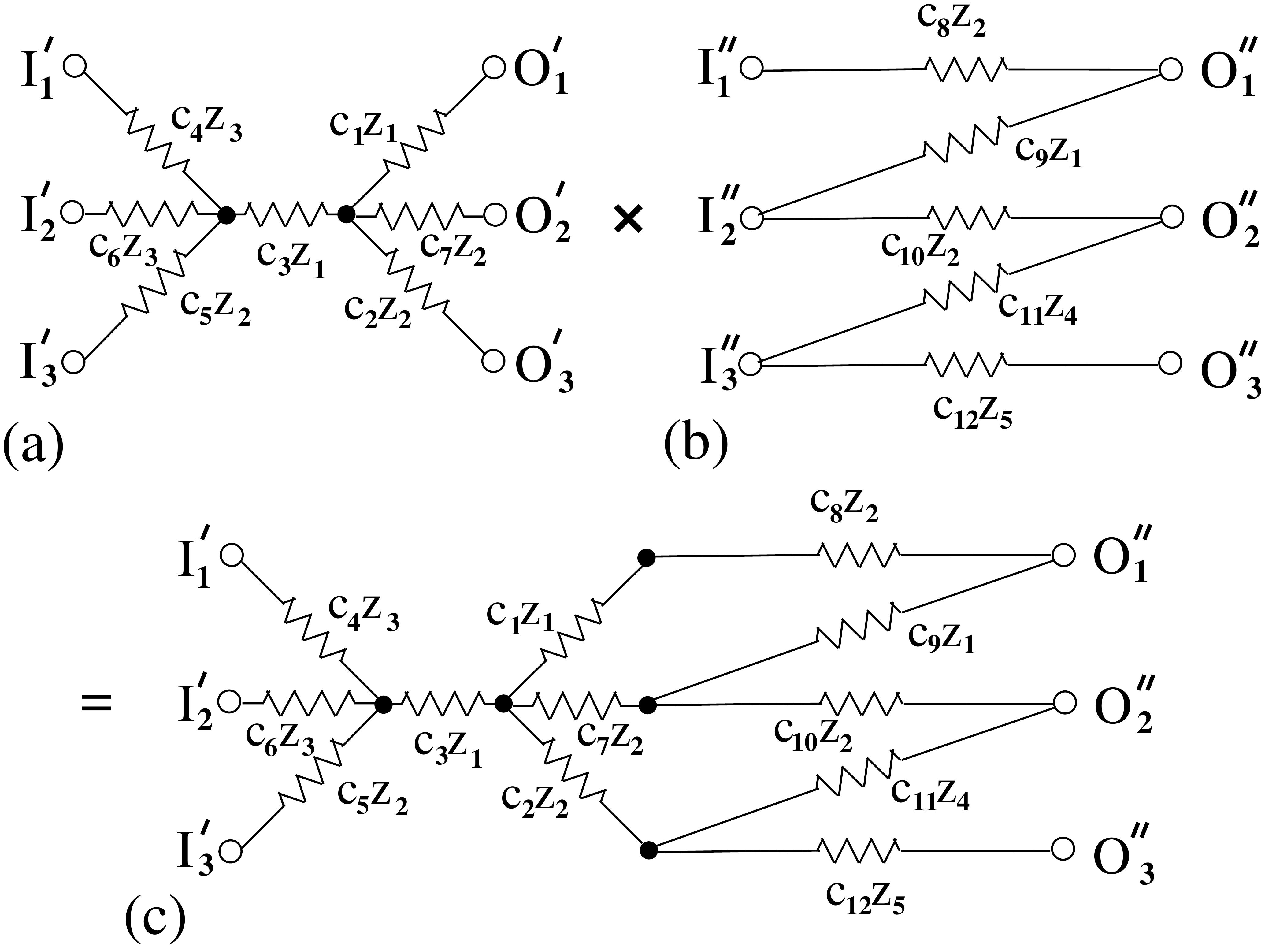}
\caption[Multiplying superfunctions is like hooking networks, with an equal number of input and output terminals, together in series.]{Multiplying superfunctions is like hooking networks, with an equal number of input and output terminals, together in series. Shown in (a) and (b) are 6 terminal electrical networks, each (along with their respective tree-like 
battery configurations on the opposite side of the circuit board that are not shown here)
represent a superfunction as the terminals 
have been divided into input terminals ($I_1'$, $I_2'$, and $I_3'$ for the circuit (a), and $I''_1$, $I''_2$, and $I''_3$ 
for the circuit (b)) and output terminals ($O_1'$, $O_2'$, and $O_3'$ for the circuit (a), and $O''_1$, $O''_2$, and $O''_3$ 
for the circuit (b)). The product superfunction is the  6 terminal electrical network (along with its tree-like 
battery configurations on the opposite side of the circuit board ) shown in (c). Note there is some flexibility
in how one multiplies superfunctions: instead of connecting the terminals $O'_i$ with $I''_i$ for $i=1,2,3$, one could
for example, connect $O_1'$,$O_2'$, and $O_3'$ with any permutation of $I_1''$, $I_2''$ and $I_3''$. This is why, when taking a
product, one needs to specify the maps ($\BM^E$ and $\BM^J$) one is using between the output space of one superfunction, and the input
space of the second superfunction by which one is multiplying it.}\label{IS2}
\end{figure}

\setcounter{equation}{0}
Multiplying superfunctions%
\index{superfunctions!multiplication}
is similar the way electrical circuits, each with $2m$ terminal can be combined. An example
is shown in Figure~\ref{IS2}.

Suppose we have two superfunctions, $(F^s)'$ and $(F^s)''$:
\beqa \CK' & = & \CE'\oplus\CJ'=(\CV^I)'\oplus(\CV^O)'\oplus\CH'\quad{\rm with}~~
\CH'=\CP_1'\oplus\CP_2'\oplus\cdots\oplus\CP_j', \nonum
 \CK'' & = & \CE''\oplus\CJ''=(\CV^I)''\oplus(\CV^O)''\oplus\CH''\quad{\rm with}~~
\CH''=\CP_1''\oplus\CP_2''\oplus\cdots\oplus\CP_k'',
\eeqa{I0.1}
where the spaces $(\CV^I)',(\CV^O)',(\CV^I)'',(\CV^O)''$ all have the same dimension $m$. To take their product
one needs to first find two nonsingular linear operators $\BM^E$ and $\BM^J$ which map $(\CV^O)'$ to $(\CV^I)''$.
The resulting product superfunction
\beq F^s=(F^s)'\times_{\BM}(F^s)'', \eeq{I0.2}
is the subspace collection
\beq \CK=\CE\oplus\CJ=(\CV^I)'\oplus(\CV^O)''\oplus\CH, \eeq{I0.2a}
where
\beq \CH=\CP_1'\oplus\CP_2'\oplus\cdots\oplus\CP_j'\oplus\CP_1''\oplus\CP_2''\oplus\cdots\oplus\CP_k''',
\eeq{I0.3}
and the operator $\BL$ acting on $\CH$ is
\beq \BL=\sum_{i=1}^j z_i'\BGL_i'+\sum_{\ell=1}^k z_\ell'' \BGL_\ell'', \eeq{I0.3a}
in which $\BGL_i'$ and $\BGL_\ell''$ are the projections onto $\CP_i'$ and $\CP_\ell''$.
A vector $\BE$ is in $\CE$ if and only if we can find vectors  
\beqa \BE' & = & (\BE^I)' +(\BE^O)'+\BE_2'\in\CE', \nonum
      \BE'' & = & (\BE^I)'' +(\BE^O)''+\BE_2''\in\CE'',
\eeqa{I0.4}
such that 
\beq
(\BE^I)''=\BM^E(\BE^O)',\quad \BE=(\BE^I)'+(\BE^O)''+\BE_2'+\BE_2'',
\eeq{I0.5}
with
\beq
 (\BE^I)'\in(\CV^I)',\quad (\BE^O)'\in(\CV^O)',\quad\BE_2'\in\CH',\quad (\BE^I)''\in(\CV^I)'',\quad (\BE^O)''\in(\CV^O)'',\quad\BE_2''\in\CH''.
\eeq{I0.6}
A vector  $\BJ$ is in $\CJ$ if and only if we can find vectors 
\beqa \BJ' & = & (\BJ^I)' +(\BJ^O)'+\BJ_2'\in\CJ', \nonum
      \BJ'' & = & (\BJ^I)'' +(\BJ^O)''+\BJ_2''\in\CJ'',
\eeqa{I0.7}
such that 
\beq
(\BJ^I)''=\BM^J(\BJ^O)',\quad \BJ=(\BJ^I)'+(\BJ^O)''+\BJ_2'+\BJ_2'',
\eeq{I0.8}
with   
\beq 
(\BJ^I)'\in(\CV^I)',\quad (\BJ^O)'\in(\CV^O)',\quad\BJ_2'\in\CH',\quad (\BJ^I)''\in(\CV^I)'',\quad (\BJ^O)''\in(\CV^O)'',\quad\BJ_2''\in\CH''.
\eeq{I0.9}

To ensure that the two spaces $\CE$ and $\CJ$ are independent we need to make the technical assumption that $\BM^E$ and $\BM^J$ are chosen so 
that the operator $\BA$ mapping $(\CV^O)'$ to $(\CV^I)''$, defined by
\beq \BA=\BM^E(\BGP^O)'\BGG_1'-(\BGP^I)''\BGG_1''[\BM^E(\BGP^O)'\BGG_1'+\BM^J(\BGP^O)'\BGG_2'], \eeq{I0.9a}
is nonsingular (i.e. the null--space of the operator contains only the zero vector).  Our aim is to show that if $\BA$ is nonsingular and 
\beq \BE=(\BE^I)'+(\BE^O)''+\BE_2'+\BE_2''=\BJ=(\BJ^I)'+(\BJ^O)''+\BJ_2'+\BJ_2'',\quad {\rm with}\,\,\BE\in\CE,\quad\BJ\in\CJ
\eeq{I0.10}
then $\BE=\BJ=0$. First note that by resolving \eq{I0.10} into components in the spaces  $(\CV^I)'$, $(\CV^I)''$, $\CH'$, and $\CH''$ we obtain
\beq  (\BE^I)'= (\BJ^I)',\quad (\BE^O)''= (\BJ^O)'',\quad \BE_2'=\BJ_2',\quad \BE_2''=\BJ_2''. \eeq{I0.11}
Also since $\BE\in\CE$ and $\BJ\in\CJ$ there exist vectors  $(\BE^O)', (\BJ^O)' \in (\CV^O)'$ and $(\BE^I)'', (\BJ^I)''\in (\CV^I)''$ such that 
\eq{I0.4} and \eq{I0.7} hold. Since $\CE'\cap\CJ'=\{0\}$ and  $\CE''\cap\CJ''=\{0\}$ it follows that 
\beq \BP\equiv (\BE^O)'-(\BJ^O)'=\BE'-\BJ'\ne 0\quad{\rm or}\quad\BE'=\BJ'=0, \eeq{I0.11a1}
and
\beq \BQ\equiv (\BE^I)''-(\BJ^I)''=\BE''-\BJ''\ne 0\quad{\rm or}\quad\BE''=\BJ''=0.  \eeq{I0.11a2}
Now we have
\beqa (\BGP^O)'\BGG_1'\BP & = & (\BGP^O)'\BE'= (\BE^O)',\quad (\BGP^O)'\BGG_2'\BP=-(\BGP^O)'\BJ'=-(\BJ^O)' \nonum
 (\BGP^I)''\BGG_1''\BQ & = & (\BGP^I)''\BE''= (\BE^I)'',\quad (\BGP^I)''\BGG_2''\BQ=-(\BGP^I)''\BJ''=-(\BJ^I)''.
\eeqa{I0.11b}
Since $(\BE^I)''=\BM^E(\BE^O)'$ and $(\BJ^I)''=\BM^J(\BJ^O)'$ we get from the first pair of equations in \eq{I0.11b} the result that
\beq (\BE^I)''=\BM^E(\BGP^O)'\BGG_1'\BP,\quad (\BJ^I)''=-\BM^J (\BGP^O)'\BGG_2'\BP, \eeq{I0.11c}
which implies
\beq \BQ=[\BM^E(\BGP^O)'\BGG_1'+\BM^J (\BGP^O)'\BGG_2']\BP. \eeq{I0.11d}
Substituting this back in the second pair of equations in \eq{I0.11b}, and using \eq{I0.11c}, gives
\beqa  (\BGP^I)''\BGG_1''[\BM^E(\BGP^O)'\BGG_1'+\BM^J (\BGP^O)'\BGG_2']\BP=\BM^E(\BGP^O)'\BGG_1'\BP \nonum
 (\BGP^I)''\BGG_2''[\BM^E(\BGP^O)'\BGG_1'+\BM^J (\BGP^O)'\BGG_2']\BP=\BM^J(\BGP^O)'\BGG_2'\BP.
\eeqa{I0.11e}
These two equations are not independent since by adding them we obtain
\beq [\BM^E(\BGP^O)'\BGG_1'+\BM^J (\BGP^O)'\BGG_2']\BP=\BM^E(\BGP^O)'\BGG_1'\BP+\BM^J(\BGP^O)'\BGG_2'\BP \eeq{I0.11f}
which is obviously satisfied. Also the first equation in \eq{I0.11e} says $\BP$ is in the null space of $\BA$, which by our assumption
implies $\BP=0$. Then \eq{I0.11d} implies $\BQ=0$ and this rules out the first possibilities in \eq{I0.11a1} and \eq{I0.11a2},
implying $\BE'=\BJ'=0$ and $\BE''=\BJ''=0$. We conclude that  $\BE=\BJ=0$.

To check that the space $\CE\oplus\CJ$ spans $(\CV^I)'\oplus(\CV^O)''\oplus\CH$ , we just need to count dimensions. The dimension of the space
on the right is $2m$+dim($\CH$). The dimension of $\CE$ according to \eq{I0.4} is dim($\CE'$)+dim($\CE''$) less $m$ because of the $m$ constraints
$(\BE^I)''=\BM^E(\BE^O)'$. Similarly the dimension of $\CJ$ is dim($\CJ'$)+dim($\CJ''$)-$m$. Adding these up, we get the dimension of
$\CE\oplus\CJ$ is dim$\CK'$+dim$\CK''$-$2m$=$2m$+dim($\CH'$)+dim($\CH'$)=$2m$+dim($\CH$).

Let $\BF'$ and $\BF'$ be the functions associated with the superfunctions $(F^s)'$ and $(F^s)''$.
Given operators
\beq \BL'=\sum_{i=1}^jz_i'\BGL_i',\quad \BL''=\sum_{i=1}^kz_i''\BGL_i'', \eeq{I0.15}
where $\BGL_i'$ projects onto $\CP_i'$ and $\BGL_i''$ projects onto $\CP_i''$, and given
input fields $(\BE^I)'$ and $(\BJ^I)'$ we can calculate
\beqa  \begin{pmatrix} (\BE^O)'  \\ (\BJ^O)' \end{pmatrix}& = & \BF'\begin{pmatrix} (\BE^I)'  \\ (\BJ^I)' \end{pmatrix}, \nonum
  (\BE^I)' & = & \BM^E(\BE^O)'', \quad  (\BJ^I)'=\BM^J(\BJ^O)'',     \nonum
       \begin{pmatrix} (\BE^O)''  \\ (\BJ^O)'' \end{pmatrix}& = & \BF''\begin{pmatrix} (\BE^I)''  \\ (\BJ^I)'' \end{pmatrix}.
\eeqa{I0.16}
From the knowledge of $(\BE^O)'$ and $(\BE^I)'$, and of $(\BE^O)''$ and $(\BE^I)''$,  
we can calculate the fields $\BE'$, $\BE''$, $\BJ'$, and $\BJ''$ of the form \eq{I0.4} and \eq{I0.7} 
solving the $Y'$ problem and the $Y''$ problem:
\beqa 
\BE'\in\CE',\quad \BJ'\in\CJ',\quad\BJ'_1=\BL'\BE'_1,\nonum
\BE''\in\CE'',\quad \BJ''\in\CJ'',\quad\BJ''_1=\BL''\BE''_1.
\eeqa{I0.17}
Then the fields $\BE$ and $\BJ$ given by \eq{I0.5} and \eq{I0.8} solve the $Y$ problem in the space $\CK$, and the function associated to 
the superfunction $F^s$ is given by the product rule%
\index{superfunctions!product rule}
\beq \BF=\BF'\begin{pmatrix} \BM^E & 0  \\ 0 & \BM^J \end{pmatrix}\BF''.
\eeq{I0.18}
Let us choose a basis $(\Bv_1^I)'',(\Bv_2^I)'',\ldots (\Bv_{m}^I)''$ for $(\CV^I)''$, choose a basis
$(\Bv_1^O)'',(\Bv_2^O)'',\ldots (\Bv_{m}^O)''$ for $(\CV^O)''$, take 
$\BM^E(\Bv_1^O)'',\BM^E(\Bv_2^O)'',\ldots \BM^E(\Bv_{m}^O)''$ as our basis for $(\CV^I)'$, and 
choose a basis $(\Bv_1^O)',(\Bv_2^O)',\ldots (\Bv_{m}^O)'$ for $(\CV^O)'$. Then the operator $\BM^E$ is represented as the
identity matrix in the basis. Let us also choose the operator $\BM^J$ so it is represented by {\it minus} the identity matrix in this basis. Then
in this basis the relation \eq{I0.18} takes the form
\beq \BF=\BF'\begin{pmatrix} \BI & 0  \\ 0 & -\BI \end{pmatrix}\BF''.
\eeq{I0.19}
Note that we could have avoided this slightly awkward multiplication rule if we had replaced the definition \eq{I2.7d} of the associated function by
\beq  \begin{pmatrix} \BE^O  \\ -\BJ^O \end{pmatrix}=\BF\begin{pmatrix} \BE^I  \\ \BJ^I \end{pmatrix}.
\eeq{I0.19a}
Then the multiplication rule (with this choice of $\BM^E$ and $\BM^J$) would have become simply $\BF=\BF'\BF''$. We chose not to do this in the interest of preserving 
the symmetric roles of the spaces $\CE$ and $\CJ$ in the definition of the function associated with the superfunction.

In passing, let us suppose there is an inner product on the vector spaces $\CK'$ and $\CK"$, and that the sets of spaces
$\{\CE',\CJ'\}$, $\{(\CV^I)',(\CV^O)',\CP_1',\CP_2',\ldots,\CP_j'\}$, $\{\CE'',\CJ''\}$, $\{(\CV^I)'',(\CV^O)'',\CP_1'',\CP_2'',\ldots, \CP_k''\}$
all contain mutually orthogonal spaces.%
\index{orthogonal spaces}
For any two fields
\beq \BP=\BP^I+\BP^O+\BP'+\BP'',\quad \BQ=\BQ^I+\BQ^O+\BQ'+\BQ'', \eeq{I0.20}
in the vector space $\CK$, with 
\beq \BP^I,\BQ^I \in (\CV^I)',\quad \BP^O,\BQ^O \in (\CV^O)'', \quad \BP',\BQ'\in \CH',\quad \BP'',\BQ''\in \CH'', \eeq{I0.21}
let us define the inner product%
\index{inner product}
of them to be
\beq (\BP,\BQ)=(\BP^I,\BQ^I)'+(\BP^O,\BQ^O)''+(\BP',\BQ')'+(\BP'',\BQ'')'', \eeq{I0.22}
in which $(~,~)'$ and $(~,~)''$ denote the inner product on the spaces $\CK'$ and $\CK''$ respectively. It is immediately clear from this definition
that the subspaces $(\CV^I)'$, $(\CV^O)''$, $\CP_1'$, $\CP_2'$, $\ldots$,$\CP_j'$, $\CP_1''$, $\CP_2''$, $\ldots$, $\CP_k''$ are mutually orthogonal in the new superfunction. 
Now take a field $\BE\in\CE$ and $\BJ\in\CJ$. By the definition of these subspaces there must exist fields $\BE'\in\CE'$ and $\BE''\in\CE''$ 
such that \eq{I0.4} to \eq{I0.6} hold, and fields $\BJ'\in\CJ'$, $\BJ''\in\CJ''$ such that \eq{I0.7} to \eq{I0.9} hold. Consequently we have
\beqa (\BJ,\BE) & = & (\BJ'+\BJ''-(\BJ^O)'-(\BJ^I)'',\BE'+\BE''-(\BE^O)'-(\BE^I)'') \nonum
 & = & ((\BJ^O)',(\BE^O)')'+((\BJ^I)'',(\BE^I)'')''-(\BJ',(\BE^O)')'-(\BJ'',(\BE^I)'')''-((\BJ^O)',\BE')'-((\BJ^I)'',\BE'')'' \nonum
& = & -((\BJ^O)',(\BE^O)')'-((\BJ^I)'',(\BE^I)'')'' \nonum
& = & -((\BJ^O)',(\BE^O)')'-(\BM^J(\BJ^O)',\BM^E(\BE^O)')'' \nonum
& = & -((\BJ^O)',(\BE^O)')'-((\BM^E)^\dagger\BM^J(\BJ^O)',(\BE^O)')',
\eeqa{I0.23}
in which $(\BM^E)^\dagger$ is the adjoint of $\BM^E$.
So we see that the spaces $\CJ$ and $\CE$ will be orthogonal if we choose
\beq (\BM^E)^\dagger\BM^J=-\BI. \eeq{I0.24}
Note that the orthogonality 
of the spaces $\CJ$ and $\CE$ immediately implies that they have no nonzero vector in their intersection.

\begin{sloppypar}
In the case of nonorthogonal subspace collections, we are free to choose the maps  $\BM^E$ and $\BM^J$ that map $(\CV^O)'$ to $(\CV^I)''$,
so long as they and the map $\BA$ are nonsingular. However, in view of \eq{I0.24}, it would be quite natural to restrict our
definition of multiplication by requiring that $\BM^J=-\BM^E$, i.e. one can pick a nonsingular map $\BM$ mapping $(\CV^O)'$ to $(\CV^I)''$
and set
\beq \BM^E=\BM,\quad \BM^J=-\BM. \eeq{I0.24a}
With this choice, subtracting the equations in \eq{I0.11e} gives
\beq  (\BGP^I)''(\BGG_1''-\BGG_2'')\BM(\BGP^O)'(\BGG_1'-\BGG_2')\BP=\BM\BP
\eeq{I0.24b}
Returning to the case where the subspaces are orthogonal, \eq{I0.24} is satisfied if $\BM\BM^\dagger=-\BI$. An alternative way
to see that $\CJ$ and $\CE$ have no nonzero vector in their intersection is as follows. Choose an orthonormal basis $(\Bv_1^O)',(\Bv_2^O)',\ldots (\Bv_{m}^O)'$ 
for $(\CV^O)'$ and take $\BM^E=-\BM^J$ as a map such that
$\BM^E(\Bv_1^O)',\BM^E(\Bv_2^O)',\ldots \BM^E(\Bv_{m}^O)'$ form an orthonormal basis for $(\CV^I)''$. Then the operator $\BM^E$ is represented as the
identity matrix in the basis, and$\BM^J$ is represented by $-\BI$. Now, recalling the definition of the norm $|\BQ|=(\BQ,\BQ)^{1/2}$
of a vector $\BQ$ recall that the action of the operators $(\BGP^O)'$, $(\BGP^I)''$ cannot increase the norm of a vector, while 
$\BGG_1'-\BGG_2'$ and $\BGG_1''-\BGG_2''$ preserve the norm (as can be seen if we take a basis where these are diagonal). Hence \eq{I0.24b} 
can be satisfied only when there is a $\BP\in (\CV^O)'$ such that
\beq (\BGG_1'-\BGG_2')\BP\in (\CV^O)' \quad (\BGG_1''-\BGG_2'')\BM\BP\in (\CV^I)''. \eeq{I0.24b0}
Then as $\BGG_1'+\BGG_2'=\BI$ and $\BGG_1''+\BGG_2''=\BI$ we obtain
\beq (\BGG_1'+\BGG_2')\BP\in (\CV^O)' \quad (\BGG_1''+\BGG_2'')\BM\BP\in (\CV^I)''. \eeq{I0.24b1}
Adding and substracting \eq{I0.24b0} and \eq{I0.24b1} then implies
\beq \BGG_1'\BP\in (\CV^O)',\quad \BGG_2'\BP\in (\CV^O)',\quad \BGG_1''\BP\in (\CV^I)'',\quad \BGG_2''\BP\in (\CV^I)''
\eeq{I0.24b2}
which is excluded by our assumption that $\CV'$ has no vector in common with $\CE'$ or $\CJ'$ and that 
$\CV''$ has no vector in common with $\CE''$ or $\CJ''$.
\end{sloppypar}

\section{Multiplicative identity superfunctions}
\setcounter{equation}{0}
Suppose we are given  nonsingular maps $\BM^E$ and $\BM^J$ which map the $m$-dimensional space $(\CV^O)'$ to  the $m$-dimensional space
$(\CV^I)''$. Let $\CK''$ denote the $2m$-dimensional space 
\beq \CK''=(\CV^O)'\oplus(\CV^I)''. \eeq{I0.24ba}
Within this space define $\CE''$ as the subspace consisting of all vectors of the form $\BE=\Bv+(\BM^E)^{-1}\Bv$ 
with $\Bv\in(\CV^I)''$ and define $\CJ''$ as the subspace consisting of all vectors of the form $\BJ=\Bw+(\BM^J)^{-1}\Bw$ 
with $\Bw\in(\CV^I)''$. If these subspaces have a vector in common then
\beq \Bv+(\BM^E)^{-1}\Bv=\Bw+(\BM^J)^{-1}\Bw,\quad {\rm i.e.,}~\Bv-\Bw=(\BM^J)^{-1}\Bw-(\BM^E)^{-1}\Bv.
\eeq{I0.24c}
In this last equation the fields on the left and on the right lie respectively in $(\CV^I)''$ and $(\CV^O)'$. As the
intersection of these subspaces consists of only the zero vector, we conclude that both sides must be zero, i.e. $\Bw=\Bv$
and  
\beq \Bu\equiv=(\BM^E)^{-1}\Bv=(\BM^J)^{-1}\Bv \eeq{I0.24d}
Thus, $\BM^E\Bu=\Bv=\BM^j\Bu$ and if we assume that $\BM^J-\BM^E$ is nonsingular, then $0=\Bu=\Bv=\Bw$. So under this assumption the
subspaces have only the zero vector in their intersection. Then, since they each have dimension  $m$ we conclude that 
\beq \CK''=(\CV^O)'\oplus(\CV^I)''=\CE''\oplus\CJ'', \eeq{I0.24e}
which defines a superfunction $(F^s)''$ in which $\CH$ is empty.

We now look at the associated superfunction problem. As the space $\CH$ is empty, if we are given vectors $\BE^I$ and $\BJ^I$ in
the input space $(\CV^I)''$ the superfunction problem then consists of finding vectors  $\BE^O$ and $\BJ^O$ in
the output space $(\CV^O)'$ such that
\beq \BE^I+\BE^O\in\CE'',\quad \BJ^I+\BJ^O\in\CJ''. \eeq{I0.24f}
From our definition of the subspaces $\CE''$ and $\CJ''$ we immediately see that the superfunction problem is solved with
fields
\beq \BE^0=(\BM^E)^{-1}\BE^I,\quad\BJ^0=(\BM^J)^{-1}\BJ^I, \eeq{I0.24g}
implying, through \eq{I2.7d}, that the associated function is
\beq \BF''=\begin{pmatrix} (\BM^E)^{-1} & 0  \\ 0 & {\BM^J}^{-1} \end{pmatrix}. \eeq{I0.24h}

So if we take another superfunction $(F^s)'$ and multiply it by this superfunction $(F^s)''$, the product rule \eq{I0.18}
implies that the resulting superfunction $F^s$ has the associated function
\beq \BF=\BF'. \eeq{I0.24i}
We conclude that this superfunction $(F^s)''$ is the multiplicative identity, when multiplication is defined with the maps
$\BM^E$ and $\BM^J$.

\section{Addition of $Y$-subspace collections and embeddings}
\begin{figure}[!ht]
\centering
\includegraphics[width=0.75\textwidth]{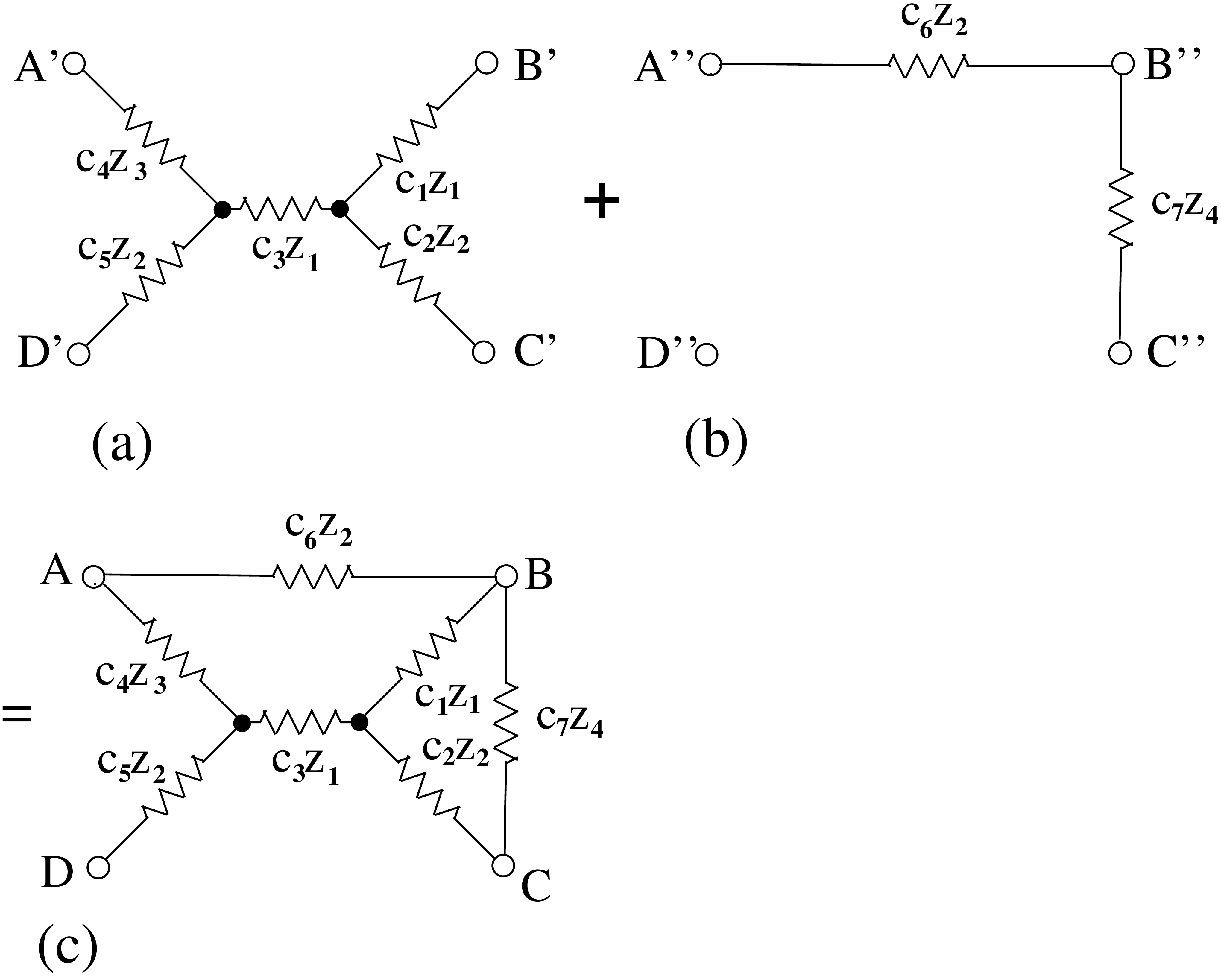}
\caption[Adding $Y$-subspace collections is like hooking networks together in parallel.]{Adding $Y$-subspace collections is like hooking networks together in parallel.
The 4 terminal networks in (a) and (b), each representing (along with their respective tree-like 
battery configurations on the opposite side of the circuit board that are not shown here) $Y(3)$ and $Y(2)$ subspace
collections, are added together to form the 4 terminal network in (c) which is a $Y(4)$ subspace collection. 
Note that the circuit in (b) is really only a 3 terminal network, so it has been embedded in a 4 terminal network
(with no electrical connections to the 4th terminal). Also note there is some flexibility in how one adds together
these subspace collections: we connected the terminals $A',B',C',$ and $D'$, to respectively the terminals $A'',B'',C'',$ and $D''$,
but we could have connected them to any permutation of these terminals. This flexibility is reflected in the need to introduce
nonsingular operators $\BS'$ and $\BS''$ which respectively map $\CV'$ and $\CV''$ to $\CV$, before addition can defined.}\label{IS3}
\end{figure}

\labsect{IadditionY}
\setcounter{equation}{0}
Adding superfunctions%
\index{superfunctions!addition}
is similar the way electrical circuits, each with $n$ terminals can be combined. An example
is shown in Figure~\ref{IS3}.

Suppose we have $Y(j)$ and $Y(k)$ subspace collections:
\beqa \CK' & = & \CE'\oplus\CJ'=\CV'\oplus\CH'\quad{\rm with}~~
\CH'=\CP_1'\oplus\CP_2'\oplus\cdots\oplus\CP_j', \nonum
 \CK'' & = & \CE''\oplus\CJ''=\CV''\oplus\CH''\quad{\rm with}~~
\CH''=\CP_1''\oplus\CP_2''\oplus\cdots\oplus\CP_k'',
\eeqa{I0.25}
where the spaces $\CV'$ and $\CV''$ have the same dimension $n$. To define the sum of the 
subspace collections we need to introduce another $n$-dimensional space $\CV$ and nonsingular operators $\BS'$ and
$\BS''$ which respectively map $\CV'$ and $\CV''$ to $\CV$. Then the sum of the subspace collections
\beq \CK=\CK'+_{\{\BS',\BS''\}}\CK'' \eeq{I0.26}
is the subspace collection
\beq \CK=\CE\oplus\CJ=\CV\oplus\CH, \eeq{I0.27}
where
\beq \CH=\CH'\oplus\CH''=\CP_1'\oplus\CP_2'\oplus\cdots\oplus\CP_j'\oplus\CP_1''\oplus\CP_2''\oplus\cdots\oplus\CP_k''.
\eeq{I0.28}
Here a field $\BE=\BE_1+\BE_2$, with $\BE_1\in\CV$ and $\BE_2\in\CH$,  is in $\CE$ if and only if there exist fields
\beq \BE'=\BE_1'+\BE_2'\in\CE',\quad \BE''=\BE_1''+\BE_2''\in\CE'', \eeq{I0.29}
with
\beq \BE_1'\in\CV',~~\BE_2'\in\CH',~~\BE_1''\in\CV'',~~\BE_2''\in\CH'', \eeq{I0.30}
such that
\beq \BS'\BE_1'=\BS''\BE_1''=\BE_1. \eeq{I0.31}
Also a field $\BJ=\BJ_1+\BJ_2$, with $\BJ_1\in\CV$ and $\BJ_2\in\CH$,  is in $\CJ$ if and only if there exist fields
\beq \BJ'=\BJ_1'+\BJ_2'\in\CE',\quad \BJ''=\BJ_1''+\BJ_2''\in\CE'', \eeq{I0.32}
with
\beq \BJ_1'\in\CV',~~\BJ_2'\in\CH',~~\BJ_1''\in\CV'',~~\BJ_2''\in\CH'', \eeq{I0.33}
such that
\beq \BS'\BJ_1'+\BS''\BJ_1''=\BJ_1. \eeq{I0.34}
So given $\BE_1\in\CV$, we let $\BE_1'=(\BS')^{-1}\BE$ and $\BE_1''=(\BS'')^{-1}\BE_1$, and we solve the $Y$-problem in each of the
two subspace collections $Y(j)$ and $Y(k)$, finding fields satisfying \eq{I0.29}, \eq{I0.30}, \eq{I0.32}, and \eq{I0.33} with
\beq \BJ_2'=\BL'\BE_2,\quad \BJ_2''=\BL''\BE_2'', \eeq{I0.35}
where  
\beq \BL'=\sum_{i=1}^jz_i'\BGL_i',\quad \BL''=\sum_{i=1}^kz_i''\BGL_i'', \eeq{I0.36}
and $\BGL_i'$ projects onto $\CP_i'$ while $\BGL_i''$ projects onto $\CP_i''$. Hence we have
\beq \BJ_2=\BJ_2'+\BJ_2''=\BL(\BE_2'+\BE_2''),~~~{\rm with}~~\BL=\BL'+\BL''. \eeq{I0.37}
Then \eq{I0.34} implies
\beq \BJ_1=\BS'\BJ_1'+\BS''\BJ_1''=\BS'\BY'\BE_1'+\BS''\BY''\BE_1''=\BY\BE_1, \eeq{I0.38}
where
\beq \BY=\BS'\BY'(\BS')^{-1}+\BS''\BY''(\BS'')^{-1}. \eeq{I0.39}
If we have a basis $\Bv_1, \Bv_2, \ldots, \Bv_n$ for $\CV$, then it is natural to take
$(\BS')^{-1}\Bv_1$, $(\BS')^{-1}\Bv_2$, $\ldots$, $(\BS')^{-1}\Bv_n$ as a basis for $\CV'$, and
to take $(\BS'')^{-1}\Bv_1$, $(\BS'')^{-1}\Bv_2$, $\ldots$, $(\BS'')^{-1}\Bv_n$ as a basis for
$\CV''$. Then the operators $\BS'$ and $\BS''$ are represented by identity matrices, and 
in these bases \eq{I0.39} becomes $\BY=\BY'+\BY''$.

In the case where either or both of the subspaces $\CV'$ and $\CV''$ 
have dimension less than the dimension $n$ of the subspace $\CV$ we can first do an embedding.%
\index{subspace collections!embedding} 
For example suppose $\CV'$ has dimension $n'<n$. Then let
$\CW'$ be a space of dimension $n-n'$. Construct the subspace collection
\beq \widetilde{\CK}'= \widetilde{\CE}'\oplus\CJ'=\widetilde{\CV}'\oplus\CH', \eeq{I0.40}
where
\beq \widetilde{\CV}'=\CV'\oplus\CW',\quad \widetilde{\CE}'=\CE'\oplus\CW'. \eeq{I0.41}
Then given a field $\widetilde{\BE}_1'\in\widetilde{\CV}'$ we can express it as a sum $\BE_1'+\BW'$
with $\BE_1'\in\CV'$ and $\BW'\in\CW'$. We write $\BE_1'=\GY\widetilde{\BE}_1'$ where $\GY$ is the projection onto $\CV'$. 
Given this $\BE_1'$ and solving the $Y$-problem%
\index{Y@$Y$-problem}
associated with $\CK'$ we obtain fields $\BE'$ and $\BJ'$ satisfying 
\beqa \BE' & = &\BE_1'+\BE_2'\in\CE',\quad \BE_1'\in\CV',~~\BE_2'\in\CH', \nonum
 \BJ' & = &\BJ_1'+\BJ_2'\in\CJ',\quad \BJ_1'\in\CV',~~\BJ_2'=\BL\BE_2 \in\CH.
\eeqa{I0.42}
It follows that the $Y$-problem in the space $\widetilde{\CK}'$ is solved with fields
\beq \widetilde{\BE}'=\BW+\BE'=\BW+\BE_1'+\BE_2',\quad {\rm and}~~~  \BJ'=\BJ_1'+\BJ_2'~~{\rm with}~~\BJ_2'=\BL\BE_2,  \eeq{I0.43}
implying that
\beq \BJ_1'=-\BY\BE_1'=-\BY\GY\widetilde{\BE}_1'. \eeq{I0.44}
We conclude that the new $Y$-problem has an operator $\widetilde{\BY}=\BY\GY$, i.e. its range is not the whole space $\widetilde{\CV}'$
but only at most the subspace $\CV'$. After making such embeddings to ensure that $\CV'$ and $\CV''$ (or
rather $\widetilde{\CV}'$ and $\widetilde{\CV}''$ have the same dimension as the dimension $n$ of the subspace $\CV$, we are 
then free to add them. 

The additive zero is easy to find. Let us consider the degenerate subspace collection
\beq \CK''=\CE''=\CV'' \eeq{I0.44a}
Clearly $\CH''$ contains only the zero vector, and we are forced to choose $\BL''=0$. Given $\BE_1\in\CV''$. The $Y$-problem is solved with vectors
\beq \BE''=\BE_1,\quad \BE_1=\BJ_1=\BJ_2=\BJ=0. \eeq{I0.44b} 
Implying the associated $Y$-operator $\BY$ is zero: thus the subspace collection \eq{I0.44a} is the additive zero. Note that this subspace collection does not satisfy
the property $\CE''\bigcap\CV''=0$ which is needed for the inverse of $\BY$ to exist, which is not surprising since $\BY=0$ has no inverse.

Now suppose we have a subspace collection
\beq \CK'  =  \CE'\oplus\CJ'=\CV'\oplus\CH'\quad{\rm with}~~
\CH'=\CP_1'\oplus\CP_2'\oplus\cdots\oplus\CP_j', \eeq{I0.44-0}
with associated operator $\BY'(z_1', z_2',\ldots,z_n')$ when
\beq \BL'=\sum_{i=1}^jz_i'\BGL_i'.  \eeq{I0.44-1}
It is clear that if we replace $\BL'$ by
\beq \BL'=-\sum_{i=1}^jz_i'\BGL_i', \eeq{I0.44-2}
then the solution to the $Y$-problem will give the $\BY$-operator
\beq \BY(-z_1',-z_2',\ldots,-z_n')= -\BY(z_1',z_2',\ldots,z_n'), \eeq{I0.44-3}
where to obtain this last identity we have used the homogeneity of the function. Since
adding \eq{I0.44-3} to the associated operator $\BY'(z_1', z_2',\ldots,z_n')$ we started with gives zero,
it is tempting to conclude that we have found the additive inverse. However the function \eq{I0.44-3} is not
the $\BY$-operator valued function of $z_1', z_2',\ldots,z_n'$ associated with the subspace collection \eq{I0.44-0}, whose 
definition does not allow us to choose $\BL'$ of the form \eq{I0.44-2}. This is made more clear in the case where we have an orthogonal subspace collection
since then the imaginary part of $(\BV,\BY(z_1', z_2',\ldots,z_n')\BV)$ is generally positive when $z_1', z_2',\ldots,z_n'$ all have positive imaginary parts, and 
$-\BY(z_1', z_2',\ldots,z_n')$ then does not share this Herglotz property. So the additive inverse of an orthogonal subspace collection should typically not
be an orthogonal subspace collection. We will find the proper additive inverse in section \sect{addinv}.

\section{Substitution of subspace collections}%
\index{subspace collections!substitution}
\labsect{Isubsubstit}
\setcounter{equation}{0}

\begin{figure}[!ht]
\centering
\includegraphics[width=0.75\textwidth]{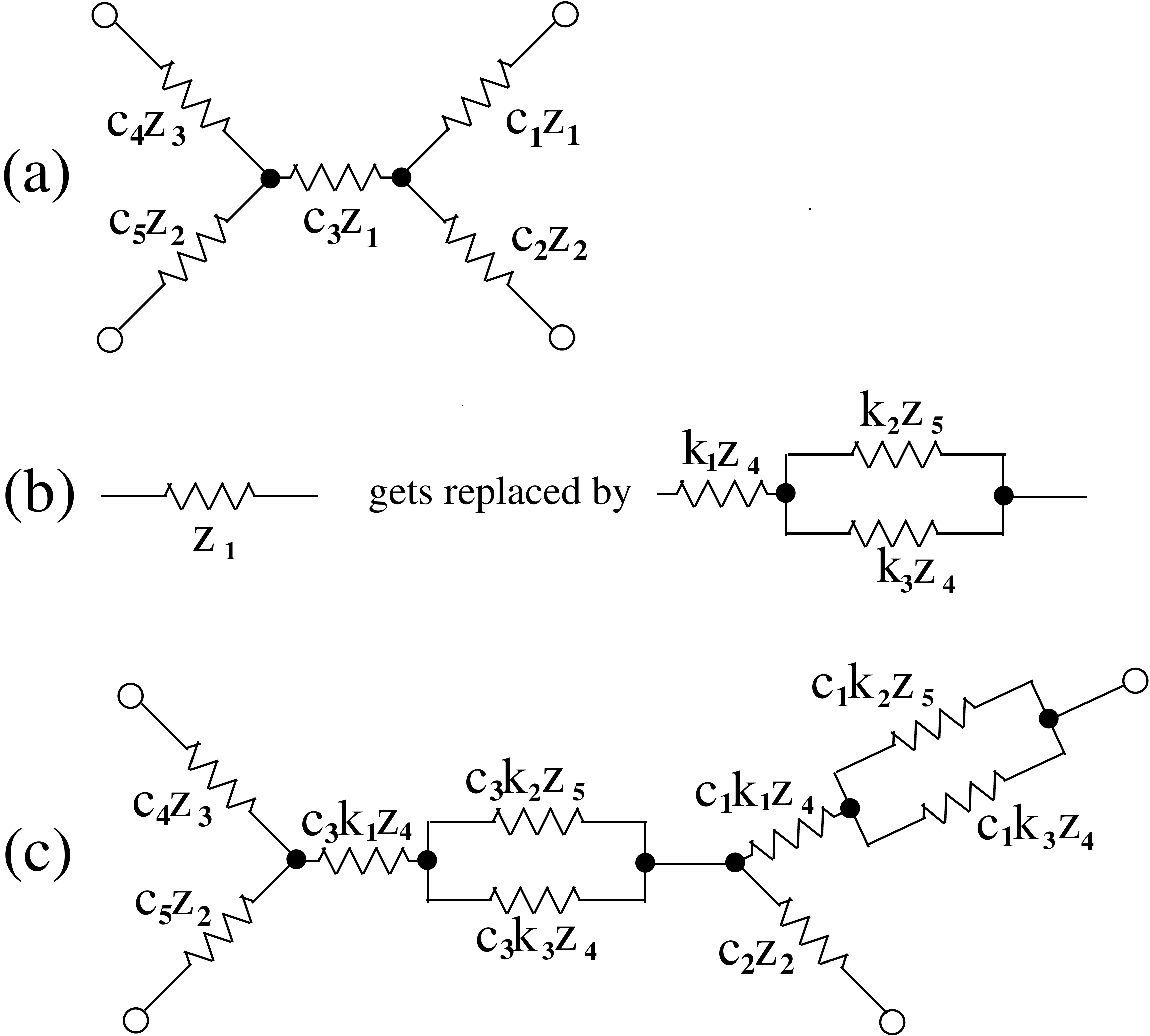}
\caption[Substitution of  $Y$- and $Z$-subspace collections is like replacing all resistors of one type by a compound network.]{Substitution of  $Y$- and $Z$-subspace collections is like replacing all resistors of one type by a compound network.
If one takes a subspace collection, as, for example, represented by the 4-terminal network in (a) 
and replaces $z_1$ by the network in (b), where $k_1+(1/k_2+1/k_3)^{-1}=1$, to ensure this replacement does
effect the resistance when $z_1=z_4=z_5=1$, one obtains the subspace collection as represented by the 4-terminal 
network in (c).}\label{IS4}
\end{figure}
Another familiar operation that we can do with rational functions is to make substitutions. Substitution of one subspace collection in another is similar to the way it can be done
in electrical circuits. An example
is shown in figure \ref{IS4}.
Thus if $\BY(z_1,z_2, \ldots,z_n)$ is a $m\times m$ matrix-valued homogeneous
function of degree one and $Z'(z_1',z_2', \ldots,z_p')$ is a scalar-valued
function, say normalized with
\beq Z'(1,1, \ldots,1)=1, \eeq{I8.0}
then
\beq \BY''(z_1',z_2', \ldots,z_p',z_2, \ldots,z_n)
=\BY(Z(z_1',z_2', \ldots,z_p'),z_2, \ldots,z_n) \eeq{I8.1}
will be another $m\times m$ matrix-valued homogeneous 
function of degree one. What is the analogous operation on subspace collections? It is natural to expect
there should be one, just as in a network of $n$ types of resistors one can replace each resistor of type $1$
with a network of $p$ other resistors.

Extending the treatment given in Section 29.1 of \citeAPY{Milton:2002:TOC} let us suppose that we are given
a $Y(n)$-subspace collection
\beq \CK=\CE\oplus\CJ=\CV\oplus\CP_1\oplus\CP_2\oplus\cdots\oplus\CP_n, \eeq{I8.2}
and a $(3,p)$-subspace collection
\beq \CH'=\CU'\oplus\CE'\oplus\CJ'=\CP'_1\oplus\CP'_2\oplus\cdots\oplus\CP'_p, \eeq{I8.2b}
in which $\CV$ is $m$-dimensional and
$\CU'$ is one-dimensional. Let  $\BY(z_1,z_2, \ldots,z_n)$ and
$Z'(z_1',z_2', \ldots,z_p')$ denote the functions associated with these subspace collections.
We take as our new $(2,n+p)$-subspace collection,
\beq \CK''=\CE''\oplus\CJ''=\CV''\oplus\CP_1''\oplus\CP_2''\oplus\cdots\oplus\CP_n'',
\eeq{I8.2c}
where
\beq  \CE''=(\CE\otimes\CU')\oplus(\CP_1\otimes\CE'),~~~
\CJ''=(\CJ\otimes\CU')\oplus(\CP_1\otimes\CJ'), \eeq{I8.3}
and
\beqa  \CV''& = & \CV\otimes\CU', \nonum
\CP_i''& = &\CP_1\otimes\CP'_i~~~{\rm for~~}1\leq i \leq p, \nonum
             & = &\CP_{i+1-p}\otimes\CU'~~~{\rm for~~}p+1\leq i \leq n+p-1.
\eeqa{I8.4}
in which $\otimes$ denotes the operation of taking the tensor product%
\index{tensor product}
of two subspaces. 
Vectors in the space
\beq \CK''=\CE''\oplus\CJ''
        = (\CK\otimes\CU')\oplus(\CP_1\otimes(\CE'\oplus\CJ')) \eeq{I8.5}
spanned by these subspaces are represented as a pair $[\BP,~\Bu']$
added to a linear combination of pairs of the form $[\BP_1,~\BP']$, where
$\BP\in\CK$, $\Bu'\in \CU'$, $\BP_1\in \CP_1$, and $\BP'\in\CE'\oplus\CJ'$.

Now define
\beqa \CH & = &\CP_1\oplus\CP_2\oplus\cdots\oplus\CP_n, \nonum
\CH''& = &\CP_1''\oplus\CP_2''\oplus\cdots\oplus\CP_n'',
\eeqa{I8.5a}
and suppose that we are given solutions to the equations
\beqa \BJ_2 & = &\sum_{i=1}^nz_i\BGL_i\BE_2~~{\rm with}~\BE_1+\BE_2\in\CE,~~
\BJ_1+\BJ_2\in\CJ,\quad \BE_1,\BJ_1\in\CV,\quad \BE_2,\BJ_2\in\CH, \nonum
 \Bj'+\BJ' & = &\sum_{j=1}^nz_j'\BGL'_j(\Be'+\BE')~~{\rm with}~\Be',\Bj'\in\CU',~~\BE'\in\CE',~~
\BJ'\in\CJ',\nonum
& ~& 
\eeqa{I8.6}
where
\beq z_1=Z(z_1',z_2', \ldots,z_p'), \eeq{I8.7}
while $\BGL_i$ and $\BGL'_j$ are the projections onto $\CP_i$ and $\CP'_j$.
Let us introduce
\beq \BP_i=\BGL_i\BE_2,\quad \BP_j'=\BGL_j'(\Be'+\BE'), \eeq{I8.7a}
and set
\beqa z_i'' & = & z_i'~~~{\rm for~~}1\leq i \leq p, \nonum
	      & = & z_{i+1-p}~~~{\rm for~~}p+1\leq i \leq n+p-1. \eeqa{I8.8}
Then, in the new subspace collection, the vectors
\beqa \BE''_1 & = & [\BE_1,~\Be']\in\CV'', \quad \BE''_2=[\BE_2,~\Be']+[\BP_1,~ \BE'],\nonum
\BJ''_1 & = & [\BJ_1,~\Be']\in\CV'', \quad \BJ''_2=[\BJ_2,~\Be']+[\BP_1,~\BJ']
\eeqa{I8.9}
satisfy
\beq \BE''_1+\BE''_2\in\CE'',\quad \BJ''_1+\BJ''_2\in\CJ''. \eeq{I8.10}
Additionally, we have
\beq \BE''_2=[\sum_{i=1}^n\BP_i,~\Be']-[\BP_1,~\Be']+[\BP_1,~\sum_{i=1}^p\BP'_i]=\sum_{i=1}^{n+p-1}\BP''_i\in\CH'',
\eeq{I8.11}
where
\beqa \BP''_i& =&  [\BP_1,~\BP'_i]~~{\rm for}~~1\leq i\leq p, \nonum
             & = & [\BP_{i+1-p},~\Be']~~{\rm for}~~p+1\leq i \leq n+p-1
\eeqa{I8.12}
satisfies $\BP''_i\in\CP''_i$.      
Similarly, and using the fact implied by \eq{I8.7} that $\Bj'=Z\Be'=z_1\Be'$, we have
\beq  \BJ''_2=[\sum_{i=1}^nz_i\BP_i,~\Be']-[\BP_1,~\Bj']+[\BP_1,\sum_{i=1}^pz'_i\BP'_i]=\sum_{i=1}^{n+p-1}z''_i\BP''_i\in\CH''.
\eeq{I8.13}
Given a basis $\Bv_1,\Bv_2,\ldots,\Bv_m$ for $\CV$ and a vector $\Bu'$ in $\CU'$ it is natural to take 
$(\Bv_1,\Bu')$,$(\Bv_2,\Bu')$, $\ldots$, $(\Bv_m,\Bu')$ as our basis for $\CV''$. Choosing $\Be'$ so that $\Be'=\Bu'$, it is evident
that $\BY(Z'(z_1',z_2', \ldots,z_p'),z_2, \ldots,z_n)$ is the matrix-valued function associated the new subspace collection, represented in these bases.

There is a similar subspace operation corresponding to substituting the $Z$-function $Z'(z_1',z_2', \ldots,z_p')$ into
another $Z$-function $\BZ(z_1,z_2, \ldots,z_n)$ to obtain
\beq \BZ''(z_1',z_2', \ldots,z_p',z_2, \ldots,z_n)
=\BZ(Z(z_1',z_2', \ldots,z_p'),z_2, \ldots,z_n). \eeq{I8.14}
Given a $Z(n)$-subspace collection
\beq \CH=\CU\oplus\CE\oplus\CJ=\CP_1\oplus\CP_2\oplus\cdots\oplus\CP_n, \eeq{I8.15}
and a $(3,p)$-subspace collection
\beq \CH'=\CU'\oplus\CE'\oplus\CJ'=\CP'_1\oplus\CP'_2\oplus\cdots\oplus\CP'_p, \eeq{I8.16}
in which $\CU$ is $m$-dimensional and
$\CU'$ is one-dimensional, we take as our new $(3,n+p-1)$-subspace collection,
\beq \CK''=\CU''\oplus\CE''\oplus\CJ''=\CP_1''\oplus\CP_2''\oplus\cdots\oplus\CP_n'',
\eeq{I8.17}
where
\beq  \CU'' =  \CU\otimes\CU',~~~ \CE''=(\CE\otimes\CU')\oplus(\CP_1\otimes\CE'),~~~
\CJ''=(\CJ\otimes\CU')\oplus(\CP_1\otimes\CJ'), \eeq{I8.18}
and
\beqa 
\CP_i''& = &\CP_1\otimes\CP'_i~~~{\rm for~~}1\leq i \leq p, \nonum
             & = &\CP_{i+1-p}\otimes\CU'~~~{\rm for~~}p+1\leq i \leq n+p-1.
\eeqa{I8.19}
Suppose that we are given solutions to the equations
\beqa \Bj+\BJ & = &\sum_{i=1}^nz_i\BGL_i(\Be+\BE)~~{\rm with}~\Be,\Bj\in\CU,~~\BE\in\CE,~~\BJ\in\CJ,\nonum
 \Bj'+\BJ' & = &\sum_{j=1}^nz_j'\BGL'_j(\Be'+\BE')~~{\rm with}~\Be',\Bj'\in\CU',~~\BE'\in\CE',~~\BJ'\in\CJ',\nonum
& ~& 
\eeqa{I8.20}
where $z_1=Z(z_1',z_2', \ldots,z_p')$, 
while $\BGL_i$ and $\BGL'_j$ are the projections onto $\CP_i$ and $\CP'_j$.
Let us introduce
\beq \BP_i  =  \BGL_i(\Be+\BE),\quad \BP_j'=\BGL_j'(\Be'+\BE'), \nonum
\eeq{I8.21}
and define $z_i''$ by \eq{I8.8}, and $\BP''_i\in\CP''_i$ by \eq{I8.12}.
In the new subspace collection, the vectors
\beqa \Be'' & = & [\Be,~\Be']\in\CU'',\quad \BE'' =  [\BE,~\Be']+[\BP_1,~\BE']\in\CE'',\nonum
\Bj'' & = & [\Bj,~\Be']\in\CU'', \quad \BJ''=[\BJ,~\Be']+[\BP_1,~\BJ']\in\CJ''
\eeqa{I8.22}
satisfy
\beqa \Be''+\BE'' & = & [\sum_{i=1}^n\BP_i,~\Be']+[\BP_1,~\sum_{j=1}^p\BP_j']-[\BP_1,~\Be'] \nonum
& = &  [\sum_{i=2}^n\BP_i,~\Be']+[\BP_1,~\sum_{j=1}^p\BP_j'] \nonum
 & = &\sum_{i=1}^{n+p-1}\BP_i'',
\eeqa{I8.23}
and, using \eq{I8.7},
\beqa \Bj''+\BJ'' & = & [\sum_{i=1}^nz_i\BP_i,~\Be']+[\BP_1,~\sum_{j=1}^pz_i'\BP_j']-[\BP_1,~\Bj'] \nonum
& = &  [\sum_{i=2}^nz_i\BP_i,~\Be']+[\BP_1,~\sum_{j=1}^pz'_i\BP_j'] \nonum
& = &\sum_{i=1}^{n+p-1}z_i''\BP_i''.
\eeqa{I8.24}
Given a basis $\Bu_1,\Bu_2,\ldots,\Bu_m$ for $\CU$ and a vector $\Bu'$ in $\CU'$ it is natural to take 
$(\Bu_1,\Bu')$,$(\Bu_2,\Bu')$, $\ldots$, $(\Bu_m,\Bu')$ as our basis for $\CU''$. Choosing $\Be'$ so that $\Be'=\Bu'$, it is evident from \eq{I8.22}
that $\BZ(Z'(z_1',z_2', \ldots,z_p'),z_2, \ldots,z_n)$ is the matrix-valued function associated the new subspace collection, represented in these bases.

\section{Some other elementary operations on subspace collections}
\setcounter{equation}{0}

A further operation we can do on functions $\BY(z_1,z_2, \ldots,z_n)$ while retaining
the homogeneity of degree 1 in the variables $z_1,z_2, \ldots,z_n$ is to replace
the function by $[\BY(1/z_1,1/z_2, \ldots,1/z_n)]^{-1}$. The analogous operation 
on the associated $Y(n)$-subspace collection is to interchange the subspaces $\CE$ and $\CJ$.
Similarly in a $Z(n)$ subspace collection, interchanging the subspaces $\CE$ and $\CJ$
corresponds to replacing $\BZ(z_1,z_2, \ldots,z_n)$ by $[\BZ(1/z_1,1/z_2, \ldots,1/z_n)]^{-1}$, as noted in
Section 29.1 of \citeAPY{Milton:2002:TOC}.
We call such a transformation a duality transformation.%
\index{duality trasformation}
As a consequence of the duality transformation \eq{I4.10} immediately implies the formula
\index{Zo@$\BZ$-operator formula}
\beq \BZ^{-1}=\BGG_0[(\BGG_0+\BGG_1)\BL(\BGG_0+\BGG_1)]^{-1}\BGG_0.
\eeq{I8.25}
One simple thing we can do in a function
is set $z_j=z_k$:  the analogous operation in a subspace collection is to replace $\CP_j\oplus\CP_k$ by a single subspace. 

Another elementary operation we can do on a $Z(n)$ subspace collection is as follows. Let $\CU$ be expressed as the direct sum
\beq \CU=\CU'\oplus\CW, \eeq{I8.25a} 
which defines the projection $\BGF$ onto $\CU'$. We now take as our subspace collection
\beq \CH=\CU'\oplus\CE\oplus\CJ'=\CP_1\oplus\CP_2\oplus\cdots\oplus\CP_n,  \eeq{I8.25b}
where 
\beq \CJ'=\CJ\oplus\CW. \eeq{I8.25c}
Then any solution to the $Z$-problem \eq{I2.11} with $\Be\in\CU'$ immediately generates a solution to the $Z$-problem associated with the
subspace collection \eq{I8.25b}:
\beq \Bj'\in\CU',~~~~\BE\in\CE,~~~~\BJ'\in\CJ',~~~~
\Bj'+\BJ'=\BL(\Be+\BE), \eeq{I8.25d}
where
\beq \BL=\sum_{i=1}^nz_i\BGL_i,\quad \Bj'=\BGF\Bj,\quad \BJ'=\BJ+(\BI-\BGF)\Bj, \eeq{I8.25e}
which ensures that
\beq \Bj+\BJ=\Bj'+\BJ'\quad{\rm and}\quad (\BI-\BGF)\Bj\in\CW. \eeq{I8.25f}
Hence the new subspace collection has a $\BZ$-operator
\beq \BZ'=\BGF\BZ, \eeq{I8.25g}
when applied to fields in $\CU'$. 

\section[Realizing any rational $Y$-matrix]{Realizing any $Y$-matrix with elements that are rational functions of degree $1$}
\setcounter{equation}{0}

Given any homogeneous rational function of degree 1,
\beq Z(z_1,z_2,\ldots,z_n)=\frac{p(z_1,z_2,\ldots,z_n)}{q(z_1,z_2,\ldots,z_n)},
\eeq{I10.1}
satisfying the normalization $Z(1,1,\ldots,1)=1$
where $p(z_1,z_2,\ldots,z_n)$ and $q(z_1,z_2,\ldots,z_n)$ are homogeneous 
polynomials of degree $k$ and $k-1$ respectively, where $k$ is a positive integer,
our first goal is to find a $Z(n)$ subspace collection
\beq \CH=\CU\oplus\CE\oplus\CJ=\CP_1\oplus\CP_2\oplus\cdots\oplus\CP_n, \eeq{I10.2}
where $\CU$ is one-dimensional which has $Z(z_1,z_2,\ldots,z_n)$ as its associated
function. Without loss of generality we could set $z_n$ =1, and then 
$p(z_1,z_2,\ldots,z_{n-1},1)$ and $q(z_1,z_2,\ldots,z_{n-1},1)$ are just arbitrary
polynomials of the $n-1$ variables $z_1,z_2,\ldots,z_{n-1}$. Also without loss of generality
we can assume
\beq p(1,1,\ldots,1)=q(1,1,\ldots,1)=1. \eeq{I10.2a}
 The first step is to realize $Z(z_1,z_2,1)=z_1z_2$ as an associated $Z$-function. Note that \eq{I3.16} implies we can realize
\beq Z(z_1,1)=z_1^2, \eeq{I10.3}
and \eq{I3.17} implies we can realize
\beq Z(z_1,z_2)=cz_1+(1-c)z_2, \eeq{I10.4}
for any constant $c$. Hence, by substitution we can realize
\beq Z(z_1,z_2,1)=9(2z_1/3+z_2/3)^2/8-(2z_1-z_2)^2/8=z_1z_2. \eeq{I10.5}
Making further substitutions, we can realize any product of the variables
\beq Z(z_1,z_2,\ldots,z_{n-1},1)=z_1^{a_1}z_2^{a_2}\ldots z_{n-1}^{a_{n-1}},  \eeq{I10.6}
where the $a_i$ are nonnegative integers. By 
repeated substitution in \eq{I10.4} we can realize
any linear combination of such terms with coefficients summing to 1. Thus we can realize the polynomials%
\index{polynomials!realizing}
$p(z_1,z_2,\ldots,z_{n-1},1)$ and $q(z_1,z_2,\ldots,z_{n-1},1)$.

Furthermore \eq{I3.16}, with the roles of $z_1$ and $z_2$ interchanged, implies we can
realize
\beq Z(z_1,1)=1/z_1, \eeq{I10.6a}
which by substitution into \eq{I10.5} implies we can realize
\beq Z(z_1,z_2,1)=z_2/z_1. \eeq{I10.7}
Substituting $p(z_1,z_2,\ldots,z_{n-1},1)$ for $z_2$ and $q(z_1,z_2,\ldots,z_{n-1},1)$
for $z_1$ we see we can find a subspace collection which realizes
\beq Z(z_1,z_2,\ldots,z_{n-1},1)=\frac{p(z_1,z_2,\ldots,z_{n-1},1)}{q(z_1,z_2,\ldots,z_{n-1},1)}
\eeq{I10.8}
as its associated $Z$-function when $z_n=1$. When $z_n$ is not 1, the subspace collection will by homogeneity
realize the function \eq{I10.1}.

Now from \eq{I3.6g} we can realize 
\beq \BY(z_1)= \begin{pmatrix} a_{11}z_1  &  0 \\ 0  &  0 \end{pmatrix}, \eeq{I10.8a}
and realize
\beq 
     \BY(z_2)=\begin{pmatrix} 0  &  a_{12}z_2 \\ 0  &  0 \end{pmatrix}.
\eeq{I10.8b}
By substitution of subspace collections, we can realize any $Y$-matrix where in the above formulae $z_1$ and $z_2$ are replaced
by any normalized rational homogeneous functions of degree $1$ (normalized in the sense that they take the value $1$ when all variables take the value $1$). 
Finally, by making suitable embeddings and adding subspace collections we can realize any $Y$-matrix with elements that are 
homogeneous rational functions%
\index{rational functions!realizing}
of degree $1$: \eq{I10.8a} with the appropriate substitutions realizes each diagonal element, while \eq{I10.8b}  
with the appropriate substitutions realizes each off-diagonal element.

\section{Extension operations on subspace collections}
\index{subspace collections!extension operations}
\setcounter{equation}{0}

Let us suppose we have a $Z(n)$ subspace collection
\beq \CH=\CU\oplus\CE\oplus\CJ=\CP_1\oplus\CP_2\oplus\cdots\oplus\CP_n,
\eeq{I5.1}
where $\CU$ is $m$-dimensional. Let $\CV$ be another $m$-dimensional space, and
consider the space
\beq \CK=\CV\oplus\CH. \eeq{I5.2}
Suppose there is a nonsingular mapping $\BT$ from $\CU$ to $\CV$. Define the 
subspace $\widetilde{\CE}$ to consist of all vectors spanned by $\Bu+\BT\Bu$ as 
$\Bu$ varies in $\CU$. Define $\widetilde{\CJ}$ to consist of all vectors spanned by
$\Bu-\BT\Bu$ as $\Bu$ varies in $\CU$. Clearly we have
\beq \CV\oplus\CU=\widetilde{\CE}\oplus\widetilde{\CJ}, \eeq{I5.3}
and consequently we obtain the $Y(n)$ subspace collection
\beq \CK=\CE'\oplus\CJ'=\CV\oplus\CP_1\oplus\CP_2\oplus\cdots\oplus\CP_n,
\eeq{I5.4}
in which
\beq \CE'=\widetilde{\CE}\oplus\CE,\quad  \CJ'=\widetilde{\CJ}\oplus\CJ.
\eeq{I5.5}
Furthermore given vectors satisfying
\beq \Bj+\BJ=\BL(\Be+\BE),\quad \BE\in\CE,~~\BJ\in\CJ,~~\Be,\Bj\in\CU, \eeq{I5.6}
where
\beq \Bj=\BZ\Be,\quad\BL=\sum_{\ell=1}^{n} z_i\BGL_i, \eeq{I5.7}
we can set
\beq \BE_2=\Be+\BE\in\CH,\quad \BE_1=\BT\Be\in\CV,
\quad \BJ_2=\Bj+\BJ\in\CH,\quad \BJ_1=-\BT\Bj.
\eeq{I5.8}
Then we have
\beq \BE_1+\BE_2=\BT\Be+\Be+\BE\in\CE',\quad \BJ_1+\BJ_2=-\BT\Bj+\Bj+\BJ\in\CJ',
\eeq{I5.9}
and
\beq \BJ_1=-\BY\BE_1,\quad{\rm with}~~\BY=\BT\BZ\BT^{-1}.
\eeq{I5.9a}
Given a basis $\Bu_1, \Bu_2,\ldots,\Bu_m$ for $\CU$, with respect to which the
matrix valued function $\BZ(z_1,z_2,\ldots,z_n)$ is defined,
it is natural to take $\BT\Bu_1, \BT\Bu_2,\ldots,\BT\Bu_m$ as our basis
for $\CV$. Then $\BT$ is represented by the identity matrix, and 
the functions $\BZ(z_1,z_2,\ldots,z_n)$ and  $\BY(z_1,z_2,\ldots,z_n)$ are identical.
We call the subspace collection \eq{I5.4} the extension of the subspace collection \eq{I5.1}.

\section{Reference transformations and additive inverses}
\index{reference transformations}
\labsect{addinv}
\setcounter{equation}{0}

Given the impedance network illustrated in Figure~\ref{IS1} we are free the change the scaling constants $c_i$ assigned to
each bond to new constants $c_i'$ and accordingly replace $z_i$ with $z_i'=z_ic_i/c_i'$ without changing the overall electrical response of the network. 
Analogously, given a homogeneous rational function $\BY(z_1,z_2, \ldots,z_n)$ of degree one, an operation 
which preserves the homogeneity is obviously to multiply the variables by constants to obtain 
the function
\beq \BY'(z'_1,z'_2, \ldots,z'_n)= \BY(d_1z_1',d_2z_2', \ldots,d_nz_n').
\eeq{I6.0}
The associated operation on the $Y(n)$ subspace collection  $(\CE,\CJ)$
and  $(\CV,\CP_1,\CP_2,\ldots,\CP_n)$ is found by generalizing the analysis given after (29.3) in \citeAPY{Milton:2002:TOC} and is
as follows. Given nonzero (possibly complex) constants
$c_i^E$ and $c_i^J$, $i=1,\ldots,n$ we  introduce
the linear transformations
\beq \Gy^{E}(\BP)=\BGP_1\BP+\sum_{i=1}^n c_i^{E}\BGL_i\BP,~~~
   \Gy^{J}(\BP)=\BGP_1\BP+\sum_{i=1}^n c_i^{J}\BGL_i\BP, \eeq{I6.1}
on fields $\BP\in\CK$, where $\BGL_1$ is the projection onto $\CP_1$. These transformations leave
the subspaces $\CV$ and $\CP_i$ invariant. Define the spaces
\beq \CE'=\Gy^{E}(\CE)~~{\rm and}~~\CJ'=\Gy^{J}(\CJ). \eeq{I6.2}
These will have the same dimension as $\CE$ and $\CJ$ respectively. To see this, suppose $\Gy^{E}(\BE)=\Gy^{E}(\BE')$
for some $\BE, \BE'\in \CE$. Then $\Gy^{E}(\BE-\BE')=0$ and since \eq{I6.1} implies $\Gy^{E}(\BP)=0$ only when $\BP=0$ we conclude
that $\BE=\BE'$. We need to make the technical assumption that 
\beq \Gy^{E}(\BE)\ne \Gy^{J}(\BJ),\quad {\rm for~all~nonzero}~\BE\in\CE,\BJ\in\CJ, \eeq{I6.2a}
to ensure $\CE'$ and $\CJ'$ have no nonzero vector in common. A more insightful meaning to the condition \eq{I6.2a} is given in the next section.

Let  $(\CE',\CJ')$ and  $(\CV,\CP_1,\CP_2,\ldots,\CP_n)$ be our new subspace
collection. Given a solution to the equations
\beq \BE\in\CE,~~~\BJ\in\CJ,~~(\BI-\BGP_1)\BJ=\sum_{i=1}^{n}z_i\BGL_i\BE, \eeq{I6.3}
in the original subspace collection, in which $\BGP_1$ is the projection
onto $\CV$, the fields $\BE'=\Gy^{E}(\BE)$ and $\BJ'=\Gy^{J}(\BJ)$ will
be a solution to the equations
\beq \BE'\in\CE',~~~\BJ'\in\CJ',~~(\BI-\BGP_1)\BJ'=\sum_{i=1}^{n}z'_i\BGL_i\BE', \eeq{I6.4}
in the new subspace collection with
\beq z'_i=z_i c_i^{J}/c_i^{E}. \eeq{I6.5}
Since $\BGP_1\BE'=\BGP_1\BE$ and $\BGP_1\BJ'=\BGP_1\BJ$, it follows that the
$\BY$-tensor functions of the two subspace collections are related by \eq{I6.0}
where 
\beq d_i= c_i^{E}/c_i^{J}. \eeq{I6.7}

In particular, if we choose $c_i^{E}=-c_i^{J}$  for all $i$ we obtain $d_i=-1$. Then using the homogeneity of the function $\BY(z_1,z_2, \ldots,z_n)$ we see that
\beq \BY'(z'_1,z'_2, \ldots,z'_n)= \BY(-z_1',-z_2', \ldots,-z_n')=-\BY(z_1',z_2', \ldots,z_n'). \eeq{I6.7a}
So if to another subspace collection, with an associated function $\BY''(z_1,z_2, \ldots,z_n)$,
we add this new subspace collection according to the prescription given in Section \sect{IadditionY}, then
it produces a subspace collection with an associated $Y$-function which is obtained by subtracting $\BY(z_1,z_2, \ldots,z_n)$ from
$\BY''(z_1,z_2, \ldots,z_n)$. In other words, when $c_i^{E}=-c_i^{J}$  for all $i$, the subspace collection with
subspaces $(\CE',\CJ')$ and  $(\CV,\CP_1,\CP_2,\ldots,\CP_n)$ is the additive inverse%
\index{subspace collections!additive inverse}
of the original subspace collection,
having subspaces $(\CE,\CJ)$ and  $(\CV,\CP_1,\CP_2,\ldots,\CP_n)$, where  $(\CE',\CJ')$ and $(\CE,\CJ)$ are linked by \eq{I6.2}. If the technical
condition \eq{I6.2a} is not satisfied it appears that the subspace collection has no arithmetic inverse.

\section{Operations on subspace collections leaving the associated function invariant}
\labsect{subinv}
\setcounter{equation}{0}

Note from \eq{I6.7} that if we choose $c_i^{J}=c_i^{E}$ for all $i$ then the associated function remains 
invariant. More generally, if we are interested in leaving the associated function invariant, we could take
\beq \CE'=\BC\CE,\quad \CJ'=\BC \CJ ,\quad \CV'=\BC\CV ,\quad \CH'=\BC\CH,\quad \CP'_i=\BC\CP_i, \eeq{I6.8}
where $\BC$ is a nonsingular linear operator which maps $\CK$ to itself.
Then the fields $\BE'=\BC\BE$ and $\BJ'=\BC\BJ$ will
be a solution to the equations
\beq \BE'\in\CE',~~~\BJ'\in\CJ',~~(\BI-\BGP'_1)\BJ'=\sum_{i=1}^{n}z_i\BGL'_i\BE', \eeq{I6.10}
where
\beq \BGP_1'=\BC\BGP_1\BC^{-1},\quad\BGL'_i=\BC\BGL_i\BC^{-1} \eeq{I6.11}
are the projections onto $\CV'$ and $\CP'_i$. If $\Bv_1,\Bv_2,\ldots\Bv_m$ is a basis for $\CV$ then setting $\Bv'_i=\BC\Bv_i$ we can 
take $\Bv'_1$,$\Bv'_2$,$\ldots$,$\Bv'_m$ as a basis for $\CV'$. Since multiplying by $\BC$ is a linear operation the coefficients in the expansions
\beq  \BGP'_1\BE'=\sum_{i=1}^m E'_i\Bu'_i,\quad \BGP_1\BE=\sum_{i=1}^m E_i\Bu_i,\quad \BGP'_1\BJ'=\sum_{i=1}^m J'_i\Bu'_i,\quad \BGP_1\BJ=\sum_{i=1}^m J_i\Bu_i
\eeq{I6.12}
can be equated:
\beq E'_i=E_i,\quad J'_i=J_i, \eeq{I6.13}
and as a consequence the same matrix $\BY$ whose coefficients govern the relation
\beq J_i=\sum_{i=1}^kY_{ik}E_k, \eeq{I6.14}
also govern the relation
\beq J'_i=\sum_{i=1}^kY_{ik}E'_k. \eeq{I6.15}

Due to this equivalence it suffices in the preceeding section to limit attention to the transformations \eq{I6.1} having $c_i^{J}=1$ for all $i$:
it is only the ratio  $d_i= c_i^{E}/c_i^{J}$ that has any real significance. Then $\Gy^{J}(\BP)=\BP$, and the technical condition \eq{I6.2a}
is violated only when there are nonzero vectors $\BE\in\CE$ and $\BJ\in\CJ$ such that
\beq \BE=\BE_1+\BE_2,\quad \BJ+\BJ_1+\BJ_2,\quad \BE_1=\BJ_1\in\CV,\quad \BJ_2=\BL\BE_2\in\CH,\quad{\rm with~~}\BL=\sum_{i=1}^n c_i^{E}\BGL_i.
\eeq{I6.15a}
Thus either $\BY(c_1^{E},c_2^{E},\ldots,c_n^{E})$ has an eigenvalue of $-1$, or the $Y$-problem with $z_i=c_i^{E}$ for all $i$ has a nonunique solution
(with a nontrivial solution having $\BE_2\ne 0$ for the homogeneous problem with $\BE_1=0$ and also $\BJ_1=0$, the latter not being needed for 
nonuniqueness but being needed if \eq{I6.15a} holds). If we are looking for the arithmetic inverse we take $c_i^{E}=-1$ for all $i$, and the inverse
exists except when $\BY(-1,-1,\ldots,-1)=-\BY(1,1,\ldots,1)$ has eigenvalue $-1$ or when the $Y$-problem with $z_i=1$ for all $i$ has a nonunique solution
(with the homogeneous problem having a nontrivial solution with both $\BE_1$ and $\BJ_1$ being zero).

There is a similar invariance of matrix functions associated with $Z(n)$ subspace collections under the linear
transformations,
\beq \quad \CU'=\BC\CU, \quad  \CE'=\BC\CE,\quad \CJ'=\BC \CJ , \quad \CP'_i=\BC\CP_i. \eeq{I6.16}

These invariances are quite natural, as they are isomorphic to changing the basis in the vector spaces $\CH$ or $\CK$. Thus, up to these
trivial equivalences, the arithmetic inverse defined in the previous section is unique. 

\section{Multiplicative Inverses of superfunctions}
\labsect{multsubinv}
\setcounter{equation}{0}
To find the multiplicative inverse%
\index{superfunctions!multiplicative inverse}
of a superfunction  $(F^s)'$ we let $\CK''$ be a vector space with the same dimension as $\CK'$, and
we take $\BC$ as a nonsingular map from $\CK'$ to $\CK''$. We then let 
\beqa \CJ''& = & \BC(\CJ'), \quad \CH''=\BC\CH', \nonum
(\CV^I)''=\BC(\CV^O)', \quad (\CV^O)''=\BC(\CV^I)', \quad\CP''_i=\BC\CP'_i\,\,{\rm for}\,\,i=1,2,\ldots j. 
\eeqa{I0.24j}
Introduce the transformation
\beq \psi(\BP)=\BGP_1''\BP-\BGP_2'', \eeq{I0.24ja}
where $\BGP_1''$ is the projection onto $(\CV^I)''\oplus(\CV^O)''$ and $\BGP_2''$ is the projection onto $\CH''$. This
is a special case of the transformations in \eq{I6.1}. We let $\CE''=\Gy(\BC\CE')$.
Note that the output space $(\CV^O)'$ gets mapped to the input space $(\CV^I)''$, and  the input space $(\CV^I)'$ gets mapped to the output space $(\CV^O)''$,
and apart from these switchings we have essentially made an additive inverse in the $Y$-problem. We still require the technical condition mentioned in the last section,
to ensure that this additive inverse exists: the operator $\BY(1,1,\ldots,1)-\BI$ is nonsingular and the $Y$-problem with $z_i=1$ for all $i$ has a unique solution
(or more precisely the homogeneous problem does not have a nontrivial solution with both $\BE_1$ and $\BJ_1$ being zero).

Now suppose we are given a solution to the superfunction problem associated with $(F^s)'$,
\beq \BE'  =  (\BE^I)' +(\BE^O)'+\BE_2'\in\CE',\quad  \BJ' =  (\BJ^I)' +(\BJ^O)'+\BJ_2'\in\CJ',\quad \BJ_2'=\BL'\BE_2',
\eeq{I0.24k}
where
\beq \BL'=\sum_{i=1}^jz_i'\BGL_i',\eeq{I0.24m}
in which $\BGL_i'$ is the projection onto $\CP_i'$, and
\beq  (\BE^I)',\, (\BJ^I)'\in (\CV^I)',\quad  (\BE^O)',\, (\BJ^O)'\in (\CV^O)',\quad \BE_2',\,\BJ_2'\in\CH'. \eeq{I0.24n}
Now take vectors
\beqa \BE'' & = &\Gy(\BC\BE'),\quad \BJ''=-\BC\BJ', \quad \BE_2''=-\BC\BE_2',\quad \BJ_2''=-\BC\BJ_2' \nonum
 (\BE^I)'' & = & \BC(\BE^O)',\quad (\BE^O)'' =\BC(\BE^I)', \quad  (\BJ^I)'' =-\BC(\BJ^O)',\quad (\BJ^O)'' =-\BC(\BJ^I)'.
\eeqa{I0.24o}
These solve the superfunction problem associated with $(F^s)''$,
\beq  \BE''  =  (\BE^I)'' +(\BE^O)''+\BE_2''\in\CE'',\quad  \BJ'' =  (\BJ^I)'' +(\BJ^O)''+\BJ_2''\in\CJ'',\quad \BJ_2''=\BL''\BE_2'',
\eeq{I0.24p}
where
\beq \BL''=\sum_{i=1}^jz_i''\BGL_i'',\quad z_i''=z_i', \eeq{I0.24q}
in which $\BGL_i''$ is the projection onto $\CP_i''$, and
\beq  (\BE^I)'',\, (\BJ^I)''\in (\CV^I)'',\quad  (\BE^O)'',\, (\BJ^O)''\in (\CV^O)'',\quad \BE_2'',\,\BJ_2''\in\CH''. \eeq{I0.24r}
Next let $\BM_1$ denote the restriction of $\BC$ to the subspace $(\CV^O)'$, i.e., that operator mapping  $(\CV^O)'$ to 
 $(\CV^I)''$, such that $\BM_1\BP=\BC\BP$ for all $\BP\in (\CV^O)'$. Then from \eq{I0.24o} we have $(\BE^I)''=\BM_1(\BE^O)'$ and 
$(\BJ^I)'' =-\BM_1(\BJ^O)'$. To see that  $(F^s)''$ is the inverse of the superfunction $(F^s)'$ when
\beq   \BL'=\sum_{i=1}^jz_i'\BGL_i',\quad \BL''=\sum_{i=1}^jz_i'\BGL_i'', \eeq{I0.24s}
we introduce the operator $\BM_2$ which is the restriction of $\BC^{-1}$ to the subspace $(\CV^O)''$, i.e., that operator mapping  $(\CV^O)''$ to 
 $(\CV^I)'$, such that $\BM_2\BP=\BC\BP$ for all $\BP\in (\CV^I)'$. Then upon taking the product of the superfunctions \eq{I0.24o} implies
\beq   \begin{pmatrix} (\BE^O)''  \\ (\BJ^O)'' \end{pmatrix}=\BF\begin{pmatrix} (\BE^I)'  \\ (\BJ^I)' \end{pmatrix}, \eeq{I0.24t}
where 
\beq \BF=\begin{pmatrix} (\BM_2)^{-1} & 0  \\ 0 & -{\BM_2}^{-1} \end{pmatrix} \eeq{I0.24u}
is the multiplicative identity operator.
From this analysis it looks like there are many multiplicative inverses, paramerized by $\BC$, but in fact all are 
equivalent: this follows from the previous section.

\section{Pruning the subspace collections}%
\index{subspace collections!pruning}
\labsect{Iprune}
\setcounter{equation}{0}
If an $m$ terminal resistor network has a cluster of resistors which is not connected to the rest of the network,
and that cluster does not have any terminals, only internal nodes, then clearly we can discard it without affecting the fields
in the rest of the network and its response matrix. The analogous operation on subspace collections is called pruning.

When $\BL$ is close to $z_0\BI$ we can expand the inverses in \eq{I4.12} and \eq{I4.14} to
obtain the series expansions%
\index{series expansions}
\beq \BE  =  \sum_{j=1}^\infty[\BGG_1(\BL-z_0\BI)/z_0]^j\Be, \eeq{I4a.1a} 
\beq
\BZ  =  z_0\BGG_0+\sum_{j=1}^\infty\BGG_0(\BL-z_0\BI)[\BGG_1(\BL-z_0\BI)/z_0]^j\BGG_0.
\eeq{I4a.1}
From these expansions it is evident that is only those fields in $\CH$ that arise from
products of the operators $\BGG_1$, $\BGL_1$, $\BGL_2$, \ldots, $\BGL_n$ applied to fields
in $\CU$ have any role in determining $\BE$ and the associated function
$\BZ(z_1,z_2,\ldots,z_n)$ (also $\Bj$ and $\BJ$): so we may as well prune away any other fields from the
vector space $\CH$. Thus we can redefine $\CH$ as the smallest subspace containing
$\CU$ that is closed under the action of $\BGG_1$, $\BGL_1$, $\BGL_2$, \ldots, $\BGL_n$
and redefine
\beq \CE=\BGG_1\CH,\quad \CJ=\BGG_2\CH,\quad \CP_j=\BGL_j\CH,\quad j=1,2,\ldots,n.
\eeq{I4a.2}
This imposes constraints on the dimensions of these subspaces, as noted in 
Section 29.2 of \citeAPY{Milton:2002:TOC} where the results are given in the case 
where $\CU$ has dimension 1 and where the spaces are orthogonal. Let $p_j$ be the
dimension of $\CP_j$, $j=1,2,\ldots,n$,  and let $m$, $q_1$ and $q_2$ represent the dimensions
of $\CU$, $\CE$ and $\CJ$. The total dimension of the vector space $\CH$ is therefore
\beq h=m+q_1+q_2=p_1+p_2+\ldots +p_n. \eeq{I4a.3}
Now the space
\beq [\BGL_1(\CU\oplus\CE)]\oplus[\BGL_2(\CU\oplus\CE)]\oplus\ldots\oplus[\BGL_n(\CU\oplus\CE)]
\eeq{I4a.4}
certainly contains $\CU$, and is closed under $\BGG_1$ (because it contains $\CE$) and is closed
under $\BGL_j$ for each $j$. It therefore must be $\CH$ and $\BGL_j(\CU\oplus\CE)$ which has at most
dimension $m+q_1$ must be $\CP_j$. Therefore for each $j$ we have the inequality
\beq p_j\leq m+q_1, \eeq{I4a.5}
and by summing these over $j$ we see that
\beq q_2\leq (n-1)(m+q_1). \eeq{I4a.6}
Similarly the subspace
\beq [\BGL_1(\CU\oplus\CJ)]\oplus[\BGL_2(\CU\oplus\CJ)]\oplus\ldots\oplus[\BGL_n(\CU\oplus\CJ)]
\eeq{I4a.7}
can also be identified with $\CH$ and we obtain the inequalities
\beq p_j\leq m+q_2,\quad q_1\leq (n-1)(m+q_2). \eeq{I4a.8}
In the particular case when $n=2$ the constraints \eq{I4a.6} and \eq{I4a.8} imply that the dimensions
of the subspaces $\CE$ and $\CJ$ can differ by at most $m$. Also in the case $n=2$ we have
\beq p_1=(m+q_2-p_2)+q_1=(m+q_1-p_2)+q_2\geq \max\{q_1,q_2\}, \eeq{I4a.8a}
and similarly for $p_2$. 

Likewise we can redefine $\CK$ as the smallest subspace containing
$\CV$ that is closed under the action of $\BGG_1$, $\BGL_1$, $\BGL_2$, \ldots, $\BGL_n$
and redefine
\beq \CE=\BGG_1\CK,\quad \CJ=\BGG_2\CK,\quad \CP_j=\BGL_j\CK,\quad j=1,2,\ldots,n.
\eeq{I4a.9}
Let $v$ be the dimension of $\CV$, $p_j$ be the
dimension of $\CP_j$, $j=1,2,\ldots,n$,  and let $q_1$ and $q_2$ represent the dimensions
of $\CE$ and $\CJ$. The total dimension of the vector space $\CK$ is therefore
\beq h=q_1+q_2=v+p_1+p_2+\ldots +p_n. \eeq{I4a.10}
The space
\beq \CV\oplus[\BGL_1(\CE)]\oplus[\BGL_2(\CE)]\oplus\ldots\oplus[\BGL_n(\CE)]
\eeq{I4a.11}
certainly contains $\CV$, and is closed under $\BGG_1$ (because it contains $\CE$) and is closed
under $\BGL_j$ for each $j$. It therefore must be $\CK$ and $\BGL_j(\CE)$ which has at most
dimension $q_1$ must be $\CP_j$. Thus for each $j$ we have the inequality
\beq p_j\leq q_1, \eeq{I4a.12}
and summing these over $j$ we obtain
\beq q_2\leq v+(n-1)q_1. \eeq{I4a.13}
Similarly since 
\beq \CK=\CV\oplus[\BGL_1(\CJ)]\oplus[\BGL_2(\CJ)]\oplus\ldots\oplus[\BGL_n(\CJ)],
\eeq{I4a.14}
we obtain the inequalities
\beq p_j\leq q_2,\quad q_1\leq v+(n-1)q_2. \eeq{I4a.15}
When  $n=2$ the constraints \eq{I4a.13} and \eq{I4a.15} imply that the dimensions
of the subspaces $\CE$ and $\CJ$ can differ by at most $v$. Also in the case $n=2$ we have
\beq p_1=(q_2-p_2)+q_1-v=(q_1-p_2)+q_2-v\geq \max\{q_1,q_2\}-v, \eeq{I4a.16}
with a similar inequality for $p_2$.

\section{Expressions for the numerator and denominator in the rational function}%
\index{rational functions!numerator and denominator}
\labsect{Inumden}
\setcounter{equation}{0}

Assume that a $Z(n)$ subspace collection, with $m=1$ has been pruned. Let $\Bw_1$,$\Bw_2$,\ldots,$\Bw_{q_1+1}$ be a basis for $\CU\oplus\CE$ with $\Bw_1$ in $\CU$ and $\Bw_2$,$\Bw_3$,\ldots,$\Bw_{q_1+1}$ in $\CE$. In this
basis $(\BGG_0+\BGG_1)\BGL_i(\BGG_0+\BGG_1)$ is represented by a $(q_1+1)\times(q_1+1)$ matrix $\BA_i$, and since the $\BGL_i$ sum up to the identity operator it follows
that
\beq \sum_{i=1}^n\BA_i=\BI. \eeq{I10.9}
Also, because the subspace is pruned, $\BGL_i(\CU\oplus\CE)$ can be identified with $\CP_i$ which implies the matrix $\BA_i$ must have at most rank $p_i$. It is exactly 
$p_i$ if $\CP_i\cap\CJ=0$. The formula \eq{I8.25} for the
$Z$-function implies
\beq 1/Z(z_1,z_2,\ldots,z_n)=\Be_1\cdot[\sum_{i=1}^nz_i\BA_i]^{-1}\Be_1,  \eeq{I10.10}
where $\Be_1$ is the $q_1+1$ component unit vector $[1,0,0,\ldots 0]^T$. Hence, following the argument given in Section 29.2 of \citeAPY{Milton:2002:TOC}, $Z(z_1,z_2,\ldots,z_n)$ can be expressed in the form \eq{I10.1} with numerator
\beq p(z_1,z_2,\ldots,z_n)=\det[\sum_{i=1}^nz_i\BA_i]=\sum_{a_1,a_2,\ldots,a_n}\Ga_{a_1a_2\ldots a_n}z_1^{a_1}z_2^{a_2}\ldots z_{n}^{a_{n}}, \eeq{I10.11}
of degree $1+q_1$, in which the sum extends over all $a_1,a_2,\ldots,a_n$ with
\beq \sum_{i=1}^na_i=1+q_1,\quad 0\leq a_i\leq p_i\quad{\rm for}~i=1,2,\ldots,n. \eeq{I10.12}
Typically one expects that the maximum power of $z_i$ in this polynomial will be the rank of $\BA_i$. However, for example, note that for the matrices
\beq \BM_1= \begin{pmatrix} 0 & 0 & 0 \\ 
                      1 & 1 & 1 \\ 
                      0 & 1 & 1 \\
\end{pmatrix},\quad \BM_2=\BI-\BM_1,
\eeq{I10.12a}
the maximum power of $z_1$ in
\beq  \det[z_1\BM_1+z_2\BM_2]=\det[(z_1-z_2)\BM_1+z_2\BI]=z_2[z_2^2+2z_2(z_1-z_2)] \eeq{I10.12b}
is $1$ while $\BM_1$ has rank 2.

Next let $\Bw_1$,$\Bw_{q_1+2}$,\ldots,$\Bw_{h}$ be a basis for $\CU\oplus\CJ$ with $\Bw_1$ in $\CU$ and $\Bw_{q_1+2}$, $\Bw_{q_1+3}$\ldots,$\Bw_{h}$ in $\CJ$. 
In this basis $(\BGG_0+\BGG_2)\BGL_i(\BGG_0+\BGG_2)$ is represented by a $(q_2+1)\times(q_2+1)$ matrix $\BB_i$, and since the $\BGL_i$ sum up to the identity operator it follows
that
\beq \sum_{i=1}^n\BB_i=\BI. \eeq{I10.13}
Also, because the subspace is pruned, $\BGL_i(\CU\oplus\CJ)$ can be identified with $\CP_i$ which implies the matrix $\BB_i$ must have rank at most $p_i$. It is
exactly $p_i$ if $\CP_i\cap\CE=0$.
The formula \eq{I4.10} for the
$Z$-function implies
\beq Z(z_1,z_2,\ldots,z_n)=\Be_2\cdot[\sum_{i=1}^n\BB_i/z_i]^{-1}\Be_2, \eeq{I10.14} 
where $\Be_2$ is the $q_2+1$ component unit vector $[1,0,0,\ldots 0]^T$. The denominator of this expression, as a polynomial in the variables $1/z_i$, is
\beq \det[\sum_{i=1}^n\BB_i/z_i]=\sum_{b_1,b_2,\ldots,b_n}\Gb_{b_1b_2\ldots b_n}/z_1^{b_1}z_2^{b_2}\ldots z_{n}^{b_{n}}, \eeq{I10.15}
in which the sum extends over all $b_1,b_2,\ldots,b_n$ with
\beq \sum_{i=1}^nb_i=1+q_2,\quad 0\leq b_i\leq p_i\quad{\rm for}~i=1,2,\ldots,n. \eeq{I10.16}
Consequently, for the denominator in the expression \eq{I10.1} for $Z(z_1,z_2,\ldots,z_n)$, we can make the identification
\beq q(z_1,z_2,\ldots,z_n)=\sum_{b_1,b_2,\ldots,b_n}\Gb_{b_1b_2\ldots b_n}z_1^{p_1-b_1}z_2^{p_2-b_2}\ldots z_{n}^{p_n-b_{n}}, \eeq{I10.17}
which is a polynomial of degree $h-(1+q_2)=q_1$. Furthermore the identities \eq{I10.9} and \eq{I10.13} imply the polynomial $p$ and $q$ satisfy the normalization \eq{I10.2a},
i.e.
\beq   \sum_{a_1,a_2,\ldots,a_n}\Ga_{a_1a_2\ldots a_n}=1,\quad \sum_{b_1,b_2,\ldots,b_n}\Gb_{b_1b_2\ldots b_n}=1. \eeq{I10.18}



\section[Rational functions of one variable and $Z(2)$ subspace collections]{The correspondence between rational functions of one variable and $Z(2)$ subspace collections with $m=1$}%
\index{subspace collections!correspondence with rational functions}%
\index{rational functions!correspondence with subspace collections}
\setcounter{equation}{0}

In the case $m=1$ and $n=2$ there are two cases to consider. When the dimension of $h$ is even, $h=2d$, then in order to satisfy the inequalities \eq{I4a.5}, \eq{I4a.6} and \eq{I4a.8}
the subspaces $\CE$ and $\CJ$ must have dimension $d$ and $d-1$ or vice versa and the
subspaces $\CP_1$ and $\CP_2$ must have dimension $d$. 
Without loss of generality, by making a duality transformation if necessary, let us suppose $\CE$ has dimension $d-1$. Given $\Bu\in\CU$ let us take as our basis for $\CH$ the
vectors 
\beq \Bv_{2j-1}=(\BGG_1\BGL_1)^{j-1}\Bu,\quad\Bv_{2j}=(\BGL_1\BGG_1)^{j-1}\BGL_1\Bu,\quad j=1,2,\ldots,d, \eeq{I10.32}
so that
\beq \Bv_1=\Bu,\quad \Bv_{2j}=\BGL_1\BGv_{2j-1},\quad j=1,2,\ldots,d,\quad  \Bv_{2j+1}=\BGG_1\BGv_{2j-1},\quad j=1,2,\ldots,d-1.
\eeq{I10.33}
These fields are independent since if they were not we could prune the subspace collection.%
\index{subspace collections!pruning} 
The vectors $\Bv_{2j+1}, j=1,2,\ldots,d-1$, which number $d-1$, must form a basis for
$\CE$ and so it follows that
\beq \BGG_1\Bv_{2d}=\sum_{i=1}^{d-1}\Gg_i\Bv_{2i+1}. \eeq{I10.34}
Also we have
\beq \BGG_0\Bv_1=\Bv_1,\quad \BGG_0\Bv_{2j}=\Gd_j\Bv_1,\quad j=1,2,\ldots,d,\quad  \BGG_0\Bv_{2j+1}=0,\quad j=1,2,\ldots,d-1. \eeq{I10.35}
The $2d-1$ constants $\Gg_1,\ldots,\Gg_{d-1}$ and $\Gd_1,\ldots,\Gd_d$ characterize the geometry of the subspace collection. The field $\Be+\BE$ must have the expansion
\beq \Be+\BE=\sum_{i=1}^{d}a_i\Bv_{2i-1}, \eeq{I10.36}
and consequently, setting $z_2=1$
\beq \Bj+\BJ=[\BI+(z_1-1)\BGL_1](\Be+\BE)=\sum_{i=1}^{d}a_i\Bv_{2i-1}+(z_1-1)\sum_{i=1}^{d}a_i\Bv_{2i}. \eeq{I10.37}
Since $\BGG_1(\Bj+\BJ)=0$ we arrive at the equations
\beqa 0 & = & \sum_{i=2}^{d}a_i\Bv_{2i-1}+(z_1-1)\sum_{i=1}^{d-1}a_i\Bv_{2i+1}+(z_1-1)\sum_{i=1}^{d-1}a_{d}\Gg_i\Bv_{2i+1} \nonum
 & = & \sum_{i=1}^{d-1}[a_{i+1}+a_i(z_1-1)+\Gg_i a_{d}(z_1-1)]\Bv_{2i+1}.
\eeqa{I10.38}
implying
\beq a_{i+1}+a_i(z_1-1)+\Gg_ia_{d}(z_1-1)=0,\quad i=1,\ldots,d-1. \eeq{I10.39}
Choosing a normalization with $a_d=(1-z_1)^{d-1}$ these equations are solved with
\beq a_i=(1-z_1)^{i-1}-\sum_{j=i}^{d-1}\Gg_{d-1+i-j}(1-z_1)^j. \eeq{I10.40}
Since
\beq \BGG_0(\Be+\BE)=a_1\Bv_1,\quad \BGG_0(\Bj+\BJ)=[a_1+(z_1-1)\sum_{i=1}^d\Gd_i a_i]\Bv_1, \eeq{I10.40a}
we obtain
\beq Z(z_1,1)=1+\frac{(z_1-1)\sum_{i=1}^d\Gd_i a_i}{a_1}. \eeq{I10.41}
Conversely suppose we are given a rational function $Z(z_1,1)$ with a denominator of degree at most $d-1$ and a numerator of degree at most $d$ satisfying
$Z(1,1)=1$. It can be expressed in the 
form 
\beq Z(z_1,1)=\frac{p(z_1,1)}{q(z_1,1)}=1-\frac{\sum_{j=0}^{d-1}t_j(1-z_1)^{j+1}}{1-\sum_{j=1}^{d-1}s_j(1-z_1)^j}. \eeq{I10.41a}
Comparing this with \eq{I10.41} we can make the identifications
\beqa 1-\sum_{j=1}^{d-1}s_j(1-z_1)^j & = & a_1 = 1-\sum_{j=1}^{d-1}\Gg_{d-j}(1-z_1)^j, \nonum
      -\sum_{j=0}^{d-1}t_j(1-z_1)^{j+1} & = & (z_1-1)\sum_{i=1}^d\Gd_i a_i \nonum
   & = & -\sum_{j=0}^{d-1}\Gd_{j+1}(1-z_1)^{j+1}+\sum_{j=0}^{d-1}\sum_{i=1}^j\Gd_{i}\Gg_{d-1+i-j}(1-z_1)^{j+1},
\eeqa{I10.42}
which imply
\beq s_j=\Gg_{d-j},\quad t_0=\Gd_1,\quad t_j=\Gd_{j+1}-\sum_{i=1}^j\Gd_{i}\Gg_{d-1+i-j}\quad j=1,\ldots,d-1.
\eeq{I10.43}
Given the coefficients $s$ and $t$ we can inductively uniquely determine the coefficients $\Gg$ and $\Gd$:
\beq \Gg_j=s_{d-j}, \quad \Gd_1=t_0, \quad \Gd_{j+1}=t_j+\sum_{i=1}^j\Gd_{i}s_{1+j-i}\quad j=1,\ldots,d-1.
\eeq{I10.44}

On the other hand when the dimension of $h$ is odd, $h=2d-1$, then in order to satisfy the inequalities \eq{I4a.5}, \eq{I4a.6} and \eq{I4a.8} 
the subspaces $\CE$ and $\CJ$ must have dimension $d-1$ and the
subspaces $\CP_1$ and $\CP_2$ must have dimension $d-1$ and $d$ or vice versa. Without loss of generality let us suppose $\CP_1$ has dimension $d-1$.  
Given $\Bu\in\CU$ let us take as our basis for $\CH$ the
vectors 
\beq  \Bv_{2j-1}=(\BGG_1\BGL_1)^{j-1}\Bu,\quad j=1,2,\ldots,d-1,\quad\Bv_{2j}=(\BGL_1\BGG_1)^{j-1}\BGL_1\Bu,\quad j=1,2,\ldots,d, \eeq{I10.45}
which satisfy
\beq \Bv_1=\Bu,\quad \Bv_{2j}=\BGL_1\BGv_{2j-1},\quad \Bv_{2j+1}=\BGG_1\BGv_{2j-1},\quad j=1,2,\ldots,d-1. \eeq{I10.46}
Again these fields are independent since if they were not we could prune the subspace collection.%
\index{subspace collections!pruning}
The vectors $\Bv_{2j}, j=1,2,\ldots,d-1$, which number $d-1$, must form a basis for
$\CP_1$ and so it follows that
\beq \BGL_1\Bv_{2d-1}=\sum_{i=1}^{d-1}\Gg_i\Bv_{2i}. \eeq{I10.47}
Also we have
\beq \BGG_0\Bv_1=\Bv_1,\quad \BGG_0\Bv_{2j}=\Gd_j\Bv_1,\quad j=1,2,\ldots,d-1,\quad  \BGG_0\Bv_{2j+1}=0,\quad j=1,2,\ldots,d-1. \eeq{I10.48}
The $2d-2$ constants $\Gg_1,\ldots,\Gg_{d-1}$ and $\Gd_1,\ldots,\Gd_{d-1}$ characterize the geometry of the subspace collection.
The field $\Be+\BE$ has the expansion \eq{I10.36} and so, with $z_2=1$, 
\beq \Bj+\BJ=[\BI+(z_1-1)\BGL_1](\Be+\BE)=\sum_{i=1}^{d}a_i\Bv_{2i-1}+(z_1-1)\sum_{i=1}^{d-1}a_i\Bv_{2i}+(z_1-1)\sum_{i=1}^{d-1}a_d\Gg_i\Bv_{2i}.
\eeq{I10.49}
Since $\BGG_1(\Bj+\BJ)=0$ we arrive at the equations
\beqa 0 & = & \sum_{i=2}^{d}a_i\Bv_{2i-1}+(z_1-1)\sum_{i=1}^{d-1}a_i\Bv_{2i+1}+(z_1-1)\sum_{i=1}^{d-1}a_{d}\Gg_i\Bv_{2i+1} \nonum
 & = & \sum_{i=1}^{d-1}[a_{i+1}+a_i(z_1-1)+\Gg_i a_{d}(z_1-1)]\Bv_{2i+1},
\eeqa{I10.50}
implying \eq{I10.39} which has the solution \eq{I10.40}. Since
\beq \BGG_0(\Be+\BE)=a_1\Bv_1,\quad \BGG_0(\Bj+\BJ)=[a_1+(z_1-1)\sum_{i=1}^{d-1}\Gd_i(a_i+a_d\Gg_i)]\Bv_1=[a_1-\sum_{i=1}^{d-1}\Gd_ia_{i+1}]\Bv_1,
\eeq{I10.51}
we obtain
\beq Z(z_1,1)=1-\frac{\sum_{i=1}^{d-1}\Gd_ia_{i+1}}{a_1}. \eeq{I10.52}
Conversely suppose we are given a rational function $Z(z_1,1)$ with a denominator of degree at most $d-1$ and a numerator of degree at most $d-1$. It can be expressed in the 
form 
\beq Z(z_1,1)=1-\frac{\sum_{j=1}^{d-1}t_j(1-z_1)^{j}}{1-\sum_{j=1}^{d-1}s_j(1-z_1)^j}. \eeq{I10.52a}
Comparing this with \eq{I10.52} we can make the identifications
\beqa 1-\sum_{j=1}^{d-1}s_j(1-z_1)^j & = & a_1 = 1-\sum_{j=1}^{d-1}\Gg_{d-j}(1-z_1)^j, \nonum
      \sum_{j=1}^{d-1}t_j(1-z_1)^{j} & = & \sum_{i=1}^{d-1}\Gd_i a_{i+1} \nonum
   & = & \sum_{j=1}^{d-1}\Gd_j(1-z_1)^{j}-\sum_{j=2}^{d-1}\sum_{i=1}^{j-1}\Gd_{i}\Gg_{d+i-j}(1-z_1)^j,
\eeqa{I10.53}
which imply
\beq s_j=\Gg_{d-j},\quad j=1,\ldots,d-1,\quad t_1=\Gd_1,\quad t_j=\Gd_{j}-\sum_{i=1}^{j-1}\Gd_{i}\Gg_{d+i-j}\quad j=2,\ldots,d-1.
\eeq{I10.54}
 Given the coefficients $s$ and $t$ we can inductively uniquely determine the coefficients $\Gg$ and $\Gd$:
\beq \Gg_j=s_{d-j},\quad j=1,\ldots,d-1 \quad \Gd_1=t_1, \quad \Gd_{j}=t_j+\sum_{i=1}^j\Gd_{i}s_{j-i}\quad j=2,\ldots,d-1.
\eeq{I10.55}

One can see from this analysis that there can be more than one pruned subspace collection associated with a rational function $Z(z_1,1)$.
It may happen that one pruned $Z(n)$ subspace collection gives rise to polynomials $p(z_1,1)=f(z_1,1)r(z_1,1)$ and $q(z_1,1)=g(z_1,1)r(z_1,1)$ 
while another pruned Z(n) subspace collection gives rise to polynomials $p'(z_1,1)=f(z_1,1)r'(z_1,1)$ and $q'(z_1,1)=t(z_1,1)r'(z_1,1)$, so that
both give rise to the same function $Z(z_1,1)$. However there is a one-to-one correspondence if the pruned subspace collection is such
that the polynomials $p(z_1,z_2)$ and $q(z_1,z_2)$ have no factor in common, and this correspondence is given by the above algorithm

\section[Rational functions of two variables and $Z(3)$ subspace collections]{On the correspondence between certain rational functions of two variables and $Z(3)$ subspace collections with $m=1$}%
\index{subspace collections!correspondence with rational functions}%
\index{rational functions!correspondence with subspace collections}
\labsect{Icount}
\setcounter{equation}{0}

In the case $m=1$ and $n=3$ can we uniquely recover a generic subspace collection (modulo the linear transformations \eq{I6.16}) from 
knowledge of the rational function $Z(z_1,z_2,1)$? The answer is no, but let us first provide a counting argument which suggests that, at least in the
generic case, we can recover the subspace collection up to a finite number of possibilities. The counting argument%
\index{counting argument}
is similar to that given in
Section 29.2 of \citeAPY{Milton:2002:TOC} but here we do not assume that the subspaces are orthogonal.

How many independent coefficients $\Ga_{a_1a_2a_3}$ are there in a polynomial
\beq p(z_1,z_2,1)=\sum_{a_1,a_2,a_3}\Ga_{a_1a_2a_3}z_1^{a_1}z_2^{a_2}, \eeq{I10.19}
that satisfies
\beq a_1+a_2+a_3=1+q_1, \quad 0\leq a_i\leq p_i\leq 1+q_1,\quad i=1,2,3\,? \eeq{I10.20}
Without loss of generality, following Section 29.2 of \citeAPY{Milton:2002:TOC}, let
us suppose that $p_1\geq p_2 \geq p_3$. With $a_1$ fixed in the regime
$0\leq a_1 < 1+q_1-p_2$, the constant $a_2$ can take integer values
from $a_2= 1+q_1-a_1-p_3$ (where $a_3=p_3$) to $a_2=p_2$, that is,  a total of
$p_2+p_3+a_1-q_1$ different values. With $a_1$ fixed in the regime
$1+q_1-p_2 \leq a_1 < 1+q_1-p_3$, the constant $a_2$ can take integer values
from $a_2= 1+q_1-a_1-p_3$ (where $a_3=p_3$) to $a_2=1+q_1-a_1$ (where $a_3=0$) that is,  a total of
$p_3+1$ different values. Finally, with $a_1$ fixed in the regime
$1+q_1-p_3 \leq a_1 \leq p_1$, the constant $a_2$ can take integer values
from $a_2= 0$ to $a_2=1+q_1-a_1$ (where $a_3=0$), that is,  a total of
$2+q_1-a_1$ different values. Therefore the total number of coefficients
in the polynomial is
\beqa &~&\sum_{a_1=0}^{q_1-p_2}(p_2+p_3+a_1-q_1)
   +\sum_{a_1=1+q_1-p_2}^{q_1-p_3}(p_3+1)
   +\sum_{a_1=1+q_1-p_3}^{p_1}(2+q_1-a_1) \nonum &~& \quad =(q_1-p_2+1)(p_2+p_3-q_1)+\frac{1}{2}(q_1-p_2+1)(q_1-p_2)
+(p_2-p_3)(p_3+1) \nonum &~& \quad\quad+ (p_1+p_3-q_1)(2+p_3)-\frac{1}{2}((p_1+p_3-q_1)(p_1+p_3-q_1+1)\nonum &~& \quad
=k_1+1, \eeqa{I10.21}
where
\beq  k_1=[2(1+q_1)q_2-p_1^2-p_2^2-p_3^2+h]/2, \eeq{I10.22}
in which $h=p_1+p_1+p_3$ and $q_2=h-1-q_1$. These coefficients are not all independent since, from \eq{I10.18} the
$\Ga_{a_1a_2a_3}$ must sum to one. Subtracting this constraint gives $k_1$ independent coefficients.

Similarly in a polynomial 
\beq q(z_1,z_2,1)=\sum_{b_1,b_2,b_3}\Gb_{b_1b_2b_3}z_1^{p_1-b_1}z_2^{p_2-b_2}, \eeq{I10.23}
that satisfies
\beq b_1+b_2+b_3=1+q_2, \quad 0\leq a_i\leq p_i\leq 1+q_2,\quad i=1,2,3,\quad \sum_{b_1,b_2,b_3}\Gb_{b_1b_2b_3}=1, \eeq{I10.24}
there are a total of 
\beq  k_2=[2(1+q_2)q_1-p_1^2-p_2^2-p_3^2+h]/2 \eeq{I10.25}
independent coefficients. Hence the total number of independent coefficients in the rational function
\beq Z(z_1,z_2,1)=\frac{p(z_1,z_2,1)}{q(z_1,z_2,1)} \eeq{I10.26}
is
\beq k_1+k_2=(1+q_1)q_2+(1+q_2)q_1-p_1^2-p_2^2-p_3^2+h=h^2-p_1^2-p_2^2-p_3^2-q_1^2-q_2^2. \eeq{I10.27}

Now how many parameters describe a $Z(n)$ subspace collection, when the spaces $\CU$, $\CE$, $\CJ$, $\CP_1$, $\CP_2$, and $\CP_3$ have dimensions
1, $q_1$, $q_2$, $p_1$, $p_2$, and $p_3$, with $1+q_1+q_2=p_1+p_2+p_3=h$? Let $\Bw_1$,$\Bw_2$,\ldots,$\Bw_{h}$ be a basis for $\CH$ with $\Bw_1$ in $\CU$, $\Bw_2$,$\Bw_3$,\ldots,$\Bw_{q_1+1}$ in $\CE$, and $\Bw_{q_1+2}$, $\Bw_{q_1+3}$\ldots,$\Bw_{h}$ in $\CJ$. Recall that it requires $s(d-s)$ parameters to describe the orientation of a subspace of dimension
$s$ in a space of dimension $d$. Therefore, it requires 
\beq p_1(h-p_1)+(h-p_2)p_2+(h-p_3)p_3=h^2-p_1^2-p_2^2-p_3^2 \eeq{I10.28}
parameters to describe the orientation of the subspaces $\CP_1$, $\CP_2$, and $\CP_3$ with respect to this basis. However some of these subspace collections are equivalent, linked
through transformations of the form \eq{I6.16}. If respect to this basis $\BC$ is represented by a matrix with block form
\beq  \BC=\begin{pmatrix} c & 0 & 0 \\ 
                      0 & \BC_1 & 0 \\ 
                      0 & 0 & \BC_2 \\
\end{pmatrix}, 
\eeq{I10.29}
where $c$ is a scalar, while $\BC_1$ and $\BC_2$ are $q_1\times q_1$ and $q_2\times q_2$ matrices,
then it will leave the subspaces $\CU$, $\CE$ and $\CJ$ unchanged. The transformation $\BC=a\BI$ leaves all subspaces unchanged for any scalar $a\ne 0$, 
and so to factor out such trivial transformations we should choose $c=1$. The number of remaining independent parameters in $\BC$ is then $q_1^2+q_2^2$. 
Subtracting these from \eq{I10.28} we see that the number of parameters describing the $Z(n)$ subspace collection is
\beq h^2-p_1^2-p_2^2-p_3^2-q_1^2-q_2^2=k_1+k_2. \eeq{I10.30}
The precise agreement between the number of coefficients in the rational function and the number of parameters describing the $Z(n)$ subspace collection is
curious (since it holds for all  $q_1$, $q_2$, $p_1$, $p_2$, and $p_3$, with $1+q_1+q_2=p_1+p_2+p_3=h$). Despite this coincidence we now show that it is
not possible to uniquely recover 
a generic subspace collection (modulo the linear transformations \eq{I6.16}) from knowledge of the associated rational function $Z(z_1,z_2,1)$.

Let us consider a subspace collection with $h=5, q_1=q_2=2, p_1=p_2=1, p_3=3$ giving $k_1+k_2=6$ according to the formula \eq{I10.27}. Given $\Bu\in\CU$ we choose as our basis
the vectors
\beq \Bv_0=\Bu,\quad\Bv_1=\BGL_1\Bu,\quad\Bv_2=\BGL_2\Bu,\quad\Bv_3=\BGG_1\BGL_1\Bu,\quad\Bv_4=\BGG_1\BGL_2\Bu,  \eeq{I10.30a}
with the closure relations
\beqa \BGL_1\Bv_3 & = & \Gg_1\Bv_1,\quad\BGL_2\Bv_3=\Gg_2\Bv_2,\quad \BGL_1\Bv_4=\Gg_3\Bv_1,\quad\BGL_2\Bv_4=\Gg_4\Bv_1,\nonum
 \BGG_0\Bv_1 & = & \Gd_1\Bv_0,\quad\BGG_0\Bv_2=\Gd_2\Bv_0,
\eeqa{I10.30b}
expressed in terms of the $6$ parameters $\Gg_1,\Gg_2,\Gg_3,\Gg_4,\Gd_1$, and $\Gd_2$ which describe the subspace collection. The question is: can one uniquely recover
these six parameters from $Z(z_1,z_2,1)$? Although the following analysis extends easily to the case  of arbitrary $\Gg_1$ and $\Gg_4$ let us assume, for simplicity,
that $\Gg_1=\Gg_4=0$ and ask whether one can recover the remaining four parameters. The field $\Be+\BE$ must have the expansion
\beq \Be+\BE=a_0\Bv_0+a_1\Bv_3+a_2\Bv_4, \eeq{I10.30c}
and consequently, setting $z_3=1$,
\beqa \Bj+\BJ & = & [\BI+(z_1-1)\BGL_1+(z_2-1)\BGL_2](\Be+\BE) \nonum
& = & a_0\Bv_0+a_1\Bv_3+a_2\Bv_4+(z_1-1)(a_0+a_2\Gg_3)\Bv_1+(z_2-1)(a_0+a_1\Gg_2)\Bv_2. 
\eeqa{I10.30d}
Since $\BGG_1(\Bj+\BJ)=0$ we arrive at the equations
\beq 0 = a_1\Bv_3+a_2\Bv_4+(z_1-1)(a_0+a_2\Gg_3)\Bv_3+(z_2-1)(a_0+a_1\Gg_2)\Bv_4,
\eeq{I10.30e}
implying
\beq a_1+(z_1-1)(a_0+a_2\Gg_3)=0,\quad a_2+(z_2-1)(a_0+a_1\Gg_2)=0. \eeq{I10.30f}
These equations have as a solution,
\beqa a_0 & = & 1-(z_1-1)(z_2-1)\Gg_2\Gg_3,\nonum
      a_1 & = & \Gg_3(z_1-1)(z_2-1)-(z_1-1),\nonum
      a_2 & = & \Gg_2(z_1-1)(z_2-1)-(z_2-1). \eeqa{I10.30g}
Since
\beq \BGG_0(\Be+\BE)=a_0\Bv_0,\quad \BGG_0(\Bj+\BJ)=[a_0+(z_1-1)(a_0+a_2\Gg_3)\Gd_1+(z_2-1)(a_0+a_1\Gg_2)\Gd_2]\Bv_0, \eeq{I10.30h}
we obtain
\beqa Z(z_1,z_2,1) & = & 1+\frac{(z_1-1)(a_0+a_2\Gg_3)\Gd_1+(z_2-1)(a_0+a_1\Gg_2)\Gd_2}{a_0}\nonum 
&=& 1+\frac{\Gd_1(z_1-1)-\Gg_3\Gd_1(z_1-1)(z_2-1)+\Gd_2(z_2-1)-\Gg_2\Gd_2(z_1-1)(z_2-1)}
{1-(z_1-1)(z_2-1)\Gg_2\Gg_3}. \nonum &~&
\eeqa{I10.30i}
Given this function we can uniquely determine $\Gd_1$ and $\Gd_2$ from the coefficients of $(z_1-1)$ and $(z_2-1)$ in the numerator. Also from the
coefficients of $(z_1-1)(z_2-1)$ in the numerator and denominator we can uniquely determine 
\beq t_1=\Gg_2\Gg_3,\quad t_2=\Gg_3\Gd_1+\Gg_2\Gd_2, \eeq{I10.30j}
in terms of which there are two possible values of $\Gg_2$, namely
\beq \Gg_2=\frac{t_3\pm\sqrt{t_3^2-3t_1\Gd_1\Gd_2}}{2\Gd_1}. \eeq{I10.30k}
Thus we cannot uniquely recover the subspace collection parameters from $Z(z_1,z_2,1)$. 

It remains an open question, raised at the end of Section 29.2 of \citeAPY{Milton:2002:TOC},
as to whether in general one can uniquely recover the subspace collection parameters when, with respect to some inner product, 
the subspaces $\CU$, $\CE$ and $\CJ$ are mutually orthogonal, and the subspaces $\CP_1$, $\CP_2$ and $\CP_3$ are mutually orthogonal. 
These orthogonality constraints overdetermine the system of equations needed to recover the subspace collection parameters which provides 
some hope that we can recover them.  It 
would be useful if one could uniquely recover the subspace collection parameters (the weight and normalization 
matrices introduced in Milton, \citeyearNP{Milton:1987:MCEa}, \citeyearNP{Milton:1987:MCEb}) from say the effective 
conductivity $\Gs_*(\Gs_1,\Gs_2,\Gs_3)$ of an isotropic composite of three isotropic phases having conductivities $\Gs_1$, $\Gs_2$, and $\Gs_3$
as then one could obtain the effective response tensor for coupled field problems. We will see in Chapter 9 of this book (\citeAY{Milton:2016:ETC})
that the effective response tensor just 
depends on the weight and normalization matrices for the uncoupled conductivity problem.


\section[Visualizing the poles and zeros of functions $Z(z_1,z_2,z_3)$]
{Visualizing the poles and zeros of functions associated with orthogonal $Z(3)$ subspace collections when $m=1$}
\setcounter{equation}{0}
For scalar functions $Z(z_1,z_2,z_3)$, associated with orthogonal $Z(3)$ subspace collections, satisfying the homogeneity,
Herglotz, and normalization properties, the trajectories of their poles and zeros in $(z_1,z_2,z_3)$ space, with $z_1$, $z_2$, and $z_3$ taking real values, have a beautiful visualization as trajectories on three interlinked hexagons:%
\index{pole-zero trajectories!hexagon representation}
To obtain this visualization we follow Appendix C in \citeAPY{Nicorovici:1993:TPT}: see also figure 5 in that paper. 

First note that if we set $z_3=1$, then the poles and zeros of  $Z(z_1,z_2,1)$ lie in one of the three quadrants: 
\begin{itemize}
\item The quadrant $z_1\leq 0$, $z_2\geq 0$;
\item The quadrant $z_2\leq 0$, $z_1\geq 0$;
\item The quadrant $z_1\leq 0$, $z_2\leq 0$.
\end{itemize}
Of course we can visualize the pole and zero trajectories by plotting them in this plane, but this has the disadvantage that the three variables
$z_1$, $z_2$ and $z_3$ are not treated in a symmetric way, and the disadvantage that its hard to see what is happening when $z_1$ and/or $z_2$ is large,
and it is hard to see  what is happening near the origin $z_1=z_2=0$ since the trajectories can bunch up there.  To get around this we map each
of the three quadrants to a hexagon. Given a quadrant, the point $z_1=z_2=0$ gets blown up to form one edge of the hexagon; 
the two edges of the quadrant where $z_1$ or $z_2$ is zero, but not the other, get mapped to
two other edges of the hexagon; the two ``boundaries'' of the quadrant where $|z_1|$ or $|z_2|$ is infinite but other is finite get mapped to
two more edges of the hexagon; finally $z_1=z_2=\infty$ gets mapped to the final sixth edge of the hexagon. We remark that
just as a pole trajectory can cross from one quadrant to another, so too can it jump from the boundary of one hexagon to the corresponding
point on the boundary of another hexagon.

To be more precise, we introduce the three variables
\beq t_1=\frac{1}{1+|z_2/z_3|},\quad t_2=\frac{1}{1+|z_3/z_1|},\quad t_3=\frac{1}{1+|z_1/z_2|}. \eeq{I10.30k-0}
Clearly $(t_1,t_2,t_3)$ takes values in the unit cube. It is confined to a surface within the unit cube as the three
ratios $|z_2|/|z_3|$, $|z_3|/|z_1|$ and $|z_1|/|z_2|$ are not independent, but have product $1$. The next step is
to map these three variables onto three variables $s_1$, $s_2$ and $s_3$ lying in the plane $s_1+s_2+s_3=0$ using the projection
\beq s_1=2t_1-t_2-t_3,\quad s_2=2t_2-t_3-t_1 \quad s_3=2t_3-t_1-t_2. \eeq{I10.30k-1}
Finally, we map these down to the $x$--$y$ plane:
\beq x=s_1,\quad y=(s_1+2s_2)/\sqrt{3}. \eeq{I10.30k-2}
\begin{sloppypar}
Some normalization is needed, so in the hexagon where $z_1$ is negative and $z_2$ and $z_3$ are positive, we plot $Z(z_1,z_2,z_3)/\sqrt{z_2z_3}$;
in the hexagon where $z_2$ is negative and $z_1$ and $z_3$ are positive, we plot $Z(z_1,z_2,z_3)/\sqrt{z_1z_3}$; and in the hexagon where $z_3$ is negative 
and $z_1$ and $z_2$ are positive, we plot $Z(z_1,z_2,z_3)/\sqrt{z_1z_2}$.
\end{sloppypar}

\begin{figure}[ht]
\centering
\includegraphics[width=0.75\textwidth]{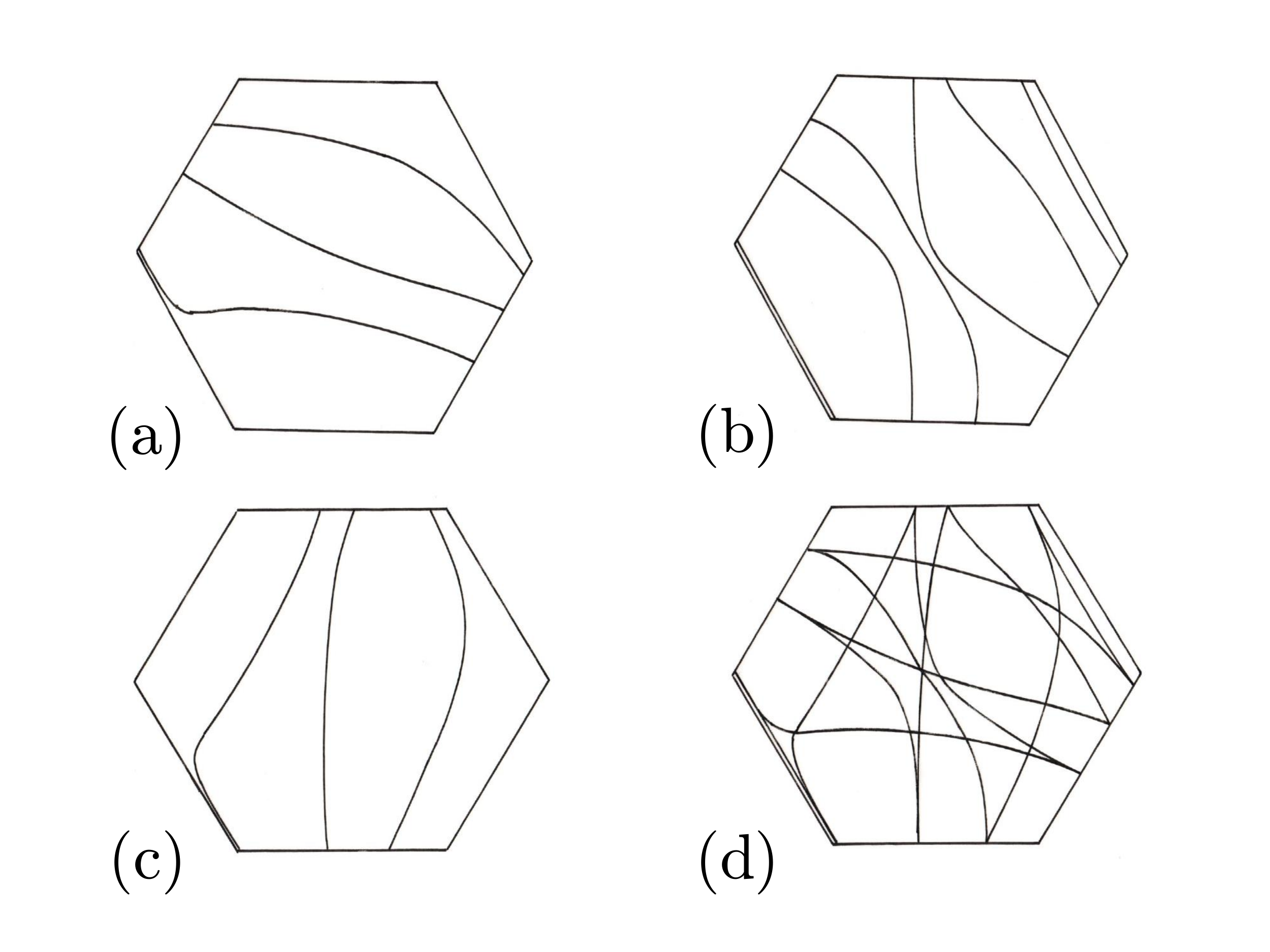}
\caption[The pole trajectory of the function $Z(z_1,z_2,z_3)$]
        {The pole trajectory of the function $Z(z_1,z_2,z_3)$ as visualized using the representation using three interlinked hexagons. 
The hexagon in (a) corresponds to real values of $(z_1,z_2,z_3)$ where $z_2$ and $z_3$ have the same sign, but $z_1$ has the opposite sign.
The hexagon in (b) corresponds to real values of $(z_1,z_2,z_3)$ where $z_1$ and $z_3$ have the same sign, but $z_2$ has the opposite sign.
The hexagon in (c) corresponds to real values of $(z_1,z_2,z_3)$ where $z_1$ and $z_2$ have the same sign, but $z_3$ the opposite sign.
By superimposing all three pictures one obtains (d) where the pole trajectory is like that of a billard ball bouncing around a hexagonal table,
following curved paths. The zero trajectory is similiar, but for clarity we chose not to include it. Note that the dimensions $3$, $6$ and $3$ of the
subspaces $\CP_1$, $\CP_2$, and $\CP_3$ can be immediately read off from the number of paths crossing the hexagons in (a), (b) and (c). }
\label{IS-10}
\end{figure}

Figure~\ref{IS-10} uses this approach to visualize the pole trajectory of a function $Z(z_1,z_2,z_3)$ associated with a $Z(3)$-subspace collection
\beq \CH=\CU\oplus\CE\oplus\CJ=\CP_1\oplus\CP_2\oplus\CP_3, \eeq{I10.30k-3} 
where in this example $\CH$ is 12-dimensional; $\CU$ is one-dimensional; $\CP_1$ is 3-dimensional;  $\CP_2$ is 6-dimensional; $\CP_3$ is 3-dimensional. 
Note that as
the subspace collection does not need pruning, the dimensions of 
$\CP_1$, $\CP_2$, and $\CP_3$ can be immediately read off from the figure by simply counting the number of pole paths on each hexagon:
figures (a), (b), and (c) have $3$, $6$ and $3$ pole paths corresponding to the dimensions of $\CP_1$, $\CP_2$, and $\CP_3$, respectively. To understand 
this, first recognize that when $z_2$ and $z_3$ are fixed, and real and positive, $Z(z_1,z_2,z_3)$ is a Herglotz function of $z_1$ taking real positive
values when $z_1>0$. Thus all its poles must be simple and located on the negative real $z_1$-axis, i.e. on  the hexagon (a). Also because the 
subspace is pruned $\BGL_1(\CU\oplus\CJ)$ can be identified with $\CP_1$ (Section \sect{Inumden}), and hence the matrix $\BC_1$
representing $\BGL_1(\BGG_0+\BGG_2)$ has rank $p_1$. Then as $\BC_1$ and $\BC_1^T\BC_1$ have equal rank (this well-known fact can easily be seen by showing
that they have the same null-space), and as the subspace collection is orthogonal, it follows that the matrix $\BB_1=\BC_1^T\BC_1$ representing
$(\BGG_0+\BGG_2)\BGL_1(\BGG_0+\BGG_2)$ has exactly rank $p_1$. Similarly the matrix $\BA_1$ representing
$(\BGG_0+\BGG_2)\BGL_1(\BGG_0+\BGG_2)$ has exactly rank $p_1$. Therefore the sum over $a_1$ in the numerator in \eq{I10.11}, goes up to $a_1=p_1$, while the sum in the
denominator in \eq{I10.17}, goes from $0$ up to $b_1=p_1$ or (when all the coefficients $\Gb_{p_1b_2b_3}$ are zero) to $b_1=p_1-1$: it cannot go only up to
$b_1=p_1-2$, since as a function of $z_1$, $Z(z_1,z_2,z_3)/\sqrt{z_2z_3}$ with fixed $z_2>0$ and fixed $z_3>0$ can only have a simple pole at $z_1=\infty$. 
When the sum over $b_1$ goes up to $b_1=p_1$, there are clearly $p$ poles of the function $Z(z_1,z_2,z_3)/\sqrt{z_2z_3}$ on the hexagon as $z_1$ varies
with fixed $z_2>0$ and fixed $z_3>0$. When the sum over $b_1$ goes up to $b_1=p_1-1$, there are still $p$ poles of the function $Z(z_1,z_2,z_3)/\sqrt{z_2z_3}$ on the hexagon as $z_1$ varies
with fixed $z_2>0$ and fixed $z_3>0$ provided we count the pole at $z_1=\infty$.

The dimension $q_2$ of the subspace $\CJ$ can also generically be read off from the pole trajectories on the three hexagons. Consider the
edge joining two of the hexagons that corresponds to the values $z_2=0$, and $z_3=1$ with $z_1<0$ varying. Then the only coefficients 
$\Gb_{b_1b_2b_3}$ that can contribute to the denominator in \eq{I10.17} are those with $p_2=b_2$. The first constraint in \eq{I10.16} then implies
\beq b_1+b_3=1+q_2-p_2. \eeq{I10.30k-4}
So $b_1$ can only range from $0$ up to the maximum of $p_1$ and $1+q_2-p_2=p_1+p_3-q_1$. Note that according to the inequality \eq{I4a.5},
$q_1\geq p_3-1$ so $1+q_2-p_2$ could be as large as $p_1+1$. If there are less than $p_1$ pole trajectories crossing this edge joining the
hexagons, the number of these crossing pole trajectories should generically allow us to determine $q_2$ and hence $q_1$, assuming $p_1$, $p_2$ and
$p_3$ have been determined from the number of pole trajectories on each hexagon. If there are exactly $p_1$ pole trajectories crossing the edge
then $q_2$ could be $p_3$ or $p_3+1$. To determine which it is (or as an additional check on the value of $q_2$) 
we could look at pole trajectories, or zero trajectories, crossing other edges where the hexagons meet.

This visualization may be useful in finding other topological features of the trajectories, which hopefully could be connected with topological features of the subspace collections.

\section{Normalization operations on subspace collections}
\index{subspace collections!normalization operations}
\setcounter{equation}{0}

Rational functions of a single variable may be expanded in continued fractions,%
\index{continued fraction expansions}
which incorporate successively
higher and higher order terms in the series expansion of the function about a point. The analogous procedure with
subspace collections is achieved through normalization and reduction operations,%
\index{subspace collections!reduction operations}
subject to some technical assumptions. The associated
functions are then linked, and provided the technical assumptions hold at each level, these links provide continued fractions for multivariate functions
$\BZ(z_1,z_2,\ldots,z_n)$ and $\BY(z_1,z_2,\ldots,z_n)$ incorporating matrices of increasingly high dimension at each level in the continued fraction.

The normalization and reduction operations are discussed in this and the next section. For more insight, in the case where the subspaces in the direct sums
are orthogonal (see Milton \citeyearNP{Milton:1987:MCEa}, \citeyearNP{Milton:1987:MCEb} and Sections 19.2, 20.6 and 29.5 in Milton \citeyearNP{Milton:2002:TOC}).

Normalization reverses extension.%
\index{subspace collections!extension operations}
Given a subspace collection
\beq \CK=\CE'\oplus\CJ'=\CV\oplus\CP_1\oplus\CP_2\oplus\cdots\oplus\CP_n, \eeq{I5.10}
define
\beqa \CH & = & \CP_1\oplus\CP_2\oplus\cdots\oplus\CP_n,\quad\CE=\CE'\cap\CH,\quad\CJ=\CJ'\cap\CH,\nonum
\quad\CU & = & \BGP_2\BGG'_1\CV=\BGP_2(\BI-\BGG'_2)\CV=\BGP_2\BGG'_2\CV, \quad
\quad \widetilde{\CE}=\BGG_1'\CV,\quad \widetilde{\CJ}=\BGG_2'\CV, 
\eeqa{I5.12}
where $\BGG_1'$ and $\BGG_2'$ are the projections onto $\CE'$ and $\CJ'$, and $\BGP_2$ is the projection onto $\CH$.

We assume that the $Y$-problem has a unique solution when $\BL=\BI$ for $\BJ_1\in\CV$ given $\BE_1\in\CV$. In other
words, we assume that the equations 
\beqa \BE_1+\BE_2 & \in & \CE',\quad \BJ_1+\BJ_2\in \CJ',\quad \BJ_2=\BE_2, \quad \BE_1,\BJ_1\in \CV,\quad \BE_2,\BJ_2\in \CH,  \nonum
\BE_1+\underline{\BE}_2 & \in & \CE',\quad \underline{\BJ}_1+\underline{\BJ}_2\in \CJ',\quad \underline{\BJ}_2=\underline{\BE}_2, \quad
\underline{\BJ}_1\in \CV,\quad \underline{\BE}_2,\underline{\BJ}_2\in \CH,
\eeqa{I5.12a}
imply $\underline{\BJ}_1=\BJ_1$. Subtracting these equations we see that
\beq \BE\equiv\BE_2-\underline{\BE_2}\in\CE',\quad \BJ\equiv \BJ_1+\BJ_2-\underline{\BJ}_1-\underline{\BJ}_2\in \CJ', 
\quad \BJ_2-\underline{\BJ_2}=\BE. \eeq{I5.12b}
These imply
\beq \BE\in\CH,\quad \BE=\BJ-\Bv,\quad{\rm where}\quad\Bv=\BJ_1-\underline{\BJ}_1. \eeq{I5.12c}
The uniqueness assumption means that these equations imply $\Bv=0$ (and if $\Bv=0$ then necessarily
$\BE=\BJ=0$ since $\CE'$ and $\CJ'$ have no vector in common). 
The relation $\BE=\BJ-\Bv$ 
with $\BE\in\CE'\cap\CH$ implies
\beq \BE=-\BGG'_1\Bv,  \eeq{I5.12d}
which will only have the trivial solution $\Bv=0$ if and only if
\beq \CH\cap\widetilde{\CE}=0\quad{\rm and}\quad \CV\cap\CJ'=0, \eeq{I5.12e}
where the latter guarantees that $\BGG'_1\Bv=0$ implies $\Bv=0$.

We also assume that the $Y$-problem has a unique solution when $\BL=\BI$ for $\BE_1\in\CV$ given $\BJ_1\in\CV$. By similar
analysis this is satisfied if and only if  
\beq \CH\cap\widetilde{\CJ}=0\quad{\rm and}\quad\CV\cap\CE'=0. \eeq{I5.12f}

We now establish that
\beq \CW\equiv\widetilde{\CE}\oplus\widetilde{\CJ}=\CV\oplus\CU. \eeq{I5.13}
First note that $\CV$ and $\CU$ have no vector in common since $\CU\subset\CH$, and similarly
$\widetilde{\CE}$ and $\widetilde{\CJ}$ have no vector in common since $\CE'\cap\CJ'=0$.
Clearly $\CW$ contains $\CV$. To show it contains $\CU$ notice that 
\beq \CU=\BGP_2\BGG'_1\CV=(\BI-\BGP_1)\BGG'_1\CV\subset\BGG'_1\CV\oplus\BGP_1\BGG'_1\CV\subset\widetilde{\CE}\oplus\CV\subset\CW.
\eeq{I5.14}
Together these imply $\CV\oplus\CU\subset\CW$. Finally we have
\beq \widetilde{\CE}=\BGG_1'\CV=(\BGP_1+\BGP_2)\BGG_1'\CV\subset\BGP_1\BGG_1'\CV\oplus\BGP_2\BGG_1'\CV\subset\CV\oplus\CU, \eeq{I5.15}
and similarly $\widetilde{\CJ}\subset\CV\oplus\CU$. Together these imply $\CW\subset\CV\oplus\CU$, establishing \eq{I5.13}.

If $\CV$ has dimension $m$ then $\widetilde{\CE}$ must also have dimension $m$ since otherwise $\BGG_1'\Bv=0$ for some nonzero $\Bv\in\CV$,
implying $\Bv=\BGG_2'\Bv$ which only has the solution $\Bv=0$ since $\CV\cap\CJ'=0$. Similarly $\widetilde{\CJ}$ must have dimension
$m$ and \eq{I5.13} then implies $\CU$ must have dimension $m$.
The first condition in \eq{I5.12e} implies
\beq \CW=\CU\oplus\widetilde{\CE}, \eeq{I5.17}
since $\CU\subset\CH$ and $\widetilde{\CE}$ have no vector in common and are $m$-dimensional spaces contained in the $2m$-dimensional space $\CW$. Now 
any vector $\BE'\in\CE'$ has the unique decomposition
\beq \BE'=\BE'_1+\BP,\quad \BE'_1\in\CV,\quad \BP\in\CH, \eeq{I5.18}
and according to \eq{I5.17} $\BE'_1$ has the unique decomposition
\beq \BE'_1=-\Be+\widetilde{\BE}, \quad \Be\in\CU,\quad \widetilde{\BE}\in\widetilde{\CE}. \eeq{I5.19}
So we have the decomposition
\beq \BE'=\widetilde{\BE}+\BE, \eeq{I5.20}
where
\beq \BE=\BP-\Be=\BE'-\widetilde{\BE}\in\CE'\cap\CH=\CE. \eeq{I5.21}
Also the first condition in \eq{I5.12e} implies $\widetilde{\CE}$ and $\CE\subset\CH$ have no vector in common, so the decomposition is unique. Therefore we conclude that
\beq \CE'=\widetilde{\CE}\oplus\CE, \eeq{I5.22}
and similarly the first condition in \eq{I5.12f} implies
\beq \CJ'=\widetilde{\CJ}\oplus\CJ. \eeq{I5.23}
These and \eq{I5.13} imply
\beq \CK=\CV\oplus\CH=\widetilde{\CE}\oplus\CE\oplus\widetilde{\CJ}\oplus\CJ=\CV\oplus\CU\oplus\CE\oplus\CJ,
\eeq{I5.24}
and since $\CU$, $\CE$ and $\CJ$ are all contained in $\CH$ we conclude that
\beq \CH=\CU\oplus\CE\oplus\CJ=\CP_1\oplus\CP_2\oplus\cdots\oplus\CP_n. \eeq{I5.25}

Now a given $\BE'_1\in\CV$ has the unique decomposition \eq{I5.19}. This defines the nonsingular operator $\BK:\CV\to\CU$ such that $\Be=\BK\BE'_1$. (It is nonsingular because $\CV$ and $\widetilde{\CE}\subset\CE'$ have no nonzero vector in common.)  Now
given $\Be$, consider the solution to 
\beq \Be,\Bj\in\CU,~~~~\BE\in\CE,~~~~\BJ\in\CJ,~~~~
\Bj+\BJ=\BL(\Be+\BE),~~~~{\rm where}~\BL=\sum_{i=1}^nz_i\BGL_i, \eeq{I5.26}
where $\BGL_i$ is the projection onto $\CP_i$, and from the definition of $\BZ$, $\Bj=\BZ\Be$. Since the second condition in \eq{I5.12e} implies $\CV$ and $\widetilde{\CJ}$ have no vector in common we have
\beq \CW=\CV\oplus\widetilde{\CJ}, \eeq{I5.27}
and consequently any $\Bj\in\CU$ has the decomposition
\beq \Bj=-\BJ'_1+\widetilde{\BJ},\quad\BJ'_1\in\CV,\quad \widetilde{\BJ}\in\widetilde{\CJ}, \eeq{I5.28}
which defines the nonsingular operator $\BM:\CU\to\CV$ such that $\BJ'_1=\BM\Bj$. Defining
\beq \BE'_2=\Be+\BE,\quad \BJ'_2=\Bj+\BJ, \eeq{I5.29}
we have
\beqa \BE'_1+\BE'_2 & = & \BE'_1+\Be+\BE=\widetilde{\BE}+\BE\in\CE', \nonum
\BJ'_1+\BJ'_2 & = & \BJ'_1+\Bj+\BJ=\widetilde{\BJ}+\BJ\in\CJ',
\eeqa{I5.30}
and
\beq \BJ'_1=\BM\Bj=\BM\BZ\Be=\BM\BZ\BK\BE'_1, \eeq{I5.31}
which by definition of the associated $Y$-function implies
\beq \BY(z_1,z_2,\ldots,z_n)=\BM\BZ(z_1,z_2,\ldots,z_n)\BK. \eeq{I5.32}
This is analogous to the relation (20.29) in \citeAPY{Milton:2002:TOC} obtained in the case where the subspaces are mutually orthogonal.

In particular by letting $z_1=z_2=\ldots=z_n=1$ we obtain
\beq \BY(1,1,\ldots,1)=\BM\BK. \eeq{I5.33}
If $\Bv_1$, $\Bv_2$,\ldots $\Bv_m$ are a basis for $\CV$, and we choose
$\BK\Bv_1$, $\BK\Bv_2$,\ldots $\BK\Bv_m$ as our basis for $\CU$ then with these
bases $\BK$ is represented by the identity matrix $\BK=\BI$ and \eq{I5.32} and
\eq{I5.33} imply
\beq  \BY(z_1,z_2,\ldots,z_n)=\BY(1,1,\ldots,1)\BZ(z_1,z_2,\ldots,z_n). \eeq{I5.34}

\section{Reduction operations on subspace collections}
\labsect{Ireduct}
\setcounter{equation}{0}

Extension is one way to go from a $Z(n)$ subspace collection to a $Y(n)$ subspace collection. Another way is through
reduction,%
\index{subspace collections!reduction operations}
which has some features in common with normalization.%
\index{subspace collections!normalization operations}
Given a $Z(n)$ subspace collection
\beq \CH=\CU\oplus\CE\oplus\CJ=\CP_1\oplus\CP_2\oplus\cdots\oplus\CP_n,
\eeq{I11.1}
let $\BGG_0$ be the projection onto $\CU$, and let $\BGL_j$ be the projection onto
$\CP_j$. Define
\beqa \CK=\CE\oplus\CJ, \quad \CP'_j=\CP_j\cap\CK~~~{\rm for}~~j=1,2,\ldots,n, \nonum
      \CV=(\BI-\BGG_0)[\BGL_1\CU\oplus\BGL_2\CU\oplus\cdots\oplus\BGL_n\CU]\subset\CK, \quad
      \widetilde{\CP}_j=\BGL_j\CU.
\eeqa{I11.2}
We now establish that
\beq \CW\equiv\widetilde{\CP}_1\oplus\widetilde{\CP}_2\oplus\cdots\oplus\widetilde{\CP}_n=\CU\oplus\CV. \eeq{I11.3}
First note that $\CV$ and $\CU$ have no vector in common since $\CV\subset\CK$, and similarly the subspaces
$\widetilde{\CP}_j$ have no vector in common since $\widetilde{\CP}_j\subset\CP_j$. Clearly $\CW$ contains
$\CU$ since the projections $\BGL_j$ sum to the identity. To show it contains $\CV$ note that
\beq  \CV\subset \BGL_1\CU\oplus\BGL_2\CU\oplus\cdots\oplus\BGL_n\CU+\BGG_0[\BGL_1\CU\oplus\BGL_2\CU\oplus\cdots\oplus\BGL_n\CU]
\subset \CW+\CU=\CW. \eeq{I11.4}
Therefore we have that $\CU\oplus\CV\subset\CW$. The converse inclusion that $\CW\subset\CU\oplus\CV$ follows from the inclusion
\beq \widetilde{\CP}_j=[\BGG_0+(\BI-\BGG_0)]\BGL_j\CU\subset\CU\oplus\CV,
\eeq{I11.5}
which establishes \eq{I11.3}. Next, to establish that for all $j$,
\beq \CP_j=\widetilde{\CP}_j\oplus\CP'_j, \eeq{I11.6}
we need to assume that for all $j$
\beq \widetilde{\CP}_j\cap\CK=0, \eeq{I11.6a}
and that
\beq \BGL_j\Bu=0,\quad \Bu\in\CU \eeq{I11.6b}
only has the trivial solution $\Bu=0$, i.e.
\beq \CU\cap(\CP_1\oplus\CP_2\oplus\ldots\oplus\CP_{j-1}\oplus\CP_{j+1}\oplus\ldots\oplus\CP_n)=0. \eeq{I11.6c}
These conditions imply that 
\beq \CU=\BGG_0\BGL_j\CU, \eeq{I11.6d}
and hence that 
\beq \CU\subset \BGL_j\CU\oplus(\BI-\BGG_0)\BGL_j\CU, \eeq{I11.6e}
which in turn implies that
\beq \CU\subset\widetilde{\CP}_j+\CV.  \eeq{I11.7}
Then any vector $\BP\in\CP_j$ has the unique decomposition
\beq \BP=\Bu+\BK,~~{\rm with}~~\Bu\in\CU,~~\BK\in\CK, \eeq{I11.8}
and according to \eq{I11.7}, $\Bu$ has the unique decomposition
\beq \Bu=\Bv+\widetilde{\BP}~~{\rm with}~~\Bv\in\CV,~~\widetilde{\BP}\in\widetilde{\CP}_j, \eeq{I11.9}
which is unique because $\CV\subset\CK$ and $\widetilde{\CP}_j$ have no nonzero vector in common.
Therefore $\BP$ has the unique decomposition
\beq \BP=\widetilde{\BP}+\BP', \eeq{I11.10}
where
\beq \BP'=\Bv+\BK=\BP-\widetilde{\BP}\in\CP_j\cap\CK=\CP'_j. \eeq{I11.11}
This decomposition and the fact that \eq{I11.6a} implies $\widetilde{\CP}_j$
and $\CP'_j\subset\CK$ have no vector in common establishes \eq{I11.6}. 

So we deduce that
\beqa \CH=\CU\oplus\CE\oplus\CJ & = & \widetilde{\CP}_1\oplus\widetilde{\CP}_2\oplus\cdots\oplus\widetilde{\CP}_n\oplus\CP'_1\oplus\CP'_2\oplus\cdots\oplus\CP'_n \nonum
&=& \CU\oplus\CV\oplus\CP'_1\oplus\CP'_2\oplus\cdots\oplus\CP'_n,
\eeqa{I11.12}
and since the $\CP'_j$, $j=1,2,\ldots,n$ are all contained in $\CK$ it follows that
\beq \CK=\CE\oplus\CJ=\CV\oplus\CP'_1\oplus\CP'_2\oplus\cdots\oplus\CP'_n. \eeq{I11.13}

Now suppose that given $\Be\in\CU$ we can solve the equations
\beq \Bj+\BJ_1=\BL(\Be+\BE_1),\quad \BJ_1=-\BY\BE_1,\quad \Be,\Bj\in\CU,\quad \BE_1,\BJ_1\in\CV, \eeq{I11.14}
where $\BY$ is the $\BY$-operator associated with the subspace collection \eq{I11.13}. From the $Y$-problem we have
\beq \BE=\BE_1+\BE_2\in\CE,\quad \BJ=\BJ_1+\BJ_2\in\CJ\quad\BJ_2=\BL\BE_2,\quad \BE_2,\BJ_2\in\CH', \eeq{I11.15}
where
\beq \CH'=\CP'_1\oplus\CP'_2\oplus\cdots\oplus\CP'_n. \eeq{I11.16}
Since
\beq \Bj+\BJ_1+\BJ_2=\BL(\Be+\BE_1+\BE_2), \eeq{I11.17}
we see that these fields solve the $Z$-problem
\beq \Be,\Bj\in\CU,\quad \BE\in\CE,\quad \BJ\in\CJ,\quad\Bj+\BJ=\BL(\Be+\BE), \eeq{I11.18}
and by definition $\Bj=\BZ\Be$. To solve \eq{I11.14} let $\BGP_1$ be the projection onto $\CV$. Then \eq{I11.14} implies
\beq -\BY\BE_1=\BGP_1\BL(\Be+\BGP_1\BE_1), \eeq{I11.19}
giving
\beq \BE_1=-\BGP_1(\BY+\BGP_1\BL\BGP_1)^{-1}\BGP_1\BL\Be, \eeq{I11.20}
where the inverse is to be taken on the subspace $\CV$. It follows that
\beq \Bj+\BJ_1=\BL\Be-\BL\BGP_1(\BY+\BGP_1\BL\BGP_1)^{-1}\BGP_1\BL\Be, \eeq{I11.21}
implying
\beq \BZ=\BGG_0\BL\BGG_0-\BGG_0\BL\BGP_1(\BY+\BGP_1\BL\BGP_1)^{-1}\BGP_1\BL\BGG_0. \eeq{I11.22}
This formula%
\index{Zo@$\BZ$-operator formula}
is analogous to that given in (29.12) of \citeAPY{Milton:2002:TOC}.

To obtain a more explicit way of writing \eq{I11.22} let us suppose we are given a basis
$\Bu_1,\Bu_2,\ldots,\Bu_m$ of $\CU$. Since \eq{I11.6b} only has the trivial solution $\Bu=0$ each space $\widetilde{\CP}_j$ has dimension $m$. It then follows from \eq{I11.3} that $\CV$ has dimension $m(n-1)$. Also, for $i=1,2,\ldots,n-1$, 
\eq{I11.3} implies $\BGL_i\Bu_j$ has the unique decomposition
\beq \BGL_i\Bu_j=\sum_{k}w_{ijk}\Bu_k+\Bv_{ij},\quad \Bv_{ij}\in\CV, \eeq{I11.23}
for some set of constants $w_{ijk}$. To show that the vectors $\Bv_{ij}$, which number $m(n-1)$, are independent, let us 
suppose
\beq 0=\sum_{i=1}^{n-1}\sum_{j=1}^m c_{ij}\Bv_{ij}=\sum_{i=1}^{n-1}\sum_{j=1}^mc_{ij}
(\BGL_i\Bu_j-\sum_{k=1}^mw_{ijk}\Bu_k). \eeq{I11.24}
By letting $\BGL_n$ act on this equation and taking into account that \eq{I11.6b} only has the trivial solution $\Bu=0$
we see that
\beq \sum_{i=1}^{n-1}\sum_{j=1}^m\sum_{k=1}^mc_{ij}w_{ijk}\Bu_k=0. \eeq{I11.25}
Then substituting this in \eq{I11.24} and letting $\BGL_i$, $i\ne n$, act on \eq{I11.24} and 
again taking into account that \eq{I11.6b} only has the trivial solution $\Bu=0$ we obtain
\beq \sum_{j=1}^mc_{ij}\Bu_j=0, \eeq{I11.26}
which shows that all the $c_{ij}$ must be zero. Therefore let us take the vectors $\Bv_{ij}$ as our basis for $\CV$. 

The identities 
\beq \BGP_1\BGL_i\BGG_0\Bu_j=\Bv_{ij}, \quad \BGG_0\BGL_i\BGG_0\Bu_j=\sum_{k}w_{ijk}\Bu_k, \eeq{I11.27}
which follow from \eq{I11.23} then gives the matrix representations for $\BGP_1\BGL_i\BGG_0$ and
$\BGG_0\BGL_i\BGG_0$ in these bases, when $i\ne m$.
Using the fact that $\BGL_n=\BI-\sum_{i\ne n}\BGL_i$ we obtain
\beq \BGG_0\BL\BGG_0=z_n\BGG_0+\sum_{i=1}^{n-1}(z_i-z_n)\BGG_0\BGL_i\BGG_0, \quad
\BGP_1\BL\BGG_0=\sum_{i=1}^{n-1}(z_i-z_n)\BGP_1\BGL_i\BGG_0. \eeq{I11.28}
Now for $p\ne n$ (and $i\ne n$) \eq{I11.23} implies (no sum over $p$)
\beqa \BGL_p\Bv_{ij} & = & \sum_k(\Gd_{pi}\Gd_{kj}-w_{ijk})\BGL_p\Bu_k\nonum
& = & \sum_k(\Gd_{pi}\Gd_{kj}-w_{ijk})(\Bv_{pk}+\sum_{q}w_{pkq}\Bu_q).
\eeqa{I11.29}
Thus we deduce
\beqa \BGG_0\BGL_p\BGP_1\Bv_{ij} & = & \sum_k(\Gd_{pi}\Gd_{kj}-w_{ijk})\sum_{q}w_{pkq}\Bu_q, \nonum
\BGP_1\BGL_p\BGP_1\Bv_{ij} & = & \sum_k(\Gd_{pi}\Gd_{kj}-w_{ijk})\Bv_{pk},
\eeqa{I11.30}
which gives the matrix representation for the operators $\BGG_0\BGL_p\BGP_1$ and $\BGP_1\BGL_p\BGP_1$ in these
bases ($p\ne n$), in terms of which we obtain the representation for the operators
\beq \BGG_0\BL\BGP_1=\sum_{p=1}^{n-1}(z_p-z_n)\BGG_0\BGL_p\BGP_1, \quad
\BGP_1\BL\BGP_1=z_n\BGP_1+\sum_{p=1}^{n-1}(z_p-z_n)\BGP_1\BGL_p\BGP_1.
\eeq{I11.31}
Thus all the matrices representing the operators entering \eq{I11.22}, aside from $\BY$, only depend on the parameters
$w_{ijk}$ and these parameters can be obtained from the representation in the basis $\Bu_1,\Bu_2,\ldots,\Bu_m$ of 
$\BZ$ when the differences $z_i-z_n$, $i=1,2,\ldots, n-1$ are small. To first order in these differences,
\eq{I11.22},\eq{I11.28}, and \eq{I11.31} imply
\beq \BZ\Bu_j\approx z_n\Bu_j+\sum_{i=1}^{n-1}(z_i-z_n)\sum_{k}w_{ijk}\Bu_k. \eeq{I11.32}
Thus knowing this expansion one can recover all the parameters $w_{ijk}$. 
\section{``Continued fraction expansions'' of subspace collections.}
\index{continued fraction expansions}
The idea to developing the continued fraction is that by a succession of reduction%
\index{subspace collections!reduction operations}
and normalization operations%
\index{subspace collections!normalization operations}
one obtains a series of recursion relations
\beqa \BZ  & = & \BGG_0\BL\BGG_0-\BGG_0\BL\BGP_1(\BY+\BGP_1\BL\BGP_1)^{-1}\BGP_1\BL\BGG_0, \label{I11.33} \\
 \BY  & = &  \BM^{(1)}\BZ^{(1)}\BK^{(1)}, \label{I11.34} \\
 \BZ^{(1)} & = &  \BGG_0^{(1)}\BL^{(1)}\BGG_0^{(1)}-\BGG_0^{(1)}\BL^{(1)}\BGP_1^{(1)}(\BY^{(1)}+\BGP_1^{(1)}\BL^{(1)}\BGP_1^{(1)})^{-1}\BGP_1^{(1)}\BL^{(1)}\BGG_0^{(1)}, \label{I11.35}\\
 \BY^{(1)}  & = &  \BM^{(2)}\BZ^{(2)}\BK^{(2)}, \label{I11.36} \\
 \BZ^{(2)} & = & \BGG_0^{(2)}\BL^{(2)}\BGG_0^{(2)}-\BGG_0^{(2)}\BL^{(2)}\BGP_1^{(2)}(\BY^{(2)}+\BGP_1^{(2)}\BL^{(2)}\BGP_1^{(2)})^{-1}\BGP_1^{(2)}\BL^{(2)}\BGG_0^{(2)}, \label{I11.37}
\end{eqnarray}
and so forth, until the dimension of the remaining space goes to zero, or until one (or more) of the assumptions necessary to proceed with the normalization or reduction operation
does not hold. By substituting \eq{I11.34} in \eq{I11.33}, then substituting \eq{I11.35} in the resulting expression, and subsequently substituting \eq{I11.36} in this expression, and so on,
one develops the continued fraction expansion for $\BZ$ incorporating the variables $z_1$, $z_2$, $\ldots$, $z_n$ and, 
as one goes down the continued fraction,
information contained in the series expansion \eq{I4a.1a} at successively higher and higher
levels of truncation. We do not address in this book whether one can go ahead with the continued fraction expansion (and if so
how) when the assumptions made to proceed with the normalization or reduction operation do not hold. In the process of developing the continued fraction through 
reduction and normalization operations, one could at those steps where one is dealing with a $Y$-problem make any desired reference transformation%
\index{reference transformation}
as described in Section 12. In this way one incorporates information at the subspace collection level that corresponds at the function level to known values of the function, and derivatives, at various points.

Such continued fraction expansions form the basis of the field equation recursion method%
\index{field equation recursion method}
for bounding the effective moduli of composites
(\citeAY{Milton:1985:TCC}; Milton \citeyearNP{Milton:1987:MCEa}, \citeyearNP{Milton:1987:MCEb}, \citeyearNP{Milton:1991:FER}; \citeAY{Clark:1994:MEC}; \citeAY{Clark:1997:CFR} and Chapter 29 of \cite{Milton:2002:TOC} in the abstract theory of composites as described in Chapter 2 of this book (\citeAY{Milton:2016:ETC}): see also 
Section 9.10 and Chapter 10 of \citeAPY{Milton:2016:ETC}). The basic idea, at least when we have an orthogonal subspace collection, 
is that crude estimates or bounds on the operator $\BZ^{(j)}$ or $\BY^{(j)}$ at some intermediate level $j$ give through the above recursion
relations good approximations or tight bounds on $\BZ$ or $\BY$ incorporating the parameters that enter the recursion relations at the different levels
up to level $j$ (obtained from series expansions up to a given order of the solutions of the $Z$-problem or $Y$-problem).

\section*{Acknowledgments}
G.W. Milton thanks his husband John K. Patton for suggesting the name superfunction. 

\bibliographystyle{mod-xchicago}
\bibliography{/home/milton/tcbook,/home/milton/newref}

\end{document}